\documentclass[11pt, english]{article} 

\usepackage[utf8]{inputenc} 

\usepackage[T1]{fontenc}

\usepackage{amsmath,amsthm,amssymb,bm,mathrsfs,MnSymbol,empheq} 

\usepackage{cite}
\usepackage{hyperref}
\hypersetup{colorlinks, linkcolor={blue}, citecolor={blue}, urlcolor={blue}}
\usepackage[nameinlink]{cleveref}

\usepackage[text={16cm, 23cm}, centering, a4paper]{geometry}
\allowdisplaybreaks 

\usepackage{subcaption}

\usepackage{aliascnt}

\newtheorem{theorem}{Theorem}[section]

\newaliascnt{lemma}{theorem}
\newtheorem{lemma}[lemma]{Lemma}
\aliascntresetthe{lemma}

\newaliascnt{proposition}{theorem}
\newtheorem{proposition}[proposition]{Proposition}
\aliascntresetthe{proposition}

\newaliascnt{definition}{theorem}
\newtheorem{definition}[definition]{Definition}
\aliascntresetthe{definition}

\newaliascnt{remark}{theorem}
\newtheorem{remark}[remark]{Remark}
\aliascntresetthe{remark}

\newaliascnt{example}{theorem}
\newtheorem{example}[example]{Example}
\aliascntresetthe{example}

\newaliascnt{assumption}{theorem}
\newtheorem{assumption}[assumption]{Assumption}
\aliascntresetthe{assumption}

\newaliascnt{rn}{theorem}
\newtheorem{rn}[rn]{Remark and Notation}
\aliascntresetthe{rn}

\crefname{theorem}{theorem}{theorems}
\Crefname{theorem}{Theorem}{Theorems}

\crefname{lemma}{lemma}{lemmas}
\Crefname{lemma}{Lemma}{Lemmas}

\crefname{proposition}{proposition}{propositions}
\Crefname{proposition}{Proposition}{Propositions}

\crefname{definition}{definition}{definitions}
\Crefname{definition}{Definition}{Definitions}

\crefname{remark}{remark}{remarks}
\Crefname{remark}{Remark}{Remarks}

\crefname{example}{example}{examples}
\Crefname{example}{Example}{Examples}

\crefname{assumption}{assumption}{assumptions}
\Crefname{assumption}{Assumption}{Assumptions}

\crefname{rn}{remark and notation}{remarks and notation}
\Crefname{rn}{Remark and Notation}{Remarks and Notation}

\numberwithin{equation}{section} 

\newcommand{\tabs}[1]{\vert #1 \vert}
\newcommand{\tnorm}[1]{\Vert #1 \Vert}
\newcommand{\abs}[1]{\left\lvert #1 \right\rvert}
\newcommand{\norm}[1]{\left\lVert #1 \right\rVert}
\newcommand{\ddd}{\textnormal{d}}
\newcommand{\R}{\mathbb{R}}
\newcommand{\N}{\mathbb{N}}
\newcommand{\Indi}{\boldsymbol{1}}
\newcommand{\dotcup}{\mathbin{\dot{\cup}}}

\DeclareMathOperator{\diag}{diag}
\DeclareMathOperator*{\esssup}{ess\,sup}
\DeclareMathOperator*{\argmax}{arg\,max}
\DeclareMathOperator*{\argmin}{arg\,min}

\title{\vspace{-1.6cm} Common Noise by Random Measures: Constructing Mean-Field Equilibria for Competitive Investment and Hedging \footnote{Stefanie Hesse thanks Deutsche Forschungsgemeinschaft (DFG, German Research Foundation, Project-ID 410208580, IRTG 2544 "Stochastic Analysis in Interaction") for funding, the audience of the 16th Bachelier Colloquium in M{\'e}tabief (01/2024) and of the AMaMeF Conference in Verona (06/2025) for feedback, and Ludovic Tangpi for inspiring discussions. Dirk Becherer acknowledges funding by DFG - CRC/TRR 388 "Rough Analysis, Stochastic Dynamics and Related Fields" - Project ID 516748464.\\
The first version (arXiv:2408.01175v1) of this preprint has been published under the slightly different title ``Common Noise by Random Measures: Mean-Field Equilibria for Competitive Investment and Hedging''; subsequent versions (with the new title) feature several new results relative to ver1.
}}
\author{Dirk Becherer \thanks{Humboldt University of Berlin, Germany (dirk.becherer+at+hu-berlin.de).}
\and Stefanie Hesse \thanks{Humboldt University of Berlin, Germany (stefanie.hesse+at+hu-berlin.de).}}
\date{November 3, 2026}

\begin{document} 

\maketitle

\vspace{-1.1cm}

\textbf{Abstract.}
We construct Nash equilibria in mean-field portfolio games of optimal investment and hedging under relative performance concerns with exponential (CARA) utility preferences. Common noise dynamics are modeled by integer-valued random measures, for instance Poisson random measures, in addition to Brownian motions. Agents differ in individual risk aversions, competition weights, and initial capital endowments, while their contingent claim liabilities depend on both common and idiosyncratic risk factors. Mean-field equilibria are characterized by solutions to McKean-Vlasov backward stochastic differential equations with jumps, for which we prove existence and uniqueness of solutions, without assuming mean-field interaction to be small. Moreover, we show how the equilibrium can be constructed from the optimal strategy of a single-agent optimization problem (without mean-field interaction) via an appropriate projection. Using successive changes of measure, our derivation provides a decomposition of the equilibrium strategy into three components with clear interpretations. Finally, we show how a limiting mean-field game of quadratic (instead of utility-based) hedging with relative performance concerns arises for vanishing risk aversion. 

     
    
\section{Introduction}

We characterize Nash equilibria for mean-field games (MFGs) of investment and hedging in incomplete markets, where each agent aims to maximize her relative utility in comparison with the mean-field (e.g., industry-wide) average within a large population. Idiosyncratic noise is given by integer-valued random measures, while common noise is described by both integer-valued random measures and Brownian motions. This brings together the two topical themes of mean-field games with jumps and the problem of competitive portfolio optimization for combined hedging and investment with relative exponential utility preferences. 

The contributions of the present paper are as follows. 
(i) We provide a formulation of common and idiosyncratic noise in terms of integer-valued random measures (see \Cref{subsec:KeyAss}), which generalize Poisson random measures. 
(ii) We establish a one-to-one correspondence between mean-field equilibria and solutions to a McKean–Vlasov backward stochastic differential equation with jumps (see \Cref{proposition:CharUnderAssFBSDEhasSolu}) and prove well-posedness of the latter (see \Cref{lemma:SoluMKVbsde}) without imposing a weak interaction condition.
(iii) We derive a decomposition of mean-field equilibrium (MFE) strategies into three components (see \hyperref[remark:InterpretationDecompMFE]{\Cref*{remark:InterpretationDecompMFE}.1}), each admitting a clear financial interpretation.
(iv) We construct the MFE strategy from the optimal strategy of a single-agent problem (without mean-field effects) combined with a projection step (see \Cref{charMFE} in \Cref{theorem:mainTheorem}).
(v) In a setting with Markov switching dynamics, we show how the MFG is characterized by the solution to a system of partial differential equations (see equations \eqref{ThetaPDEoptimalStratSingleAgentProb}-\eqref{PDEmfeGeneral}). 
(vi) In the limit of vanishing risk aversion, we obtain a new limiting MFG of quadratic (rather than utility-based) hedging under relative performance concerns (see \Cref{theorem:LimitMFG}). 

We investigate a new type of mean-field portfolio games that has both idiosyncratic and common noise components in terms of jumps. 
We address the problem of competitive investment \textit{and hedging} of liabilities against jump-driven unpredictable risks in an incomplete market. We offer an approach where both common and idiosyncratic components of jump noise are described by a general integer-valued random measure. This permits, for instance, time- and $\omega$-dependent compensators and both common and idiosyncratic times and marks of marked jump processes. Filtrations supporting such random measures, together with Brownian noises, are non-continuous in that they admit non-continuous martingales. Our application of competitive investment and hedging in incomplete markets with non-continuous filtrations is motivated, for instance, by risk management problems at the interface of finance and insurance (see \Cref{remark:ExampleMoeller}) with respect to utility-based preferences (cf.\ \cite{moller03reinsurance,becherer03,fritelliEtal11} and references therein). Because of market incompleteness, only partial hedging of liabilities is possible in general. This is a reason to consider exponential (instead of power) utility preferences for the competitive hedging and investment problem, which are finitely defined on the real line with constant absolute risk aversion (CARA) and are known for their convenient properties and for being linked to a suitable dual problem by relative entropy minimization over equivalent martingale measures (see \cite{becherer03,fritelliEtal11} and remarks after \Cref{proposition:BechUnique}). 

Our main result, \Cref{theorem:mainTheorem}, provides existence and uniqueness as well as a construction of the MFE in terms of the optimal strategy of some single-agent optimization problem \eqref{SingleAgentProbForB} without any mean-field interaction. In particular, computing the MFE essentially reduces to solving this single-agent optimization problem \eqref{SingleAgentProbForB}, whose unique solution is also given in our main result \Cref{theorem:mainTheorem}. 
Based on this characterization we give in \Cref{sec:PDE} in a setting with a Markov switching process a further characterization of the MFE by solutions of a system of reaction-diffusion partial differential equations (PDEs) to illustrate by computational examples the MFE, its decomposition, and the effect of common jump noise on the MFE. The explicit construction of the MFE in our main result, \Cref{theorem:mainTheorem}, is also crucial in \Cref{sec:AlphaTo0} to show how a new MFG \eqref{MFGAlpha0} of quadratic hedging with relative performance concerns emerges as risk aversion vanishes (i.e., as risk tolerance tends to infinity), without imposing a Markovian restriction. 

From a methodological viewpoint, our proofs avoid the need to assume a weak interaction condition, which would require mean-field interaction (i.e., competition weights) to be sufficiently small. The weak interaction condition is known from the wider literature \cite{horst2005stationary,fu2021mean} and also typically appears in papers on the relative utility maximization problem, e.g., \cite{tangpi2024optimal}, notably when general measurable coefficients for the (possibly non-Markovian) price dynamics are admitted. 
We prove our main \Cref{theorem:mainTheorem} about existence and uniqueness of the MFE, together with its characterization, via a chain of one-to-one correspondences. First, we relate mean-field equilibria of the MFG \eqref{MFG} to mean-field equilibria of an auxiliary MFG \eqref{MFG3} formulated under a new measure $\widehat{\mathbb{P}}$ (\Cref{proposition:MFGgleichMFG3}). We then show in \Cref{proposition:CharUnderAssFBSDEhasSolu} that mean-field equilibria of this auxiliary MFG \eqref{MFG3} are in one-to-one correspondence with solutions of the McKean-Vlasov backward stochastic differential equation with jumps (JBSDE) defined in \eqref{CharFBSDE}. Well-posedness (\Cref{lemma:SoluMKVbsde}) follows by reducing the McKean-Vlasov JBSDE \eqref{CharFBSDE} to an auxiliary JBSDE \eqref{auxiliaryJBSDE} with bounded terminal condition and without mean-field interaction. Using that the auxiliary JBSDE \eqref{auxiliaryJBSDE} has a unique trivial solution, we obtain even an explicit representation of the unique solution to the McKean-Vlasov JBSDE \eqref{CharFBSDE}, in terms of the optimal strategy of some single-agent problem \eqref{SingleAgentProbForB}, and thereby an explicit representation of the MFE in \Cref{charMFE} in \Cref{theorem:mainTheorem}. Moreover, our approach to prove existence and uniqueness of the MFE also yields a decomposition of the MFE strategy with a clear interpretation, namely into a Merton ratio, a hedging component, and a mean-field component (see \hyperref[remark:InterpretationDecompMFE]{\Cref*{remark:InterpretationDecompMFE}.1}).

MFGs were introduced by Lasry and Lions \cite{lasry2007mean} and Huang et al. \cite{huang2006large}. Until recently, research on MFGs in continuous time has been predominantly focused on probabilistic bases described solely by Brownian motions \cite{huang2006large,lasry2007mean,carmonadelarue18}. Recently, there is increasing interest in MFGs with jumps \cite{benazzoli2020mean,cecchin2020probabilistic,dianetti2026ergodic,bo2024mean,bauerle2023nash,hernandez2024mean,alasseur2023mfg,bayraktar2021finite}, for instance Poisson-like jumps, as a common model for the occurrence of non-predictable events. Yet, most of the literature has been concerned with the technically simpler case where jumps are entirely part of the idiosyncratic noise \cite{benazzoli2020mean,cecchin2020probabilistic,dianetti2026ergodic,bo2024mean}. MFGs where the jumps have systemic global influence, i.e., are part of the common noise, are studied, e.g., in \cite{bauerle2023nash,hernandez2024mean} and \cite{bayraktar2021finite}. Notably, \cite{alasseur2023mfg} studies a MFG about production from exhaustible resources that exhibits jumps even in both idiosyncratic and common noises. 

The problem of relative utility maximization and its investigation by backward stochastic differential equation (BSDE) methods can be traced back to Espinosa and Touzi \cite{espinosa2010stochastic,espinosa2015optimal}. Meanwhile, this problem has been studied in complete \cite{espinosa2010stochastic,espinosa2015optimal} and incomplete markets \cite{espinosa2010stochastic,espinosa2015optimal,lacker2019mean,lacker2020many,fu2023mean1}, for Markovian \cite{espinosa2010stochastic,espinosa2015optimal,lacker2019mean,dos2022forward,dianetti2024optimal,bo2024mean}, and non-Markovian \cite{fu2023mean1,fu2023mean2,tangpi2024optimal,bauerle2023nash} asset price models, and for different utility functions. With few exceptions, like \cite{fujii2023mean}, most articles on relative utility maximization do consider pure investment problems without additional hedging of terminal liabilities. Except for \cite{bo2024mean,bauerle2023nash} and \cite{bank2025much}, the previously mentioned work on relative performance concerns considers games on Brownian filtrations. In \cite{bo2024mean} jumps occur only as part of the idiosyncratic noise, and in \cite{bauerle2023nash} they are only part of the common noise. In the recent paper \cite{bank2025much}, dynamic Markovian programming results are obtained in a model where marked jump processes can reveal common and idiosyncratic information to investors, see \Cref{example:CondIndepAndWPRP} and subsequent comments. 

The present paper is organized as follows. \Cref{sec:Preliminaries} introduces the general setting and recalls basic facts about stochastic integration with respect to random jump measures. \Cref{sec:FormulationMFG} formulates the MFG of hedging and investment under relative performance concerns and the main \Cref{theorem:mainTheorem} for the characterization of mean-field equilibria. The remainder of \Cref{sec:FormulationMFG} as well as \Cref{sec:transformation} and \Cref{sec:CharacterizationMFE} is devoted to proving this theorem, which also requires a description of optimal strategies for the embedded single-agent optimization problem, being provided in \Cref{sec:FormulationMFG}. \Cref{sec:transformation} establishes a one-to-one relation to an auxiliary MFG whose characterization is obtained in \Cref{sec:CharacterizationMFE}, where finally we combine the results to prove the main theorem. Building on that, \Cref{sec:PDE} derives a PDE characterization and presents numerical illustrations of the MFE in a Markovian example,  while \Cref{sec:AlphaTo0} demonstrates how a new MFG of quadratic hedging with relative performance concerns emerges in the limit for vanishing risk aversion.

\section{Preliminaries}\label{sec:Preliminaries} 

This section provides notations and the probabilistic setup. \Cref{subsec:IntegrationRandomMeausres} introduces assumptions on the stochastic basis and recalls essential facts on stochastic integration with respect to random measures for jumps. We refer to \cite{jacod2003limit,cohen2015stochastic} for more details of the theory. \Cref{subsec:KeyAss} presents our abstract general setting for common and idiosyncratic noises originating from integer-valued random measures and describes our two key assumptions concerning the filtrations involved, along with several concrete examples. \Cref{subsec:FinacialMarketFramework} formulates the financial market model for the mean-field game (MFG) of investment and hedging (in \Cref{sec:FormulationMFG}).

\subsection{Stochastic basis and integration w.r.t.\ random measures}\label{subsec:IntegrationRandomMeausres}

We work on a stochastic basis $(\Omega, \mathcal{F}, \mathbb{F}, \mathbb{P})$  with a finite time horizon $T<\infty$ and a filtration $\mathbb{F}=(\mathcal{F}_t)_{t\in [0,T]}$ satisfying the usual conditions of right-continuity and completeness. Thus we can and do take all semimartingales to have c{\`a}dl{\`a}g paths. Let $(E,\mathcal{B}(E))$ be a measurable space where $\mathcal{B}(E)$ denotes the Borel $\sigma$-field on $E$. For simplicity, let $E:=\R^{\ell}\backslash\{0\},\ \ell \in \N$ (more generally, one may take a Blackwell space) and define $\tilde{\Omega}:=\Omega\times [0,T]\times E$. Let the stochastic basis support a $d$-dimensional Brownian motion $W=(W_t)_{t\in [0,T]}$ and an integer-valued random measure
\begin{align*}
	\boldsymbol{\mu}(\ddd t, \ddd e)=(\boldsymbol{\mu}(\omega, \ddd t, \ddd e)\vert \omega \in \Omega)
\end{align*}
on $([0,T]\times E, \mathcal{B}([0,T])\otimes \mathcal{B}(E))$ with compensator $\boldsymbol{\nu}$ ({w.r.t.\ $\mathbb{P}$ and $\mathbb{F}$}), cf.\ ~\cite{he1992semimartingale,jacod2003limit}. We call $\tilde{\boldsymbol{\mu}}=\boldsymbol{\mu}-\boldsymbol{\nu}$ the compensated measure of $\boldsymbol{\mu}$ under $\mathbb{P}$ (and $\mathbb{F}$). For sub-filtration $\mathbb{G}\subseteq \mathbb{F}$, let $\mathcal{P}(\mathbb{G})$ (resp.\ $\mathcal{O}(\mathbb{G})$) denote the predictable (resp.\ optional) $\sigma$-field on $\Omega\times[0,T]$ w.r.t.\ $\mathbb{G}$. We call a function on $\Omega\times [0,T]$ that is $\mathcal{P}(\mathbb{G})$-measurable $\mathbb{G}$-predictable. By $\tilde{\mathcal{P}}(\mathbb{G}):=\mathcal{P}(\mathbb{G})\otimes\mathcal{B}(E)$ (resp. $\tilde{\mathcal{O}}(\mathbb{G}):=\mathcal{O}(\mathbb{G})\otimes\mathcal{B}(E)$) we denote the predictable (resp. optional) $\sigma$-field on $\tilde{\Omega}$ w.r.t.\ $\mathbb{G}$. We assume the following.
 
\begin{assumption}\label{assumption:NuFinite}
 	The compensator is absolutely continuous with respect to the product measure $\lambda\otimes \ddd t$ with Radon-Nikodym density $\zeta$, such that
	\begin{equation*}
		\boldsymbol{\nu}(\omega,\ddd t, \ddd e)=\zeta(\omega, t, e)\lambda(\ddd e)\ddd t
	\end{equation*}
	holds, with $\lambda$ being a finite measure on $(E, \mathcal{B}(E))$ and  density $\zeta$ being $\tilde{\mathcal{P}}(\mathbb{F})$-measurable and bounded, such that 
	\begin{equation}
		0\le \zeta(\omega,t,e)\le c_{\boldsymbol{\nu}}<\infty, \qquad \mathbb{P}\otimes\lambda\otimes\ddd t\text{-a.e.}\label{bddZeta}
	\end{equation}
	for some constant $c_{\boldsymbol{\nu}}$. Thus, in particular, $\boldsymbol{\nu}([0,T]\times E)\le c_{\boldsymbol{\nu}} T\lambda(E)<\infty\ \text{almost surely}$.
\end{assumption}
 
\begin{example}\label{example:IntValRM}
	Let $N$ be a Poisson process with intensity $\lambda^N \in (0, \infty)$ and let $D^i$, $i\in \N$, be independent, integrable, $E$-valued random variables, identically distributed according to $\lambda^D$ on $(E,\mathcal{B}(E))$. The integer-valued random measure associated with the compound Poisson process $C=\sum_{i=1}^N D^i$ is then given by $\mu^C(\ddd t, \ddd e):=\sum_{s: \Delta C_s\neq 0}\delta_{(s,\Delta C_s)}(\ddd t, \ddd e)$ with $\Delta C_t:=C_t-C_{t-}$ denoting jumps, and the associated compensator $\nu^C(\ddd t, \ddd e)=\lambda^D(\ddd e)\lambda^N\ddd t$ satisfies \Cref{assumption:NuFinite}.
\end{example}
 
\begin{remark}
	The integer-valued random measures setup permits jump processes that are significantly more general than marked or compound Poisson processes. They allow dependence on time and $\omega$ for jump intensities and jump heights. They can accommodate, for instance, (semi-) Markov chains (appearing in regime-switching models, see \cite{dianetti2026ergodic}); finite-state jump dynamics (see \cite{cecchin2020probabilistic}); or even more general step-processes (see \cite[Ch.XI]{he1992semimartingale} and \cite[Example~2.1]{becherer2019monotone}). Notably, the Brownian motion $W$ and the integer-valued random measure $\boldsymbol{\mu}$ can be stochastically dependent, meaning that jump heights and intensities could depend on the history of Brownian trajectories. See, e.g., \cite[eq.(1)]{andreis2018mckean} for state-dependent jump intensities via thinning of a Poisson random measure. Instead of repeating examples already given in \cite{becherer2019monotone,becherer2006bounded}, we present below in \Cref{example:CondIndepAndWPRP} several variants of other examples of increasing generality, which are centered around and extend the basic example with independent compound Poisson processes being the common and idiosyncratic noise components originating from jumps.
\end{remark}

Let $U:\tilde{\Omega}\rightarrow \R$ be an $\tilde{\mathcal{O}}(\mathbb{F})$-measurable function. The integral process of $U$ with respect to the integer-valued random measure $\boldsymbol{\mu}$ is defined by
\begin{align*}
	U\ast \boldsymbol{\mu}_t(\omega)= 
	\begin{cases}
		\int_{[0,t]\times E}U(\omega,s,e)\boldsymbol{\mu}(\omega,\ddd s, \ddd e) & \text{if } \int_{[0,t]\times E}\abs{U(\omega,s,e)}\boldsymbol{\mu}(\omega,\ddd s, \ddd e)<\infty,\\
		+\infty & \text{otherwise}.
	\end{cases}
\end{align*}
The integral process for the compensator $\boldsymbol{\nu}$ is defined analogously (cf.\ \cite[Eq.II.1.5]{jacod2003limit}). We recall that for any $\tilde{\mathcal{P}}(\mathbb{F})$-measurable function $U$, by the definition of the compensator, $\mathbb{E}\left[\tabs{U}\ast \boldsymbol{\mu}\right]=\mathbb{E}\left[\tabs{U}\ast \boldsymbol{\nu}\right]$ holds (cf.\ \cite[Thm.II.1.8.(i)]{jacod2003limit}). If, moreover, $(\tabs{U}^2\ast \boldsymbol{\mu})^{1/2}$ is locally integrable, then $U$ is integrable with respect to $\tilde{\boldsymbol{\mu}}$ and the process $U\ast \tilde{\boldsymbol{\mu}}=(U\ast \tilde{\boldsymbol{\mu}}_t)_{t\in [0, T]}$ is defined as the purely discontinuous local martingale such that the jump process of $U \ast \tilde{\boldsymbol{\mu}}$ is equal to $(\int_E U_t(e)\boldsymbol{\mu}(\{t\},\ddd e))_{t\in [0,T]}$ (cf.\ \cite[Def.II.1.27]{jacod2003limit}). Furthermore, the equality $U \ast \tilde{\boldsymbol{\mu}}=U \ast \boldsymbol{\mu}-U \ast \boldsymbol{\nu}$ applies (cf.\ \cite[Prop.II.1.28]{jacod2003limit}).

An integer-valued random measure $\boldsymbol{\mu}$ is called optional w.r.t.\ $\mathbb{F}$ if for each positive $\tilde{\mathcal{O}}(\mathbb{F})$-measurable function $U$ the process $U\ast \boldsymbol{\mu}$ is $\mathbb{F}$-optional (cf.\ \cite[Def.13.2.9]{cohen2015stochastic}). The natural filtration $(\mathcal{F}_t^{\boldsymbol{\mu}})_{t\in [0,T]}$ of $\boldsymbol{\mu}$ is defined as the smallest filtration such that $\boldsymbol{\mu}$ is optional (see \cite[Sect.~13.6.1]{cohen2015stochastic}). 

Next, we define spaces of processes, common in the literature, for $\mathbb{Q}$ denoting some probability measure on $(\Omega,\mathcal{F})$: For $p\in [1,\infty]$, let $\mathbb{S}^{p}(\mathbb{Q})$ denote the space of $\R$-valued $\mathbb{F}$-adapted c{\`a}dl{\`a}g semimartingales $(Y_t)_{t\in [0,T]}$ with $\tnorm{Y}_{\mathbb{S}^p(\mathbb{Q})}:=\tnorm{\sup_{t\in [0,T]}\tabs{Y_t}}_{L^p(\mathbb{Q})}<\infty$. Let $\mathcal{L}_T^{2}(\mathbb{Q})$ denote the space of $\mathbb{F}$-predictable processes $Z$ taking values in $\R^d$ with $\tnorm{Z}_{\mathcal{L}_T^2(\mathbb{Q})}^2:=\mathbb{E}^\mathbb{Q}[\int_0^T \abs{Z_t}^2 \ddd t]<\infty$. For brevity, we write $\mathbb{E}$ for the expectation $\mathbb{E}^{\mathbb{P}}$. Let $\boldsymbol{\nu}^{\mathbb{Q}}$ be the compensator of $\boldsymbol{\mu}$ under the measure $\mathbb{Q}$. We denote by $\mathcal{L}_{\boldsymbol{\nu}^{\mathbb{Q}}}^{2}(\mathbb{Q})$ the space of $\tilde{\mathcal{P}}(\mathbb{F})$-measurable functions $U:\tilde{\Omega}\rightarrow \R$ with $\tnorm{U}_{\mathcal{L}_{\boldsymbol{\nu}^{\mathbb{Q}}}^2(\mathbb{Q})}^2:=\mathbb{E}^{\mathbb{Q}}[\smallint_0^T\smallint_E \abs{U_t(e)}^2 \boldsymbol{\nu}^{\mathbb{Q}}(\ddd t, \ddd e)]<\infty$. Let $BMO(\mathbb{Q})$ denote the space of $BMO(\mathbb{Q})$-martingales (see \cite[Def.10.6]{he1992semimartingale}). Let $\mathbb{H}_{BMO}^2(\mathbb{Q})$ denote the space of $\mathbb{F}$-predictable processes $Z$ with bounded norm $\tnorm{Z}_{\mathbb{H}_{\text{BMO}}^2(\mathbb{Q})}^2:=\sup_{t\in [0,T]}\tnorm{\mathbb{E}^\mathbb{Q}[\smallint_t^T \tabs{Z_s}^2 \ddd s \vert \mathcal{F}_t]}_{L^\infty}<\infty$.

For $U^{\mathbb{Q}}\in \mathcal{L}_{\boldsymbol{\nu}^{\mathbb{Q}}}^2(\mathbb{Q})$, $U^{\mathbb{Q}}\ast \tilde{\boldsymbol{\mu}}^{\mathbb{Q}}=U^{\mathbb{Q}}\ast (\boldsymbol{\mu}-\boldsymbol{\nu}^{\mathbb{Q}})$ is a square integrable $\mathbb{Q}$-martingale; we write $U^{\mathbb{Q}}\ast \tilde{\boldsymbol{\mu}}_t^{\mathbb{Q}}= \int_0^t \int_E U_s^{\mathbb{Q}}(e)\tilde{\boldsymbol{\mu}}^{\mathbb{Q}}(\ddd s, \ddd e)$ (cf.\ \cite[Thm.II.1.33.a)]{jacod2003limit}). For $Z^{\mathbb{Q}}\in \mathbb{H}_{\text{BMO}}^2(\mathbb{Q})$ and $W^\mathbb{Q}$ being a $\mathbb{Q}$-Brownian motion, $\int Z^{\mathbb{Q}} \ddd W^{\mathbb{Q}}$ is in $BMO(\mathbb{Q})$, see \cite[Thm.10.9.4]{he1992semimartingale}.

\subsection{Basic Assumptions on Common Noise and Filtration}\label{subsec:KeyAss}

This subsection introduces our two key assumptions, \Cref{assumption:weakPPR} and \Cref{assumption:condIndep}, about relevant filtrations and common noise, that are assumed for the analysis in the sequel. These assumptions are fairly general but abstract, and are to be explained and illustrated by concrete examples in \Cref{example:CondIndepAndWPRP}.

The first key assumption concerns martingale representation with respect to the overall filtration $\mathbb{F}$, jointly by the Brownian motion and the compensated integer-valued random measure. It is a natural assumption which 
enables applicability of solution theory for BSDEs with jumps in the sequel.

\begin{assumption}\label{assumption:weakPPR}
	$W$ and $\tilde{\boldsymbol{\mu}}:=\boldsymbol{\mu}-\boldsymbol{\nu}$ have the weak property of predictable representation w.r.t.\ the filtration $\mathbb{F}$. This means that every square integrable $\mathbb{F}$-martingale $M$ has a representation
	\begin{align*}
		M_t=M_0 + \int_0^t Z_s \ddd W_s + U \ast \tilde{\boldsymbol{\mu}}_t, \quad t\in [0,T],
	\end{align*}
    where $Z:\Omega \times [0,T] \rightarrow \R^d$ and $U: \tilde{\Omega}\rightarrow \R$ are predictable processes such that $\mathbb{E}[\int_0^T \tabs{Z_t}^2 \ddd t]<\infty$ and $\mathbb{E}[\tabs{U}^2 \ast \boldsymbol{\nu}_T]<\infty$. In particular, this means that both stochastic integrals lie in the space of the square-integrable martingales. 
\end{assumption}

Regarding notations, let $\mathbb{F}^W=(\mathcal{F}_t^W)_{t\in [0,T]}$ denote the natural filtration of the Brownian motion $W$. We denote the common noise filtration by $\mathbb{F}^0=(\mathcal{F}_t^0)_{t\in [0,T]}$. 

We assume throughout that the common noise includes the Brownian filtration, i.e., $\mathbb{F}^W\subseteq \mathbb{F}^0\subseteq \mathbb{F}$, and that $\mathbb{F}^0$ satisfies the second key \Cref{assumption:condIndep} below. In the interest of generality, we do not further specify the common noise filtration beyond the abstract assumptions. Yet, we exemplify below how those are satisfied in several more specific situations.

\begin{remark}
	We take Brownian noise as being common noise entirely (which is a simplification), since the original contributions of the present paper concern the originating of common and idiosyncratic noises for the MFG from integer-valued random measures, and we aim for generality related to the latter only.

    In contrast, for finite-state MFGs with Wright Fisher common noise (on the simplex) by \cite{bayraktar2021finite}, the collective Brownian motions form the common noise \cite[see eqns.(2.4),(2.7) and subsequent comments]{bayraktar2021finite}, while jumps occur solely idiosyncratically as part of the conditional Markov chain dynamics of the representative agent’s state \cite[cf.\ eqn.~(2.10)]{bayraktar2021finite}.
    
    While our general framework deliberately avoids concrete specifications of the common noise filtration to preserve broad applicability, it is instructive to recall the standard setting for common noise in MFGs on Brownian filtrations, as in  \cite{carmonadelarue18}. There, common and idiosyncratic noises originate from independent Brownian motions, and  martingale representation is provided by a sum of two strongly orthogonal stochastic integrals against those. 
    An analogous separation for an integer-valued random measure $\boldsymbol{\mu}$ can be obtained by decomposing its mark space into disjoint sets $E=E_0\dotcup E_1$, $E_0,E_1\in \mathcal{B}(E)$, thereby defining a common part $\mu^0$ and an idiosyncratic part $\mu^1$ by 
	\begin{align}
		\mu^0(A)= \boldsymbol{\mu}(A\cap ([0,T]\times E_0))\quad  \text{and}\quad \mu^1(A)=\boldsymbol{\mu}(A\cap ([0,T]\times E_1))\label{Mu01}
	\end{align}
	for any $A\in \mathcal{B}([0,T])\otimes \mathcal{B}(E)$. 
    Based on such a decomposition, one may take the common noise filtration $\mathbb{F}^{0}=(\mathcal{F}_t^0)_{t\in [0,T]}$, with $\mathcal{F}^0_t= \mathcal{F}_t^{\mu^0}\vee \mathcal{F}_t^{W}$, to be generated by the natural filtrations from $\mu^0$ and from the Brownian motion $W$. Provided that Assumptions \ref{assumption:weakPPR} and \ref{assumption:condIndep} are satisfied (see \Cref{example:CondIndepAndWPRP}), our later MFG analysis then applies. Yet, those assumptions also admit common noise examples beyond such a decomposition of $\boldsymbol{\mu}$ (see \Cref{example:CondIndepAndWPRP}, part~2). 
\end{remark}

\begin{assumption}\label{assumption:condIndep}
	For all $t\in [0,T]$ the $\sigma$-fields $\mathcal{F}_t$ and $\mathcal{F}_T^{0}$ are conditionally independent given $\mathcal{F}_t^{0}$. This means that the idiosyncratic information up to time $t$ provides no information for the future common-noise information, but the common-noise information up to time $t$ can provide information on the future idiosyncratic information.
\end{assumption}

The \Cref{assumption:condIndep} is known in the literature as immersion property or (H)-hypothesis. It is equivalent to requiring that any martingale in the common noise filtration $(\mathcal{F}_t^{0})$ is also a martingale in the overall filtration $(\mathcal{F}_t)$, we refer to \cite[Vol.II, Def.1.2 \& Prop.1.3]{carmonadelarue18}, also for further equivalences. \Cref{example:CondIndepAndWPRP} provides concrete examples of increasing generality, ranging from independent common and idiosyncratic jump information at distinct times (1) to simultaneous jump information (2) and time and $\omega$ dependent intensities and jump sizes (3), satisfying our key Assumptions \ref{assumption:weakPPR} and \ref{assumption:condIndep}.

\begin{example}\label{example:CondIndepAndWPRP} 
	1. Let $C^0$, $C^1$ be $\R^{\ell_0}$ and $\R^{\ell_1}$-valued compound Poisson processes with $\ell_0,\ell_1\in\N$, $W$ a Brownian motion and $\mathcal{A}$ a $\sigma$-field. We denote the corresponding natural filtrations by $(\mathcal{F}_t^{C^0}), (\mathcal{F}_t^{C^1})$ and $(\mathcal{F}_t^W)$. Let $\mathcal{F}_T^{C^0}, \mathcal{F}_T^{C^1}, \mathcal{F}_T^W$ and $\mathcal{A}$ be independent. Setting $E_0:=(\R^{\ell_0}\backslash\{0\})\times \{0\}$, $E_1:=\R^{\ell_0+\ell_1}\backslash (\R^{\ell_0}\times \{0\})$, $C:=(C^0, C^1)$ and $\boldsymbol{\mu}^C(\ddd t, \ddd e): =\sum_{s:\Delta C_s\neq 0}\delta_{(s,\Delta C_s)}(\ddd t, \ddd e)$, it follows for the natural filtration $(\mathcal{F}_t^{\boldsymbol{\mu}^C})$ of the integer-valued random measure $\boldsymbol{\mu}^C$, that $(\mathcal{F}_t^{\boldsymbol{\mu}^C})=\sigma(\mathcal{F}_t^{C^0}, \mathcal{F}_t^{C^1})$, and for the integer-valued random measures $\mu^i$, $i=0,1$, defined by \eqref{Mu01}, that the completion of the natural filtration $(\mathcal{F}_t^{\mu^i})$ of $\mu^i$ and the completion of the filtration $(\mathcal{F}_t^{C^i})$ are identical (cf.\ \cite[Sect.13.6.1]{cohen2015stochastic}). Let the basic filtration $(\mathcal{F}_t)$ be the usual filtration generated by $(\mathcal{A}, \mathcal{F}_t^W, \mathcal{F}_t^{\boldsymbol{\mu}^C})$ and the common noise filtration the filtration generated by $(\mathcal{F}_t^W,\mathcal{F}_t^{\mu^0})$, then \Cref{assumption:condIndep} is satisfied. This can be proven by verifying an equivalent condition \cite[Sect.3.2 Prop.13 (ii)]{rao2006probability} for conditional independence using the just mentioned relations of the filtrations and a suitable intersection-stable generator. Furthermore, the Brownian motion $W$ and the compensated integer-valued random measure $\tilde{\boldsymbol{\mu}}^C$ satisfy \Cref{assumption:weakPPR} according to \cite[Example 2.1.3]{becherer2019monotone}, which can be argued using general theory for so-called step-processes (cf.\ \cite[Ch.XI]{he1992semimartingale}).

	2. To extend and generalize the example in part~1, let $\bar{C}$ be a further compound Poisson process of dimension $\bar{\ell}_0+\bar{\ell}_1$, independent of the $\sigma$-fields of part~1, where $\bar{\ell}_0,\bar{\ell}_1\in\N$. Let its jump heights $\bar{D}^{0,k}\in \R^{\bar{\ell}_0}\backslash\{0\}$ and $\bar{D}^{1,\bar{k}}\in \R^{\bar{\ell}_1}\backslash \{0\}$, $k, \bar{k}\in \N$ be independent and for each fixed $i=0,1$ let $\bar{D}^{i,k}$ be identically distributed. We denote by $0<T_1\le T_2\le \dotsc \le T$ the sequence of jump times of the process $\bar{C}$. Let the basic filtration $(\mathcal{F}_t)$ be the usual filtration generated by $(\mathcal{A}, W_s, C_s^0, C_s^1, \bar{D}^{0,k}\Indi_{\{T_k\le s\}}, \bar{D}^{1,k}\Indi_{\{T_k\le s\}}\vert k\in \N,s \le t)$ and the common noise filtration the filtration generated by $(W_s, C_s^0, \bar{D}^{0,k}\Indi_{\{T_k\le s\}} \vert k\in \N,s\le t)$. Then it follows with arguments as in part~1 that \Cref{assumption:condIndep} is satisfied. We note that the $\sigma$-field generated by $(C_s^0, C_s^1, \bar{D}^{0,k}\Indi_{\{T_k\le s\}}, \bar{D}^{1,k}\Indi_{\{T_k\le s\}}\vert k\in \N,s \le t)$ can also be obtained as a $\sigma$-field defined by an integer-valued random measure $\boldsymbol{\mu}$ (see \cite[Sections~13.3,13.6]{cohen2015stochastic}). Again, the Brownian motion $W$ and the compensated integer-valued random measure $\tilde{\boldsymbol{\mu}}$ then satisfy \Cref{assumption:weakPPR} (see \cite[Example 2.1.3]{becherer2019monotone}).

	3. Let $W$ and $\tilde{\boldsymbol{\mu}}$ be such that \Cref{assumption:weakPPR} is satisfied under $\mathbb{P}$ and, in addition, let the basic filtration $\mathbb{F}$ and the common noise filtration $\mathbb{F}^0$ satisfy \Cref{assumption:condIndep}. Let $\mathbb{Q}$ be a probability measure absolutely continuous with respect to $\mathbb{P}$ with density process $Z$ adapted to the common noise filtration $\mathbb{F}^0$. Then the Brownian motion $W^{\mathbb{Q}}=W-\int (Z_{-})^{-1}\ddd \langle Z,W \rangle$ under the new measure $\mathbb{Q}$ and the $\mathbb{Q}$-compensated jump measure $\tilde{\boldsymbol{\mu}}^{\mathbb{Q}}:=\mu-\nu^{\mathbb{Q}}$ satisfy \Cref{assumption:weakPPR} (cf.\ \cite[Thm.13.22]{he1992semimartingale}) and the filtrations satisfy \Cref{assumption:condIndep} as before. The latter can be proven directly using \cite[Prop.13]{rao2006probability}. Based on the previous examples, dependencies can now also be created between the processes, whereby the main assumptions \Cref{assumption:weakPPR} and \Cref{assumption:condIndep} are still satisfied.
\end{example}

Parts 2 and 3 of the preceding example illustrate how marked jumps convey both common and idiosyncratic information, through different respective components of the jump mark, in the mean-field portfolio game. Such aspects are systematically investigated in the context of Meyer-sigma-fields in \cite{bank2025much} with a view toward jump signals, using dynamic programming methods in a Markovian setup to derive the best response map and Schauder's fixed point theorem to show that equilibria exist.

\subsection{The financial market framework}\label{subsec:FinacialMarketFramework}

The market contains a riskless numeraire asset (with unit price one) and $d$ risky assets, whose (discounted) prices evolve as an It\^{o}-process, described by the SDE
\begin{equation}
    \ddd S_t=\diag(S_t^i)_{i\in \{1,\dotsc, d\}} \sigma_t(\varphi_t\ddd t + \ddd W_t), \quad t\in [0,T], \label{MertonSDE}
\end{equation}
with $S_0\in (0,\infty)^d$, where $\diag(x)$ denotes the diagonal matrix with entries $x$ on the diagonal. The market price of risk $\varphi$ is an $\mathbb{F}^0$-predictable, $\R^d$-valued and bounded process. The volatility $\sigma$ is an $\R^{d\times d}$-valued, $\mathbb{F}^0$-predictable process such that $\sigma_t$ is invertible ($\mathbb{P}\otimes\ddd t$-a.e.) and integrable with respect to
\begin{equation}
	\widehat{W}:=W+\int \varphi_t \ddd t.\label{DefWHat}
\end{equation}
An investment strategy $\vartheta$ is taken to be an $\mathbb{F}$-predictable, $S$-integrable, $\R^d$-valued process. A strategy $\vartheta$ describes the dynamic holding of risky assets $S$ over time. The discounted gains process associated with the strategy $\vartheta$ is given by 
\begin{equation}
	\left(\int_0^t \vartheta_s \ddd S_s\right)_{t\in [0,T]}.\label{DPP}
\end{equation}
We define $\Sigma_t:=(\diag(S_t^i)_{i\in \{1,\dotsc, d\}})\sigma_t$, {write $\Sigma^\top$ for the process of transposed matrices}, and will use the parametrization $\theta=\Sigma^{\top}\vartheta$ to simplify the exposition in the following, keeping in mind that by 
\begin{equation}
	\theta(\vartheta):=\Sigma^{\top}\vartheta\quad \text{and}\quad \vartheta(\theta)=(\Sigma^{\top})^{-1}\theta \label{parametrizationVarthetaTheta}
\end{equation}
we have a bijection between the parameterizations of strategies $\theta$ and $\vartheta$. The discounted gains process \eqref{DPP} can thereby simply be written as 
\begin{equation*}
	\left(\int_0^t \theta_s \ddd \widehat{W}_s\right)_{t\in [0,T]}.
\end{equation*}

The paper's analysis benefits from the fact that stochastic integrals with respect to $\widehat{W}$ constitute gains attainable by trading while those with respect to $\widetilde{\mu}$ have vanishing quadratic covariation with such gains.

\begin{remark}
    While the price processes of tradeable assets available for investment and hedging in our model are continuous, the contingent claim liabilities and also the coefficients in the SDE for the It\^o-process $S$ could depend on the evolution of the integer-valued random measure $\boldsymbol{\mu}$ and of the Brownian motion $W$, in a general (measurable, possibly path-dependent) way. Consequently, the financial market is incomplete, the overall filtration being non-Brownian. We emphasize that this incompleteness may already arise from the dependence of the market coefficients on the jump measure, so that even European put and call options are generally not replicable.
\end{remark}

\section{The mean-field game of investment and hedging}\label{sec:FormulationMFG}

This section introduces the mean-field game (MFG) of hedging and investment with relative performance concerns under exponential utility preferences and states the main \Cref{theorem:mainTheorem}, which fully characterizes the mean-field (Nash) equilibrium. We also state and prove \Cref{proposition:BechUnique}, which solves a class of single-agent optimization problems needed to prove \Cref{theorem:mainTheorem}. 

We impose the following assumption henceforth, in addition to those of \Cref{sec:Preliminaries}.

\begin{assumption}\label{ass:AinitialSigmaField}
    Let $\mathcal{A}\subseteq\mathcal{F}_0$ be a $\sigma$-field independent of $\mathcal{F}_T^0$ and let the individual (representative) agent's characteristics $x_0$ (initial capital endowment), $\alpha$ (risk aversion) and $\rho$ (competition weight) be $\mathcal{A}$-measurable random variables. Furthermore, the initial capital $x_0\in L^2(\mathcal{A}, \mathbb{P})$ is square integrable, the risk aversion $\alpha\in L^{\infty}(\mathcal{A})$ is strictly positive, bounded and bounded away from $0$ and the competition weight $\rho\in L^\infty(\mathcal{A})$ is bounded with $\mathbb{E}[\rho]\neq 1$ ($\rho\ge 0$ represents a competitive interaction and $\rho<0$ a homophilic one). Finally, the contingent claim is a bounded $\mathcal{F}_T$-measurable random variable $B\in L^{\infty}(\mathcal{F}_T)$.
\end{assumption}

We consider an investor who aims to maximize her relative utility with respect to the mean-field (say, e.g., industry) average by finding an optimal investment \textit{and} (partial) hedging strategy given her liabilities $B$ in the financial market ($S$, respectively $\widehat{W}$), in competition with other agents who of course trade in the same market. Throughout, she optimizes over strategies in $\mathbb{H}_{\text{BMO}}^2(\mathbb{P})$. An equilibrium to our MFG of investment and hedging can be described along the following three-step scheme,
\begin{equation}
	\begin{cases}
		\text{1. Given a real-valued random variable } F,\\
		\text{2. }\text{find } \widetilde{\theta}\in \underset{\theta \in \mathbb{H}_{\text{BMO}}^2(\mathbb{P})}{\argmax} \mathbb{E}\left[-\exp(-\alpha(X_T^\theta-B-\rho F))\vert \mathcal{F}_0\right],\\
		\text{\phantom{2. }} \text{for wealth process $ X^{\theta}$ given by }\ddd X_t^{\theta}= \theta_t(\varphi_t \ddd t + \ddd W_t), X_0^\theta=x_0.\\
		\text{3. Find a fixed point such that } F=\mathbb{E}\left[X_T^{\widetilde{\theta}}-B\vert \mathcal{F}_T^{0}\right],\\
		\text{\phantom{3. }where $X_T^{\widetilde{\theta}}$ is the optimal wealth from step 2 with respect to $F$.}
	\end{cases}\label{MFG0}
\end{equation}
When we write $\widetilde{\theta}\in \argmax_{\theta \in \mathbb{H}_{\text{BMO}}^2(\mathbb{P})} \mathbb{E}\left[-\exp(-\alpha(X_T^\theta-B-\rho F))\vert \mathcal{F}_0\right]$, we mean that $\widetilde{\theta}\in \mathbb{H}_{\text{BMO}}^2(\mathbb{P})$ and that for all $\theta\in \mathbb{H}_{\text{BMO}}^2(\mathbb{P})$ we have
\begin{equation*}
	\mathbb{E}\left[-\exp(-\alpha(X_T^\theta-B-\rho F))\vert \mathcal{F}_0\right]\le \mathbb{E}\left[-\exp(-\alpha(X_T^{\widetilde{\theta}}-B-\rho F))\vert \mathcal{F}_0\right] \ \text{a.s.} 
\end{equation*}

\begin{definition}[mean-field equilibrium]\label{definition:MFE}
    A strategy $\widetilde{\theta}\in \mathbb{H}_{\text{BMO}}^2(\mathbb{P})$ is called a mean-field equilibrium (MFE) for a MFG if, for an  exogenously given random variable $F$, it solves the optimization problem in the second step of the scheme for the MFG and also satisfies the consistency condition in the 3rd step. 
\end{definition}

Equivalently, one could define a MFE as $(\widetilde{\theta},F)$, with $F$ determined by $\widetilde{\theta}$ through the 3rd step of the MFG scheme, in analogy with the measure flow appearing in classical mean-field formulations. For brevity, we use the shorter variant of \Cref{definition:MFE}.

\begin{example}
	A strategy $\widetilde{\theta}\in \mathbb{H}_{\text{BMO}}^2(\mathbb{P})$ is a MFE of the MFG \eqref{MFG0} if and only if it satisfies for the random variable $F:=\mathbb{E}[X_T^{\widetilde{\theta}}-B\vert \mathcal{F}_T^{0}]$ the equality
	\begin{equation*}
		\esssup_{\theta \in \mathbb{H}_{\text{BMO}}^2(\mathbb{P})}\mathbb{E}\left[-\exp\left(-\alpha(X_T^\theta-B-\rho F)\right)\middle\vert \mathcal{F}_0\right]
		=\mathbb{E}\left[-\exp\left(-\alpha(X_T^{\widetilde{\theta}}-B-\rho F)\right)\middle\vert \mathcal{F}_0\right].
	\end{equation*}
\end{example}

The general connection between MFGs and large-population games \cite{huang2006large,carmonadelarue18}, to be elaborated elsewhere for our specific setting, identifies the random variable $F$ as the mean field representing the population’s aggregate terminal performance.

\begin{remark}\label{remark:ExampleMoeller}
	Examples that we have in mind for motivation, are risk management applications at the interface of finance and insurance \cite{moller2002valuation,moller03reinsurance}. The contingent claim $B$ could consist of financial stop-loss contracts covering combined financial and actuarial losses, provided by a reinsurance company to insurers, as described in \cite{moller2002valuation}. Such a claim is of the schematic form
	\begin{equation*}
		B=(\text{InsuranceLoss}_T+\text{FinancialLoss}_T-K_1)^+ \land K_2,
	\end{equation*}
	for retention levels $0\le K_1, K_2<\infty$ (with $\min\{a,b\}=a\land b$); the claim covers losses above level $K_1$ and below level $K_1 + K_2$, where $K_1$ is the deductible, see \cite[Subsec.4.2.3]{moller2002valuation}). Compound Poisson processes (CPP), or generalizations thereof, are common examples of cumulative loss processes of insurance claims, with losses occurring at some intensity rate at times, when the process jumps, while jump heights describe the individual loss sizes. Likewise, It\^{o}-processes (as $S$ in \Cref{subsec:FinacialMarketFramework}) for financial asset price processes encompass standard continuous-time models for hedging and investment from classical Black-Scholes and Merton theory. The paper's assumptions allow jump heights and times to be described by stochastic (predictable) intensities and compensating jump measures, and predictable SDE coefficients in the It\^{o}-process. Such a framework permits stochastic dependencies amongst different compound Poisson processes which can be involved in the idiosyncratic and common noise from jumps for the MFG (see \Cref{example:CondIndepAndWPRP}, parts~1 \& 3), to model reinsurance-specific losses and industry-wide ones, and also between those and the price processes for assets available for optimal partial hedging in the financial market. We emphasize that mentioning (non-)independent CPPs for idiosyncratic and common cumulative noise from jump processes is just a first illustrative example. Our assumptions underlying the analysis encompass generalizations thereof, see \Cref{example:CondIndepAndWPRP}, parts~2 \& 3: For instance, a multivariate CPP (or generalizations thereof) may serve as a model for an abstract risk factor process, and company-specific and industry-wide individual loss sizes could be functions of different coordinate components of jumps in the multivariate CPP jumps (happening at the same times in contrast to part~1).
\end{remark}

For the MFG \eqref{MFG0} of relative utility maximization with optimal investment and hedging in the present work, the term $\mathbb{E}[B\vert \mathcal{F}_T^0]$ does not depend on the choice of strategy. For the sequel, we thus simplify and re-state the scheme \eqref{MFG0} as follows. 
\begin{equation}
    \begin{cases}
        \text{1. Given a real-valued random variable } F,\\
        \text{2. }\text{find } \widetilde{\theta} \text{ in }\underset{\theta \in \mathbb{H}_{\text{BMO}}^2(\mathbb{P})}{\argmax} \mathbb{E}\left[-\exp(-\alpha(X_T^\theta-(B-\rho \mathbb{E}[B\vert \mathcal{F}_T^0])-\rho F))\vert \mathcal{F}_0\right],\\
        \text{\phantom{2. }} \text{for wealth process $ X^{\theta}$ given by } \ddd X_t^{\theta}= \theta_t(\varphi_t \ddd t + \ddd W_t),\; X_0^\theta=x_0.\\
        \text{3. Find a fixed point such that } F= \mathbb{E}\left[X_T^{\widetilde{\theta}}\vert \mathcal{F}_T^{0}\right],\\
        \text{\phantom{3.} for $X_T^{\widetilde{\theta}}$ being the optimal terminal wealth from step 2 with respect to $F$.}
    \end{cases}\label{MFG}
\end{equation}

In our main theorem, \Cref{theorem:mainTheorem}, we characterize the MFE of the MFG \eqref{MFG}. 
We first introduce the equivalent martingale measure $\widehat{\mathbb{P}}^B$ (or later just $\widehat{\mathbb{P}}$) used therein, together with the related measure $\widehat{\mathbb{P}}^0$ needed for the vanishing-risk-aversion limit studied in \Cref{sec:AlphaTo0}. Although the model is formulated under the probability measure $\mathbb{P}$, our analysis also relies on the equivalent martingale measures $\widehat{\mathbb{P}}^B$ and $\widehat{\mathbb{P}}^0$.

The equivalent martingale measures are intimately related to the single-agent utility optimization problem
\begin{equation}
	\text{maximize } \mathbb{E}\left[-\exp(-\alpha(X_T^\theta-\xi))\vert \mathcal{F}_0\right] \text{ over } \theta \in \mathbb{H}_{\text{BMO}}^2(\mathbb{P})\label{MeasureCharaktOptimProb}
\end{equation} 
with bounded liability $\xi$, wherein we let $\xi=B-\rho \mathbb{E}[B\vert \mathcal{F}_T^0]$ in the case of $\widehat{\mathbb{P}}^B$ and $\xi=0$ in the case of $\widehat{\mathbb{P}}^0$. The optimal wealth process for this utility problem and the solution to the JBSDE which characterizes the optimal control, permit us to represent the Radon-Nikodym density of the equivalent martingale measures $\widehat{\mathbb{P}}^\beta$, $\beta\in\{0,B\}$, (related to a suitable dual problem) both as an ordinary and as a stochastic exponential. We will use the latter to define the measures $\widehat{\mathbb{P}}^\beta$, $\beta\in \{0,B\}$, by the Radon-Nikodym density process
\begin{equation}
	\frac{\ddd \widehat{\mathbb{P}}^\beta}{\ddd \mathbb{P}} \bigg\vert_{\mathcal{F}_t}
	=\mathscr{E}\left(- \int_0^\cdot \varphi_s \ddd W_s+ \int_0^\cdot\int_E \left(\exp(\alpha U_s^\beta (e))-1 \right) \tilde{\boldsymbol{\mu}}(\ddd s, \ddd e)\right)_t,\quad t\in [0,T], \label{changeOfMEasure} 
\end{equation} 
satisfying thus $\mathbb{P}=\widehat{\mathbb{P}}^\beta$ on $\mathcal{F}_0$ in particular, where $(Y^\beta,Z^\beta,U^\beta)\in \mathbb{S}^\infty(\mathbb{P})\times \mathbb{H}_{\text{BMO}}^2(\mathbb{P})\times \mathcal{L}_{\boldsymbol{\nu}}^2(\mathbb{P})$ with $U^\beta$ being bounded is the unique solution (see \Cref{lemma:solutionYB}) to the JBSDE
\begin{align}
	\left\{\arraycolsep=1pt \def\arraystretch{1.2}
	\begin{array}{ll}
		\ddd Y_t^\beta&= (Z_t^\beta \varphi_t + \frac{1}{2\alpha}\vert \varphi_t \vert^2) \ddd t - \int_E \frac{\exp(\alpha U_t^\beta(e))-1-\alpha U_t^\beta(e)}{\alpha}\zeta(t,e)\lambda(\ddd e) \ddd t\\
		&\quad ~ + Z_t^\beta \ddd W_t + \int_E U_t^\beta(e)\tilde{\boldsymbol{\mu}}(\ddd t, \ddd e),\\
		Y_T^\beta&=\beta-\rho \mathbb{E}[\beta\vert \mathcal{F}_T^0].
	\end{array}\label{trueCharBSDE}\right.
\end{align}
The stochastic exponential $\mathscr{E}(M^\beta)$ for $M^\beta:=- \int_0^\cdot \varphi_s \ddd W_s+ \int_0^\cdot\int_E \exp(\alpha U_s^\beta (e))-1 \tilde{\boldsymbol{\mu}}(\ddd s, \ddd e)$ in \eqref{changeOfMEasure}, for $\beta\in\{0,B\}$, is indeed a strictly positive and uniformly integrable $\mathbb{P}$-martingale and thus a density process \eqref{changeOfMEasure} which defines an equivalent measure $\widehat{\mathbb{P}}^\beta\approx \mathbb{P}$. To see this, note that because of the boundedness of $\alpha, U^\beta$ and $\varphi$, $M^\beta$ is a $BMO(\mathbb{P})$-martingale satisfying $\Delta M^\beta\ge -1+\delta$ for some $\delta$ with $0<\delta\le 1$ (using the notation $\Delta M_t ^\beta:=M_t^\beta-M_{t-}^\beta$). Thus, by \cite[remark after Lem.1]{kazamaki1979transformation} $\mathscr{E}(M^\beta)$ is a strictly positive uniformly integrable $\mathbb{P}$-martingale. One could show (yet, we do not use this later) that $\widehat{\mathbb{P}}^\beta$ is the martingale measure minimizing the entropy relative to the measure $\mathbb{P}^\xi$ being defined in terms of the claim $\xi:=\beta-\rho \mathbb{E}[\beta\vert \mathcal{F}_T^0]$: See comments after \Cref{proposition:BechUnique} and \Cref{StochExpHerleitung} in the proof of \Cref{lemma:UmformungSingleAgentProbInMFG}. 

\begin{remark}\label{remark:measureHatP0}
    The measure $\widehat{\mathbb{P}}^0$ occurs in
    \Cref{sec:AlphaTo0} for the description of 
    the limiting MFG for vanishing risk aversion $\alpha$. 
    To this end, it is helpful that this measure does not depend on $\alpha$ (cf.\ measure $\mathcal{Q}^E$ in \cite[beginning of Sec.4.3]{becherer2006bounded}), we may assume $\alpha=1$. 
    In contrast, the measure $\widehat{\mathbb{P}}^B$ plays the central role for the characterization of mean-field equilibria for the MFG \eqref{MFG} in \Cref{theorem:mainTheorem} and for the financial interpretation of the decomposition in \hyperref[remark:InterpretationDecompMFE]{\Cref*{remark:InterpretationDecompMFE}.1}.
\end{remark}

\begin{remark}\label{remark:PopertiesChangeOfMeasure} 
For $\beta\in \{B,0\}$ we denote by 
$(Y^\beta,Z^\beta,U^\beta)\in \mathbb{S}^\infty(\mathbb{P})\times \mathbb{H}_{\text{BMO}}^2(\mathbb{P})\times \mathcal{L}_{\boldsymbol{\nu}}^2(\mathbb{P})$ with $U^\beta$ the corresponding solutions of the JBSDE \eqref{trueCharBSDE}. The following statements apply.

	1. $\widehat{W}$ defined in \eqref{DefWHat} is a Brownian motion under the measure $\widehat{\mathbb{P}}^{\beta}$. 

	2. The compensator $\widehat{\boldsymbol{\nu}}^\beta$ of $\boldsymbol{\mu}$ under $\widehat{\mathbb{P}}^{\beta}$ is given by $\ddd \widehat{\boldsymbol{\nu}}^\beta=\exp(\alpha U^{\beta}) \ddd \boldsymbol{\nu}$ (see \cite[Thm.III.3.17]{jacod2003limit}). Letting $\widehat{\zeta}^\beta: = \exp(\alpha U^{\beta})\zeta$, we can write $\widehat{\boldsymbol{\nu}}^\beta(\omega, \ddd t, \ddd e)=\widehat{\zeta}^\beta(\omega, t,e)\lambda(\ddd e)\ddd t$ and the $\widehat{\mathbb{P}}^{\beta}$-compensator $\widehat{\boldsymbol{\nu}}^\beta$ satisfies \Cref{assumption:NuFinite}, since $\alpha$ and $U^{\beta}$ are bounded.

	3. The definition of $\mathbb{H}_{\text{BMO}}^2$ a priori depends on the probability measure. Yet, as $M^\beta:=- \int_0^\cdot \varphi_s \ddd W_s+ \int_0^\cdot\int_E \exp(\alpha U_s^{\beta} (e))-1 \tilde{\boldsymbol{\mu}}(\ddd s, \ddd e)$ is a $BMO(\mathbb{P})$-martingale, with $\Delta M^\beta\ge -1+\delta$ for some $\delta\in (0,1]$, a combination of \cite[Thm.10.9]{he1992semimartingale} and \cite[Thm.1]{kazamaki1979transformation} (or \cite[Rem.3.3]{kazamaki1994continuous}) yields that there are constants $c,C>0$ such that
    \begin{align*}
        c\norm{\cdot }_{\mathbb{H}_{\text{BMO}}^2(\mathbb{P})}
        \le \norm{\cdot }_{\mathbb{H}_{\text{BMO}}^2(\widehat{\mathbb{P}}^{\beta})}
        \le C \norm{\cdot}_{\mathbb{H}_{\text{BMO}}^2(\mathbb{P})}.
    \end{align*}
    Hence, $\mathbb{H}_{\text{BMO}}^2(\mathbb{P})=\mathbb{H}_{\text{BMO}}^2(\widehat{\mathbb{P}}^{\beta})$, and thus the strategy space over which we optimize is unchanged under the changes of measure, which is going to be used frequently in the sequel.
\end{remark}

Recalling that $\mathbb{E}=\mathbb{E}^{\mathbb{P}}$, we define moreover for brevity the following notations 
\begin{equation}
    \widehat{\mathbb{P}}:=\widehat{\mathbb{P}}^B, \quad
    \widehat{\mathbb{E}}:=\mathbb{E}^{\widehat{\mathbb{P}}}, \quad
    \widehat{\boldsymbol{\nu}}:=\widehat{\boldsymbol{\nu}}^B, \quad 
    \widehat{\tilde{\boldsymbol{\mu}}}:=\boldsymbol{\mu}-\widehat{\boldsymbol{\nu}}\quad \text{ and }\quad 
    \widehat{\zeta}:=\widehat{\zeta}^B. \label{DefHutNotation}
\end{equation}
To help the reader keep track of the notation, we provide a brief overview of the main conventions.  Objects carrying a hat are associated with a martingale measure. Without a superscript, they refer to the $B$-dependent quantities fixed in \eqref{DefHutNotation}. By $\theta$, we denote trading strategies in the reparameterized form of \Cref{parametrizationVarthetaTheta}. Under this reparameterization, $\int \vartheta\ddd S$ equals $\int \theta\ddd\widehat{W}$, so that trading gains with respect to $S$ are represented by integration with respect to $\widehat{W}$. Plain $\theta$ denotes a generic strategy, superscripted $\theta$ the optimal strategy of a corresponding single-agent problem,  $\widetilde{\theta}$ a MFE strategy, and $\bar{\theta}$ a fixed strategy in an auxiliary problem. Further superscripts are added where needed to distinguish the underlying problems. In particular, superscripts $B$ and $0$ refer to the problems with terminal liabilities $B-\rho \mathbb{E}[B \vert \mathcal{F}_T^0]$ and $0$. The same conventions apply to the non-reparameterized strategies $\vartheta$.

Our main \Cref{theorem:mainTheorem} below establishes existence and uniqueness of the MFE of the MFG \eqref{MFG} and provides two complementary characterizations. First, it gives a JBSDE characterization of the MFE  in terms of a decomposition whose three components admit a financial interpretation (see \hyperref[remark:InterpretationDecompMFE]{\Cref*{remark:InterpretationDecompMFE}.1}). Second, it shows that the MFE can be constructed directly from the optimal strategy of a single-agent problem through an additional projection step. Before stating the theorem, we therefore introduce the $\Pi$-projection used in this construction.

\begin{rn}\label{rn:defPi}
	Since $(\Omega\times [0,T], \mathcal{P}(\mathbb{F}), \mathbb{P}\otimes \ddd t)$ is a finite measure space, for any $z\in L^1(\Omega\times [0,T], \mathcal{P}(\mathbb{F}), \mathbb{P}\otimes \ddd t)$ the 
    conditional expectation with respect to $\mathcal{P}(\mathbb{F}^0)$ is well-defined under the measure $\mathbb{P}\otimes \ddd t$. We denote this $\mathbb{F}^0$-predictable conditional expectation process by $\Pi(z)$ (notation indicating projection). It can be interpreted as the conditional expectation of $z$ given the common information at the process level and directly provides an $\mathbb{F}^0$-predictable process suitable for stochastic integration. 
\end{rn}

Throughout the paper, uniqueness of optimal strategies in single-agent problems and of mean-field equilibria in mean-field games is  understood up to indistinguishability of the corresponding wealth processes. Next, we state our main theorem.

\begin{theorem} \label{theorem:mainTheorem}
There exists a unique mean-field equilibrium $\widetilde{\theta}$ to the mean-field game \eqref{MFG}. It is given by 
 	\begin{equation}
		\widetilde{\theta}
        =\frac{1}{\alpha}\varphi+Z^B+Z
        =\theta^B+Z
        =\theta^B+\frac{\rho}{1-\mathbb{E}[\rho]}\Pi(\theta^B),\label{charMFE} 
	\end{equation} 
    where $\theta^B$ is the unique optimal strategy for the single-agent optimization problem
	\begin{equation}
		 \text{maximize }\mathbb{E}\left[-\exp(-\alpha (X_T^{\theta}-(B-\rho \mathbb{E}[B\vert \mathcal{F}_T^0])))\vert \mathcal{F}_0\right] \text{ over } \theta \in \mathbb{H}_{\text{BMO}}^2(\mathbb{P}) \label{SingleAgentProbForB}
	\end{equation}	
    and given by $\theta^B=\frac{1}{\alpha}\varphi+Z^B$, where $(Y^B,Z^B,U^B)\in \mathbb{S}^\infty(\mathbb{P})\times \mathbb{H}_{\text{BMO}}^2(\mathbb{P})\times \mathcal{L}_{\boldsymbol{\nu}}^2(\mathbb{P})$ with $U^B$ taken bounded is the solution to the JBSDE \eqref{trueCharBSDE} for $\beta:=B$. Moreover, $(Y,Z,U)\in \mathbb{S}^2(\widehat{\mathbb{P}})\times \mathbb{H}_{\text{BMO}}^2(\widehat{\mathbb{P}})\times\mathcal{L}_{\widehat{\boldsymbol{\nu}}}^2(\widehat{\mathbb{P}})$ with $U$ bounded is the unique solution to the McKean-Vlasov JBSDE 
    \begin{align}
		\left\{\arraycolsep=1pt \def\arraystretch{1.2}
		\begin{array}{ll}
    		\ddd Y_t&= - \int_E \frac{\exp(\alpha U_t(e))-1-\alpha U_t(e)}{\alpha}\widehat{\zeta}(t,e)\lambda(\ddd e) \ddd t\\
    		&\quad + Z_t \ddd \widehat{W}_t + \int_E U_t(e)\widehat{\tilde{\boldsymbol{\mu}}}(\ddd t, \ddd e),\\
			Y_T&=\rho \mathbb{E}[x_0+\int_0^T (Z_s +\theta_s^B)\ddd \widehat{W}_s \vert \mathcal{F}_T^{0}].
        \end{array}\label{CharFBSDE}\right.
	\end{align} 
\end{theorem}

\begin{remark}\label{remark:InterpretationDecompMFE} 
    1. The first decomposition of the MFE $\widetilde{\theta}$ in \Cref{charMFE}  permits a financial interpretation as follows. If SDE coefficients in \eqref{MertonSDE} are deterministic, the first component is simply the optimal pure investment strategy (Merton ratio) for maximizing expected exponential (CARA) utility of terminal wealth without any claim and without mean-field interaction. The second component represents the hedging part for the claim $B-\rho \mathbb{E}[B\vert \mathcal{F}_T^0]$ and, when SDE coefficients are stochastic, it additionally also hedges the (intertemporal) stochastic investment opportunities. Finally, \Cref{remark:MKVfbsdeToMKVbsde} explains how the last term relates to the mean-field interaction.\\
    2. The final, third, decomposition of the MFE $\widetilde{\theta}$ in \Cref{charMFE}  is particularly useful to compute mean-field equilibria. It reduces the problem to solving just a single-agent optimization problem \eqref{SingleAgentProbForB} without mean-field interaction or fixed point search, plus an additional projection step. This is useful, e.g., in \Cref{sec:PDE}, to calculate and illustrate MFE via the solution to a PDE system in a Markov-modulated SDE setting; it is moreover crucial for  \Cref{sec:AlphaTo0}, where we show how a new mean-field game emerges in the limit of vanishing risk aversion.
\end{remark}

\begin{remark}
    Within our standing \Cref{ass:AinitialSigmaField}, the condition $\mathbb{E}[\rho]\neq 1$ for our competition weight is more than a merely technical condition ensuring that our existence-and-uniqueness proof for the MFE of the MFG \eqref{MFG} goes through. Indeed, in the case $\rho=1$ one can easily construct from a given MFE additional different equilibria; in particular there cannot be a unique MFE. 
\end{remark} 

The remainder of this section and \Cref{sec:transformation} and \Cref{sec:CharacterizationMFE} are devoted to the proof of \Cref{theorem:mainTheorem}. The proof uses the JBSDE characterization of optimal strategies of single-agent problems from \Cref{proposition:BechUnique} and proceeds through a chain of one-to-one correspondences. First, \Cref{proposition:MFGgleichMFG3} relates MFE of the MFG \eqref{MFG} to those of the auxiliary MFG \eqref{MFG3} via a change of measure to $\widehat{\mathbb{P}}$. 
Second, \Cref{proposition:CharUnderAssFBSDEhasSolu} establishes a one-to-one relation between MFE of the auxiliary MFG \eqref{MFG3} and solutions of the McKean-Vlasov JBSDE \eqref{CharFBSDE}. Finally, \Cref{lemma:SoluMKVbsde} proves well-posedness of the latter by a further one-to-one relation, and provides the explicit representation \eqref{MKVJBSDEReprPiThetaB} of its solution. Combined with the preceding correspondences, this yields the existence, uniqueness, and representations in \eqref{charMFE} of the MFE strategy of the MFG \eqref{MFG} stated in \Cref{theorem:mainTheorem}.

We next characterize optimal strategies for a single-agent utility optimization problem with an exogenously given terminal liability $\xi\in L^2(\mathcal{F}_T)$.

\begin{proposition}[Optimal strategy in single-agent optimization problem]\label{proposition:BechUnique}~\\
    For a given random variable $\xi\in L^2(\mathcal{F}_T,\mathbb{P})$ we have that if the JBSDE 
    \begin{align}
		\left\{\arraycolsep=1pt \def\arraystretch{1.2}
		\begin{array}{ll}
			\ddd Y_t&= (Z_t \varphi_t + \frac{\tabs{\varphi_t}^2}{2\alpha} - \int_E \frac{\exp(\alpha U_t(e))-1-\alpha U_t(e)}{\alpha}\zeta(t,e)\lambda(\ddd e)) \ddd t\\
			&\quad + Z_t \ddd W_t + \int_E U_t(e)\tilde{\boldsymbol{\mu}}(\ddd t, \ddd e),\\
			Y_T&= \xi
		\end{array}\label{BSDEBech}\right.
	\end{align}
    has a solution $(Y,Z,U)\in \mathbb{S}^2(\mathbb{P})\times \mathbb{H}_{BMO}^2(\mathbb{P})\times \mathcal{L}_{\boldsymbol{\nu}}^2(\mathbb{P})$, with $U$ bounded, then there exists a unique strategy $\theta^*$ for the optimization problem 
    \begin{equation}
		\text{maximize } \mathbb{E}\left[-\exp(-\alpha(X_T^\theta-\xi))\vert \mathcal{F}_0\right] \text{ over } \theta \in \mathbb{H}_{\text{BMO}}^2(\mathbb{P}),\label{singleOptimizationProb}
	\end{equation}
	where $X^{\theta}$ is the solution to $\ddd X_t^{\theta}= \theta_t(\varphi \ddd t + \ddd W_t), \; X_0^\theta=x_0\in L^2(\mathcal{F}_0, \mathbb{P})$. It is given by 
    \begin{align}
		\theta^*=Z+\frac{1}{\alpha}\varphi\in \mathbb{H}_{\text{BMO}}^2(\mathbb{P}).\label{ThetaStar}
	\end{align}
    Moreover, the optimal value function $V_t^{\xi, \alpha}(x_t)$ defined by 
	\begin{align*}
    	\esssup_{\theta\in \mathbb{H}_{\text{BMO}}^2}\mathbb{E}\left[-\exp\left(-\alpha \left(x_t+\int_t^T \theta \ddd \widehat{W} - \xi \right)\right)\bigg\vert \mathcal{F}_t \right],
 	\end{align*}
	for starting from initial capital $x_t\in L^2(\mathbb{P},\mathcal{F}_t)$ at time $t$ and having a liability $\xi$, is
	\begin{align*}
		V_t^{\xi, \alpha}(x_t)= -\exp\left(-\alpha (x_t-Y_t)\right),\quad \text{ for } x_t\in L^2(\mathbb{P},\mathcal{F}_t), t\in [0,T].
	\end{align*}
\end{proposition}

By martingale duality theory, the optimal wealth $X^{\theta^*}=x_0+\int_0^\cdot \theta^* d\widehat{W}$ for the primal exponential utility maximization problem is associated to the minimizer of a dual problem to minimize relative entropy with respect to $d\mathbb{P}^{\xi}:= const\ e^{\alpha \xi} d\mathbb{P}$ (where $const$ denotes, here and below, an $\mathcal{F}_0$-measurable normalizing constant chosen such that the resulting measure coincides with $\mathbb{P}$ on $\mathcal{F}_0$) over a suitable set of equivalent martingale measures, which coincides on $\mathcal{F}_0$ with $\mathbb{P}$. The density of this entropy minimizing martingale measure (w.r.t.\ $\mathbb{P}$) is given by $const \exp(-\alpha (X^{\theta^*}_T-\xi))$, see \cite{schweizer2010minimal,fritelliEtal11,becherer03,becherer2006bounded}. These relations to the dual side of the primal utility maximization problem are, however, not used in our analysis. 

\begin{proof}[Proof of \Cref{proposition:BechUnique}]
 	We show at first that strategy \eqref{ThetaStar} is optimal and then prove uniqueness. Optimality is obtained by the familiar martingale optimality principle, just like for continuous Brownian filtrations in \cite{hu2005utility}, in slight adaptation of \cite[Thm.4.1]{becherer2006bounded} to the present technically modified setting with jumps. Because we compare with the terminal wealth of the mean-field average in our MFG, the characterizing JBSDE requires unbounded terminal conditions a priori. Whereas \cite[Thm.4.1]{becherer2006bounded} uses boundedness of $Y$, we are going to argue with BMO-martingales to obtain analogous results.

	Let $(Y,Z,U) \in \mathbb{S}^2(\mathbb{P})\times \mathbb{H}_{BMO}^2(\mathbb{P})\times \mathcal{L}_{\boldsymbol{\nu}}^2(\mathbb{P})$, with $U$ bounded, be a solution to the JBSDE \eqref{BSDEBech}. Since $\alpha$ is bounded and $x_0, Y_0$ are $\mathcal{F}_0$-measurable and $\mathbb{P}$-a.s. finite due to $x_0\in L^2(\mathcal{F}_0), Y\in \mathbb{S}^2(\mathbb{P})$, using the notation of $\widehat{W}$ from \eqref{DefWHat}, it follows that the optimal strategy $\theta^*$ is given by 
	\begin{align}
		\theta^* \in  &\argmax_{\theta \in \mathbb{H}_{\text{BMO}}^2(\mathbb{P})}\mathbb{E}\left[-\exp\left(-\alpha(X_T^\theta-\xi)\right)\Big\vert \mathcal{F}_0\right]\notag\\
        &=\argmax_{\theta \in \mathbb{H}_{\text{BMO}}^2(\mathbb{P})}\mathbb{E}\left[-\exp\bigg(-\alpha\Big(Y_0 + \int_0^T \theta_s \ddd \widehat{W}_s-\xi\Big)\bigg)\bigg\vert \mathcal{F}_0\right]. \label{auxiliaryOptimizationProblem}
	\end{align}
	With the calculations as in the proof of \cite[Thm.4.1]{becherer2006bounded}, one obtains that
	\begin{align}
		&-\exp\left(-\alpha \bigg(Y_0+\int_0^T \theta_s \ddd \widehat{W}_s - \xi\bigg)\right)\notag\\
		&=-e^{(\alpha^2/2) \int_0^T \tabs{\theta_s-Z_s-\varphi_s/\alpha}^2 \ddd s}\mathscr{E}\left(-\alpha \int_0^\cdot \theta_s - Z_s \ddd W_s+\int_0^\cdot \int_E \exp(\alpha U_s(e))-1 \tilde{\boldsymbol{\mu}}(\ddd s, \ddd e)\right)_T\label{StochasticExponentialSingleAgent}
	\end{align}
    for any $\theta \in \mathbb{H}_{BMO}^2(\mathbb{P})$. Since $Z,\theta$ are in  $\mathbb{H}_{BMO}^2(\mathbb{P})$ and $\alpha,U$ are bounded, the martingale inside the stochastic exponential is a $BMO(\mathbb{P})$-martingale with jumps greater than $-1$, and jumps bounded away from $-1$. Thus, according to \cite[remark after Lem.1]{kazamaki1979transformation} the stochastic exponential is a strictly positive uniformly integrable $\mathbb{P}$-martingale. The exponent in the first factor in \eqref{StochasticExponentialSingleAgent} is non-negative. The essential supremum of the expected utility of \eqref{StochasticExponentialSingleAgent} conditioned on $\mathcal{F}_0$  and hence of the auxiliary optimization problem \eqref{auxiliaryOptimizationProblem} is therefore  attained when the non-negative integrand in the exponent of the first factor equals $0$, i.e., for $\theta^*:=Z+\varphi/{\alpha}$. Since $Z$ is in $ \mathbb{H}_{\text{BMO}}^2(\mathbb{P})$, $\varphi$ is bounded and $\alpha$ bounded away from $0$, $\theta^*\in \mathbb{H}_{\text{BMO}}^2(\mathbb{P})$ follows by linearity of $\mathbb{H}_{\text{BMO}}^2(\mathbb{P})$. Since the utility in the auxiliary optimization problem \eqref{auxiliaryOptimizationProblem} for the optimal strategy $\theta^*$ is given by the uniformly integrable stochastic exponential, it follows in particular that the optimization problem \eqref{singleOptimizationProb} is well-posed. This means 
	\begin{align}
		-\infty<\esssup_{\theta\in \mathbb{H}_{\text{BMO}}^2(\mathbb{P})}\mathbb{E}[-\exp(-\alpha(X_T^\theta-\xi))\vert \mathcal{F}_0]<0\label{WellPosedSingleAgent}.
	\end{align}
	Having shown existence, proving uniqueness of the optimal $\theta^*$ now is straightforward. Indeed, by strict concavity of the exponential utility function and convexity of the (linear) space of admissible strategies $\mathbb{H}_{\text{BMO}}^2(\mathbb{P})$ over which the utility maximization problem is posed, one obtains uniqueness (a.s.) of the optimal terminal wealth and thereby of the optimal wealth process, which is a $\widehat{\mathbb{P}}$-martingale thanks to the identity $\mathbb{H}_{\text{BMO}}^2(\mathbb{P})=\mathbb{H}_{\text{BMO}}^2(\widehat{\mathbb{P}})$. The claim for the optimal value function of the single-agent problem follows by the martingale-optimality-principle, just like in \cite{hu2005utility, becherer2006bounded}.
\end{proof}

\begin{remark}\label{remark:becUnterPHut}
	In the sequel, we will also use the analogue of \Cref{proposition:BechUnique} under the measure $\widehat{\mathbb{P}}$ (defined in \eqref{changeOfMEasure} for $\beta=B$, see notation \eqref{DefHutNotation}). Namely, we consider the problem
	\begin{align*}
		\text{maximize }\widehat{\mathbb{E}}[-\exp(-\alpha(X_T^\theta-\xi))\vert \mathcal{F}_0]\text{ over }\theta \in \mathbb{H}_{\text{BMO}}^2(\widehat{\mathbb{P}}).
	\end{align*}
    Recalling the hat notation introduced in \eqref{DefHutNotation}, the corresponding characterization result follows from \Cref{proposition:BechUnique} by replacing $(\mathbb{P},\mathbb{E},W,\tilde{\boldsymbol{\mu}},\boldsymbol{\nu},\zeta,\varphi)$ by $(\widehat{\mathbb{P}},\widehat{\mathbb{E}},\widehat{W},\widehat{\tilde{\boldsymbol{\mu}}},\widehat{\boldsymbol{\nu}},\widehat{\zeta},0)$ throughout its statement and proof. Consequently, the optimal strategy is given by $\theta^*=Z$.
\end{remark}

\begin{remark}\label{remark:MKVfbsdeToMKVbsde}
    To interpret the $Z$-component of the McKean-Vlasov JBSDE \eqref{CharFBSDE} and obtain an associated McKean-Vlasov jump forward-backward SDE (JFBSDE), we introduce a forward process which is to appear in the terminal condition. If $(Y,Z,U)$ in $\mathbb{S}^2(\widehat{\mathbb{P}})\times \mathbb{H}_{\text{BMO}}^2(\widehat{\mathbb{P}})\times\mathcal{L}_{\widehat{\boldsymbol{\nu}}}^2(\widehat{\mathbb{P}})$ with $U$ bounded is a solution of the McKean-Vlasov JBSDE \eqref{CharFBSDE}, then by defining $X:=x_0 + \int (Z +\theta^B) \ddd\widehat{W}$ and using \Cref{remark:PopertiesChangeOfMeasure}, we see that $(X,Y,Z,U)$ in $ \mathbb{S}^2(\widehat{\mathbb{P}})\times \mathbb{S}^2(\widehat{\mathbb{P}})\times \mathbb{H}_{\text{BMO}}^2(\widehat{\mathbb{P}})\times\mathcal{L}_{\widehat{\boldsymbol{\nu}}}^2(\widehat{\mathbb{P}})$, with $U$ bounded, is a solution of the McKean-Vlasov JFBSDE 
    \begin{align}
        \left\{\arraycolsep=1pt \def\arraystretch{1.2}
		\begin{array}{ll}
    		\ddd X_t&= (Z_t +\theta_t^B)\ddd \widehat{W}_t, \quad X_0=x_0,\\
 			\ddd Y_t&= - \int_E \frac{\exp(\alpha U_t(e))-1-\alpha U_t(e)}{\alpha}\widehat{\zeta}(t,e)\lambda(\ddd e) \ddd t \\
            &\quad + Z_t \ddd \widehat{W}_t + \int_E U_t(e)\widehat{\tilde{\boldsymbol{\mu}}}(\ddd t, \ddd e),\\
			Y_T&=\rho \mathbb{E}[X_T\vert \mathcal{F}_T^{0}].
		\end{array}\label{trueCharFBSDE}\right.
    \end{align}
    It is easy to see that a reverse statement holds also true by  omitting $X$.

    To interpret the $Z$-component, we recast the McKean-Vlasov JBSDE \eqref{CharFBSDE} in terms of its associated JFBSDE \eqref{trueCharFBSDE}. The forward component $X$ is the wealth process under the MFE  strategy $Z+\theta^B$. The backward component $(Y,Z)$ solves a JBSDE with terminal condition induced by $X$. If we treat this terminal condition as exogenously given, then, by \Cref{remark:becUnterPHut}, $Z$ coincides with the optimal trading strategy for an exponential utility maximizer whose terminal liability is $Y_T=\rho \mathbb{E}[X_T \vert \mathcal{F}_T^0]$, that is the competition weight $\rho$ times the conditional expectation of the terminal wealth from the MFE strategy given the common noise information. In this sense, $Z$ constitutes the component of the equilibrium strategy arising from the mean-field interaction.
\end{remark}

\section{Transformation to an auxiliary mean-field game}\label{sec:transformation}

Working towards the proof of our main  \Cref{theorem:mainTheorem}, we next establish a one-to-one correspondence between mean-field equilibria of the MFG \eqref{MFG} and mean-field equilibria of the auxiliary MFG problem \eqref{MFG3}. For a given mean-field term $F$ the optimization problem in step 2 of the MFG \eqref{MFG} is a hedging problem for the mean-field-adjusted claim $B-\rho\mathbb{E}[B\vert\mathcal{F}_T^0]+\rho F$, consisting of the  component $B-\rho\mathbb{E}[B\vert\mathcal{F}_T^0]$, which depends on the exogenous claim $B$, and the mean field component $\rho F$. 

In the present section, we derive a related decomposition for MFE strategies by changing to the measure $\widehat{\mathbb{P}}$ associated with the optimizer $\theta^B$ of the single-agent problem \eqref{SingleAgentProbForB} for the claim-dependent component. More precisely, \Cref{lemma:UmformungSingleAgentProbInMFG} reformulates the single-agent problem under the martingale measure $\widehat{\mathbb{P}}$. This yields \Cref{proposition:MFGgleichMFG3}, which decomposes any MFE strategy of the MFG \eqref{MFG} into $\theta^B$ and a MFE strategy of the remaining auxiliary MFG \eqref{MFG3}, referring to the claim $\rho F$. This separation is the key to the financially interpretable decomposition of the MFE in \hyperref[remark:InterpretationDecompMFE]{\Cref*{remark:InterpretationDecompMFE}.1}.

\begin{lemma}\label{lemma:UmformungSingleAgentProbInMFG}
    Let $F$ be an $\R$-valued random variable. Then the equality of sets 
	\begin{align*}
    	&\argmax_{\theta \in \mathbb{H}_{BMO}^2(\mathbb{P})} \mathbb{E}[-\exp(-\alpha (X_T^{\theta}-(B-\rho \mathbb{E}[B\vert \mathcal{F}_T^0])- \rho F))\vert \mathcal{F}_0]\\
    	&=\argmax_{\theta \in \mathbb{H}_{BMO}^2(\mathbb{P})} \widehat{\mathbb{E}}[-\exp(-\alpha (\int_0^T \theta_s-\theta_s^B \ddd \widehat{W}_s- \rho F))\vert \mathcal{F}_0]
	\end{align*}
	holds, where $\widehat{\mathbb{E}}$ denotes the expectation under the measure $\widehat{\mathbb{P}}=\widehat{\mathbb{P}}^B$ defined in \eqref{changeOfMEasure} for $\beta=B$, $\widehat{W}$ is the $\widehat{\mathbb{P}}$-Brownian motion defined in \eqref{DefWHat}, and $\theta^B$ denotes the optimal strategy for the single-agent optimization problem \eqref{SingleAgentProbForB}. This strategy is unique and given by $\theta^B=Z^B+\frac{1}{\alpha}\varphi$, where $(Y^B,Z^B,U^B)\in \mathbb{S}^\infty(\mathbb{P})\times \mathbb{H}_{\text{BMO}}^2(\mathbb{P})\times \mathcal{L}_{\boldsymbol{\nu}}^2(\mathbb{P})$ with $U^B$ bounded is the solution to the JBSDE \eqref{trueCharBSDE} for $\beta=B$.
\end{lemma} 

Note that conditional expectations of non-positive random variables are always well-defined. Therefore, no integrability assumptions on $F$ in \Cref{lemma:UmformungSingleAgentProbInMFG} are required. 
		
\begin{proof}[Proof of \Cref{lemma:UmformungSingleAgentProbInMFG}]
	Let $\theta\in \mathbb{H}_{BMO}^2(\mathbb{P})$. We have the identity
	\begin{align}
		\lefteqn{\mathbb{E}\left[-\exp(-\alpha (X_T^{\theta}-(B-\rho \mathbb{E}[B\vert \mathcal{F}_T^0])- \rho F))\vert \mathcal{F}_0\right] }\notag\\
		&=\mathbb{E}\bigg[\frac{-e^{-\alpha (X_T^{\theta}-(B-\rho \mathbb{E}[B\vert \mathcal{F}_T^0])- \rho F)}}{e^{-\alpha (X_T^{B}-(B-\rho \mathbb{E}[B\vert \mathcal{F}_T^0]))}} \exp(-\alpha (X_T^{B}-(B-\rho \mathbb{E}[B\vert \mathcal{F}_T^0]))) \bigg\vert \mathcal{F}_0\bigg],\label{UmformenStep2}
	\end{align}
	where $X^B$ is the wealth process for the optimal strategy $\theta^B$ for the single-agent optimization problem \eqref{SingleAgentProbForB}. According to \Cref{proposition:BechUnique}, the unique optimal strategy $\theta^B$ for the single-agent problem \eqref{SingleAgentProbForB} is given by $\theta^B=Z^B+\frac{1}{\alpha}\varphi$, where $(Y^B,Z^B,U^B)\in \mathbb{S}^\infty(\mathbb{P})\times \mathbb{H}_{\text{BMO}}^2(\mathbb{P})\times \mathcal{L}_{\boldsymbol{\nu}}^2(\mathbb{P})$ with $U^B$ bounded is the solution to the JBSDE \eqref{trueCharBSDE} for $\beta=B$. The utility $\exp(-\alpha (X_T^B-(B-\rho \mathbb{E}[B\vert \mathcal{F}_T^0]))$ contains our change of measure to $\widehat{\mathbb{P}}=\widehat{\mathbb{P}}^B$ (defined in \eqref{changeOfMEasure} for $\beta=B$). Indeed, analogous calculations as in the proof of \cite[Thm.4.1]{becherer2006bounded} yield
	\begin{align}
		&\exp(-\alpha (X_T^B-(B-\rho \mathbb{E}[B\vert \mathcal{F}_T^0]))\notag\\
		&= e^{-\alpha(x_0-Y_0^B)}\cdot \mathscr{E}\left(-\int_0^\cdot \varphi_s \ddd W_s+ \int_0^\cdot\int_E \exp(\alpha U_s^B (e))-1 \tilde{\boldsymbol{\mu}}(\ddd s, \ddd e)\right)_T. \label{StochExpHerleitung}
	\end{align}
	The stochastic exponential in \eqref{StochExpHerleitung} is the Radon–Nikodym density of our change of measure in \eqref{changeOfMEasure} for $\beta=B$. By inserting \eqref{StochExpHerleitung} into \eqref{UmformenStep2}, and using boundedness of $\alpha$ and $Y^B$, and that $x_0\in L^2(\mathcal{A})$ is $\mathcal{F}_0$-measurable and finite, we obtain the claim.
\end{proof}
	
\begin{proposition}\label{proposition:MFGgleichMFG3}
	There is a one-to-one relationship between mean-field equilibria $\widetilde{\theta}$ to the MFG \eqref{MFG} and mean-field equilibria $\widetilde{\theta}^{\text{aux}}$ to the auxiliary MFG 
 	\begin{align} 
		\begin{cases}
			\text{1. given a real-valued random variable } F,\\ 
			\text{2. }\text{find } \widetilde{\theta}^{\text{aux}}\in \underset{\theta \in \mathbb{H}_{\text{BMO}}^2(\widehat{\mathbb{P}})}{\argmax} \widehat{\mathbb{E}}[-\exp(-\alpha (\int_0^T \theta_s\ddd \widehat{W}_s-\rho F))\vert \mathcal{F}_0].\\
			\text{3. Find a fixed point such that } F= \mathbb{E}[x_0+\int_0^T \widetilde{\theta}^{\text{aux}}_s+\theta_s^B \ddd \widehat{W}_s \vert \mathcal{F}_T^0].
		\end{cases}\label{MFG3}
 	\end{align}
	This relationship is given by $\widetilde{\theta}= \widetilde{\theta}^{\text{aux}}+\theta^B$ with $\theta^B$ from \Cref{lemma:UmformungSingleAgentProbInMFG}.
\end{proposition}

\begin{proof}
    This is a direct consequence of \Cref{lemma:UmformungSingleAgentProbInMFG} and the fact that $\mathbb{H}_{\text{BMO}}^2(\mathbb{P})=\mathbb{H}_{\text{BMO}}^2(\widehat{\mathbb{P}})$ by \hyperref[remark:PopertiesChangeOfMeasure]{\Cref*{remark:PopertiesChangeOfMeasure}.3}.
\end{proof}

\section{MFE of auxiliary game and proof of main theorem}\label{sec:CharacterizationMFE}

This section establishes the proof of our main \Cref{theorem:mainTheorem} on existence and uniqueness of the MFE for the MFG \eqref{MFG}, together with its construction. To this end, \Cref{proposition:CharUnderAssFBSDEhasSolu} shows a one-to-one correspondence between solutions of the McKean-Vlasov JBSDE \eqref{CharFBSDE} and mean-field equilibria of the auxiliary MFG \eqref{MFG3}; with MFE given by the $Z$-component. \Cref{lemma:SoluMKVbsde} proves well-posedness of the McKean-Vlasov JBSDE \eqref{CharFBSDE} and provides a construction of its solution in terms of the $\Pi$-projection (notation $\Pi$ introduced in \Cref{rn:defPi}) of the optimal strategy of the single-agent problem \eqref{SingleAgentProbForB}. In \Cref{lemma:ExistenceUniquenessOfSingleBSDE} we establish well-posedness of the JBSDEs \eqref{BSDESingle}, thereby verifying the assumptions of \Cref{proposition:CharUnderAssFBSDEhasSolu}. Combining these results with \Cref{proposition:MFGgleichMFG3}, which links the mean-field equilibria of the MFG \eqref{MFG} and \eqref{MFG3}, yields \Cref{theorem:mainTheorem}. 

While our analysis has benefited from those by \cite{fu2023mean1,tangpi2024optimal} for the purely Brownian case without jumps, there are important differences in how well-posedness is obtained. In \cite[Prop.6.1]{tangpi2024optimal} and \cite[Thm.3.8]{fu2023mean1}, well-posedness of the relevant McKean-Vlasov (F)BSDE is established via Banach's fixed point argument, under a weak interaction assumption, which is restrictive in requiring mean-field interaction to be sufficiently small (competition weight close to $0$). In contrast, we establish a one-to-one correspondence between solutions of the McKean-Vlasov JBSDE \eqref{CharFBSDE} and solutions of the JBSDE \eqref{auxiliaryJBSDE} without mean-field. 
This provides well-posedness without imposing a weak interaction condition and, moreover, a constructive description of the solution to the McKean-Vlasov JBSDE \eqref{CharFBSDE} in \Cref{lemma:SoluMKVbsde}. 

\begin{proposition}\label{proposition:CharUnderAssFBSDEhasSolu}
	Let $\theta^B=Z^B +\frac{1}{\alpha}\varphi$ be as in \Cref{lemma:UmformungSingleAgentProbInMFG} and let for each mean-field equilibrium $\bar{\theta}$ to the auxiliary mean-field game \eqref{MFG3} the JBSDE
	\begin{align}
		\left\{\arraycolsep=1pt \def\arraystretch{1.2}
		\begin{array}{ll}
			\ddd \bar{Y}_t&= - \int_E \frac{\exp(\alpha \bar{U}_t(e))-1-\alpha \bar{U}_t(e)}{\alpha}\widehat{\zeta}(t,e)\lambda(\ddd e) \ddd t + \bar{Z}_t \ddd \widehat{W}_t + \int_E \bar{U}_t(e)\widehat{\tilde{\boldsymbol{\mu}}}(\ddd t, \ddd e),\\
			\bar{Y}_T&=\rho \mathbb{E}[x_0+\int_0^T (\bar{\theta}_s +\theta_s^B)\ddd \widehat{W}_s \vert \mathcal{F}_T^{0}] 
		\end{array}\label{BSDESingle}\right.
	\end{align} 
	have a solution $(\bar{Y},\bar{Z},\bar{U})$ in $ \mathbb{S}^2(\widehat{\mathbb{P}})\times \mathbb{H}_{\text{BMO}}^2(\widehat{\mathbb{P}})\times \mathcal{L}_{\widehat{\boldsymbol{\nu}}}^2(\widehat{\mathbb{P}})$ with $\bar{U}$ bounded. Then there exists a solution $(Y,Z, U)\in \mathbb{S}^2(\widehat{\mathbb{P}})\times \mathbb{H}_{\text{BMO}}^2(\widehat{\mathbb{P}})\times \mathcal{L}_{\widehat{\boldsymbol{\nu}}}^2(\widehat{\mathbb{P}})$, with $U$ bounded, to the McKean-Vlasov JBSDE \eqref{CharFBSDE} if and only if the auxiliary mean-field game \eqref{MFG3} has a mean-field equilibrium $\widetilde{\theta}^{\text{aux}}$. These McKean-Vlasov JBSDE solutions to \eqref{CharFBSDE} and respective mean-field equilibria are related by the identity $\widetilde{\theta}^{\text{aux}}=Z$.
\end{proposition}

\begin{proof} 
	First, we show that if we have a solution for the McKean-Vlasov JBSDE \eqref{CharFBSDE}, then we have a MFE $\widetilde{\theta}^{\text{aux}}$ and we can write $\widetilde{\theta}^{\text{aux}}$ as in the proposition. For this, let $(Y,Z,U)\in \mathbb{S}^2(\widehat{\mathbb{P}})\times \mathbb{H}_{\text{BMO}}^2(\widehat{\mathbb{P}})\times \mathcal{L}_{\widehat{\boldsymbol{\nu}}}^2(\widehat{\mathbb{P}})$ with $U$ bounded be a solution of the McKean-Vlasov JBSDE $\eqref{CharFBSDE}$. Let $F:=\mathbb{E}[x_0+\int_0^T (Z_s +\theta_s^B)\ddd \widehat{W}_s \vert \mathcal{F}_T^{0}]$. Then, the process $(Y,Z,U)$ solves the JBSDE
	\begin{align*}
		\ddd Y_t&= - \int_E \frac{\exp(\alpha U_t(e))-1-\alpha U_t(e)}{\alpha}\widehat{\zeta}(t,e)\lambda(\ddd e) \ddd t + Z_t \ddd \widehat{W}_t + \int_E U_t(e)\widehat{\tilde{\boldsymbol{\mu}}}(\ddd t, \ddd e),
	\end{align*} 
 	with terminal condition $Y_T=\rho F$. By \Cref{example:IntInL2}, $F$ is in $L^2(\mathcal{F}_T,\widehat{\mathbb{P}})$. Thus, \Cref{remark:becUnterPHut} yields that $\widetilde{\theta}^{\text{aux}}=Z\in \mathbb{H}_{\text{BMO}}^2(\widehat{\mathbb{P}})$ is an optimal strategy for the optimization problem of step~2 of the MFG \eqref{MFG3}, given $F$. It satisfies the fixed point condition in step~3 of the MFG \eqref{MFG3}. Thus, the strategy $\widetilde{\theta}^{\text{aux}}$ is a MFE to the MFG \eqref{MFG3}.

	Conversely, let $\widetilde{\theta}^{\text{aux}}$ be a MFE to the MFG \eqref{MFG3}. Then $\widetilde{\theta}^{\text{aux}}$ solves the single-agent maximization problem from step 2 in \eqref{MFG3} for exogenously given $F:=\mathbb{E}[x_0+\int_0^T \widetilde{\theta}^{\text{aux}}_s+\theta_s^B \ddd \widehat{W}_s\vert \mathcal{F}_T^0]$. According to \Cref{example:IntInL2}, $F$ is again in $L^2(\mathcal{F}_T,\widehat{\mathbb{P}})$. Thus, according to \Cref{remark:becUnterPHut}, the strategy $\widetilde{\theta}^{\text{aux}}$ is given by
	\begin{equation}
		\widetilde{\theta}^{\text{aux}}=Z \in \mathbb{H}_{\text{BMO}}^2(\widehat{\mathbb{P}}),\label{OS}
	\end{equation} 
	where $(Y,Z,U)\in \mathbb{S}^2(\widehat{\mathbb{P}})\times \mathbb{H}_{\text{BMO}}^2(\widehat{\mathbb{P}}) \times \mathcal{L}_{\boldsymbol{\nu}}^2(\widehat{\mathbb{P}})$ with $U$ bounded is the solution of the JBSDE \eqref{BSDESingle} with $\bar{\theta}=\widetilde{\theta}^{\text{aux}}$, which exists by the assumption of this proposition. By inserting the representation \eqref{OS} for the strategy $\widetilde{\theta}^{\text{aux}}$ into the terminal condition of the JBSDE \eqref{BSDESingle}, $(Y,Z,U)$ is a solution of the McKean-Vlasov JBSDE \eqref{CharFBSDE} with bounded $U$.
\end{proof}

Next, we prove \Cref{lemma:SoluMKVbsde}, which provides an explicit unique solution to the McKean-Vlasov JBSDE \eqref{CharFBSDE}, and \Cref{lemma:ExistenceUniquenessOfSingleBSDE}, which verifies the conditions for \Cref{proposition:CharUnderAssFBSDEhasSolu}.  
Finally, we combine the results to prove our main \Cref{theorem:mainTheorem}. In the following, we use the $\Pi$-projection introduced in \Cref{rn:defPi}.

\begin{lemma}\label{lemma:SoluMKVbsde}
	The McKean-Vlasov JBSDE \eqref{CharFBSDE} has a unique solution $(Y,Z,U)$ in $\mathbb{S}^2(\widehat{\mathbb{P}})\times \mathbb{H}_{\text{BMO}}^2(\widehat{\mathbb{P}}) \times \mathcal{L}_{\boldsymbol{\nu}}^2(\widehat{\mathbb{P}})$ with $U$ bounded. It is, more specifically, given by
    \begin{equation}
        (Y,Z,U)=\left(\rho \mathbb{E}[x_0]+\int_0^\cdot \frac{\rho}{1-\mathbb{E}[\rho]}\Pi(\theta^B)_s \ddd\widehat{W}_s, \frac{\rho}{1-\mathbb{E}[\rho]}\Pi(\theta^B), 0\right). \label{MKVJBSDEReprPiThetaB}
    \end{equation}
\end{lemma}

\begin{proof}
    To begin, exploiting the linear dependence on $X$ in the terminal condition of the McKean-Vlasov JBSDE \eqref{CharFBSDE}, we establish via the mutually inverse mappings \eqref{DefTildeYZU} and \eqref{DefYZU} a one-to-one correspondence between its solutions and those of the auxiliary JBSDE \eqref{auxiliaryJBSDE}, whose terminal condition is bounded. Afterwards, we show that the auxiliary JBSDE \eqref{auxiliaryJBSDE} is well-posed and has the trivial solution, which yields by the mapping \eqref{DefYZU} the explicit representation \eqref{MKVJBSDEReprPiThetaB} of the unique solution to the McKean-Vlasov JBSDE \eqref{CharFBSDE}.

	At first, let us show how solutions of the McKean-Vlasov JBSDE provide solutions to the auxiliary JBSDE by suitable transformation. For this, let $(Y,Z,U)\in \mathbb{S}^2(\widehat{\mathbb{P}})\times \mathbb{H}_{\text{BMO}}^2(\widehat{\mathbb{P}})\times \mathcal{L}_{\widehat{\boldsymbol{\nu}}}^2(\widehat{\mathbb{P}})$ with $U$ bounded be a solution to the McKean-Vlasov JBSDE \eqref{CharFBSDE}. The terminal condition 
 	\begin{equation}
		Y_T= \rho \mathbb{E}\bigg[x_0+\int_0^T (Z_s +\theta_s^B)\ddd \widehat{W}_s \vert \mathcal{F}_T^{0}\bigg]
		=\rho\left( \mathbb{E}\left[x_0\right]+ \int_0^T \Pi(Z +\theta^B)_{s} \ddd \widehat{W}_s\right),\label{EndbedingungYTL}
	\end{equation}
	is transformed by using in the second equality that $x_0$ is $\mathcal{A}$-measurable, $\mathcal{A}$ is independent of $\mathcal{F}_T^0$, that $Z\in \mathbb{H}_{\text{BMO}}^2(\widehat{\mathbb{P}})=\mathbb{H}_{\text{BMO}}^2(\mathbb{P})$ according to \Cref{remark:PopertiesChangeOfMeasure}, \Cref{lemma:La} and that $\varphi$ is $\mathbb{F}^0$-predictable. We define
    \begin{equation}
			\tilde{Y}:= Y-\rho \mathbb{E}[x_0]-\rho \int_0^\cdot \Pi(Z +\theta^B)_s \ddd \widehat{W}_s, \quad 
			\tilde{Z}:=Z-\rho \ \Pi(Z +\theta^B) \quad\text{and}\quad 
			\tilde{U}:=U.\label{DefTildeYZU}
	\end{equation}
    Using the definition \eqref{DefTildeYZU} together with \eqref{EndbedingungYTL} we obtain the auxiliary JBSDE
    \begin{align}
		\tilde{Y}_t&=0+\int_t^T \int_E \frac{\exp(\alpha \tilde{U}_s(e))-1-\alpha \tilde{U}_s(e)}{\alpha}\widehat{\zeta}(s,e)\lambda(\ddd e) \ddd s \label{auxiliaryJBSDE}\\
		&\qquad \quad -\int_t^T \tilde{Z}_s \ddd \widehat{W}_s -\int_t^T \int_E \tilde{U}_s(e)\widehat{\tilde{\boldsymbol{\mu}}}(\ddd s, \ddd e). \notag
	\end{align}
	Recall that $\mathbb{H}_{\text{BMO}}^2(\mathbb{P})=\mathbb{H}_{\text{BMO}}^2(\widehat{\mathbb{P}})$. By \Cref{lemma:zBMOozBMO} and from the linearity of $\mathbb{H}_{\text{BMO}}^2$ it follows that $\Pi(Z+\theta^B)\in \mathbb{H}_{\text{BMO}}^2$. Since $\rho$ is bounded and $\mathcal{F}_0$-measurable, the representation in \eqref{DefTildeYZU} implies $\tilde{Y}\in \mathbb{S}^2(\widehat{\mathbb{P}})$ and likewise $\tilde{Z}\in \mathbb{H}_{\text{BMO}}^2(\widehat{\mathbb{P}})$. Thus, $(\tilde{Y}, \tilde{Z}, \tilde{U})\in \mathbb{S}^2(\widehat{\mathbb{P}})\times \mathbb{H}_{\text{BMO}}^2(\widehat{\mathbb{P}})\times \mathcal{L}_{\widehat{\boldsymbol{\nu}}}^2(\widehat{\mathbb{P}})$ from \eqref{DefTildeYZU} solves the auxiliary JBSDE \eqref{auxiliaryJBSDE} with $\tilde{U}$ bounded.

	Next, we show for the reversed correspondence how, if the auxiliary JBSDE \eqref{auxiliaryJBSDE} has a solution, we then can obtain from it by suitable transformations a solution for the McKean-Vlasov JBSDE \eqref{CharFBSDE}. Let $(\tilde{Y}, \tilde{Z}, \tilde{U})\in \mathbb{S}^2(\widehat{\mathbb{P}})\times \mathbb{H}_{\text{BMO}}^2(\widehat{\mathbb{P}})\times \mathcal{L}_{\widehat{\boldsymbol{\nu}}}^2(\widehat{\mathbb{P}})$ with $\tilde{U}$ bounded now be a solution to the auxiliary JBSDE \eqref{auxiliaryJBSDE}. Inspired by the mapping  \eqref{DefTildeYZU}, we are looking to find a unique solution $z\in \mathbb{H}_{\text{BMO}}^2(\widehat{\mathbb{P}})$ of the equation
    \begin{equation}
        \tilde{Z}= z-\rho \ \Pi(z +\theta^B).\label{GleichungZTildeZ}
    \end{equation}
    According to \hyperref[remark:PopertiesChangeOfMeasure]{\Cref*{remark:PopertiesChangeOfMeasure}.3}, $\tilde{Z}\in \mathbb{H}_{\text{BMO}}^2(\widehat{\mathbb{P}})=\mathbb{H}_{\text{BMO}}^2(\mathbb{P})$. By taking the conditional expectation $\Pi$ on both sides in \eqref{GleichungZTildeZ}, using that $\rho$ is $\mathcal{A}$-measurable, $\mathcal{A}$ is independent of $\mathcal{F}_T^0$, $\mathbb{E}[\rho]\neq 1$ and linearity of $\Pi$, we obtain
	\begin{equation}
		\Pi(z)= \frac{\Pi(\tilde{Z})+ \mathbb{E}[\rho] \cdot \Pi(\theta^B)}{1-\mathbb{E}[\rho]}. \label{OptProjZ}
	\end{equation}
	Thanks to linearity of the projection $\Pi$, by \eqref{OptProjZ} the unique solution $z$ to \eqref{GleichungZTildeZ} is given by
	\begin{align}
		z
		=\tilde{Z}+ \rho \ \Pi(z +\theta^B)
		=\tilde{Z}+ \rho \ \frac{\Pi(\tilde{Z})+\mathbb{E}[\rho] \cdot \Pi(\theta^B)}{1-\mathbb{E}[\rho]} + \rho \ \Pi(\theta^B)=: G(\tilde{Z}). \label{GleichungZ}
	\end{align}
	With \Cref{lemma:zBMOozBMO}, $\mathbb{H}_{\text{BMO}}^2(\mathbb{P})=\mathbb{H}_{\text{BMO}}^2(\widehat{\mathbb{P}})$ and the linearity of $\mathbb{H}_{\text{BMO}}^2$, it follows that $G(\tilde{Z})\in \mathbb{H}_{\text{BMO}}^2(\widehat{\mathbb{P}})$. Next, we define $Y,Z$ and $U$ by
    \begin{equation}		
        Y:= \tilde{Y}+ \rho \mathbb{E}[x_0]+ \rho \int_0^\cdot \Pi(G(\tilde{Z}) +\theta^B)_{s} \ddd \widehat{W}_s\,,\quad\:
        Z:=G(\tilde{Z})\,,\quad\:
        U:=\tilde{U}
        \label{DefYZU}
	\end{equation}
	and again get $Y\in \mathbb{S}^2(\widehat{\mathbb{P}})$ and $Z\in \mathbb{H}_{\text{BMO}}^2(\widehat{\mathbb{P}})$. By definitions \eqref{DefYZU} and equality \eqref{GleichungZ} we obtain
	\begin{align*}
		Y_t&=\rho \mathbb{E}[x_0]+\int_t^T \int_E \frac {1}{\alpha}\left( {\exp(\alpha \tilde{U}_s(e))-1-\alpha \tilde{U}_s(e)} \right) \widehat{\zeta}(s,e)\lambda(\ddd e) \ddd s \notag\\
		&\qquad \qquad + \rho \int_0^T \Pi(Z +\theta^B)_{s} \ddd \widehat{W}_s -\int_t^T Z_s \ddd \widehat{W}_s -\int_t^T \int_E \tilde{U}_s(e)\widehat{\tilde{\boldsymbol{\mu}}}(\ddd s, \ddd e).
	\end{align*}
	Using now the arguments like in \eqref{EndbedingungYTL} (in reverse order) and the definition of $U$ from \eqref{DefYZU}, yields the McKean-Vlasov JBSDE \eqref{CharFBSDE}. That means, $(Y,Z,U)$ in $\mathbb{S}^2(\widehat{\mathbb{P}})\times \mathbb{H}_{\text{BMO}}^2(\widehat{\mathbb{P}})\times \mathcal{L}_{\widehat{\boldsymbol{\nu}}}^2(\widehat{\mathbb{P}})$ given by \eqref{DefYZU} (with $U$ being bounded) is a solution of the McKean-Vlasov JBSDE \eqref{CharFBSDE}. 
	
	Overall, we have shown that the mutually inverse mappings \eqref{DefTildeYZU} and \eqref{DefYZU} establish a one-to-one relation between solutions for the McKean-Vlasov JBSDE \eqref{CharFBSDE} and the auxiliary JBSDE \eqref{auxiliaryJBSDE}. Next, we prove well-posedness 
    of the auxiliary JBSDE \eqref{auxiliaryJBSDE}, and thus well-posedness of the McKean-Vlasov JBSDE \eqref{CharFBSDE}. Clearly the terminal condition of the auxiliary JBSDE \eqref{auxiliaryJBSDE} is bounded. Since $\alpha$ is bounded and greater than $0$, $u\mapsto g(u):=(\exp(\alpha u)-1-\alpha u)/\alpha$ is absolutely continuous in $u$. The density function $g'$ is strictly greater than $-1$ and locally bounded in $u$, uniformly on $\Omega$. Thus, the auxiliary JBSDE \eqref{auxiliaryJBSDE} has a unique solution $(\tilde{Y},\tilde{Z}, \tilde{U})\in \mathbb{S}^\infty(\widehat{\mathbb{P}})\times \mathcal{L}_T^2(\widehat{\mathbb{P}})\times \mathcal{L}_{\widehat{\boldsymbol{\nu}}}^2(\widehat{\mathbb{P}})$ by \cite[Prop.4.3]{becherer2019monotone}, and for $\tilde{U}$ a bounded representative can be chosen in $\mathcal{L}_{\widehat{\boldsymbol{\nu}}}^2(\widehat{\mathbb{P}})$ (by \cite[Lem.2.2]{becherer2019monotone}). According to \cite[Lem.3.4]{becherer2006bounded}, $\tilde{Z}\in \mathbb{H}_{\text{BMO}}^2(\widehat{\mathbb{P}})$ also follows. Since $\tilde{Y}$ is in $\mathbb{S}^\infty(\widehat{\mathbb{P}})$, we also have $\tilde{Y} \in \mathbb{S}^2(\widehat{\mathbb{P}})$. Hence, we have a solution $(\tilde{Y},\tilde{Z}, \tilde{U})\in \mathbb{S}^2(\widehat{\mathbb{P}})\times \mathbb{H}_{\text{BMO}}^2(\widehat{\mathbb{P}})\times \mathcal{L}_{\widehat{\boldsymbol{\nu}}}^2(\widehat{\mathbb{P}})$ of the auxiliary JBSDE \eqref{auxiliaryJBSDE} with $\tilde{U}$ being bounded. The uniqueness of the solution follows from the boundedness of the $U$-component, since we can then, by a truncation argument, regard the generator of the auxiliary JBSDE \eqref{auxiliaryJBSDE} as being Lipschitz continuous in its (only) argument $u$, whereby uniqueness follows even in the larger space $ \mathbb{S}^2(\widehat{\mathbb{P}})\times \mathcal{L}_{T}^2(\widehat{\mathbb{P}})\times \mathcal{L}_{\widehat{\boldsymbol{\nu}}}^2(\widehat{\mathbb{P}})$ with $\tilde{U}$ being bounded (see \cite[Prop.3.3]{becherer2006bounded}).

    Finally, we observe that $(0, 0, 0)$ is the solution of the auxiliary JBSDE \eqref{auxiliaryJBSDE} since the generator vanishes for $U\equiv 0$. Using the representation \eqref{DefYZU} of the solution of the McKean-Vlasov JBSDE \eqref{CharFBSDE} in terms of the solution of the auxiliary JBSDE \eqref{auxiliaryJBSDE}, we thus obtain the desired representation \eqref{MKVJBSDEReprPiThetaB} 
    of the solution of the McKean-Vlasov JBSDE \eqref{CharFBSDE}.
\end{proof}

\begin{lemma}\label{lemma:ExistenceUniquenessOfSingleBSDE}
	For any strategy $\bar{\theta}\in \mathbb{H}_{\text{BMO}}^2(\widehat{\mathbb{P}})$, the JBSDE \eqref{BSDESingle} with terminal condition depending on $\bar{\theta}$ has a unique solution $(Y,Z,U)\in \mathbb{S}^2(\widehat{\mathbb{P}})\times \mathbb{H}_{\text{BMO}}^2(\widehat{\mathbb{P}})\times \mathcal{L}_{\widehat{\boldsymbol{\nu}}}^2(\widehat{\mathbb{P}})$ with bounded $U$.
\end{lemma}

The proof is very similar to that of \Cref{lemma:SoluMKVbsde}, hence we omit the details here. The main idea is to establish a one-to-one connection between solutions $(Y,Z,U)$ of the actual JBSDE \eqref{BSDESingle} and solutions $(\tilde{Y},\tilde{Z}, \tilde{U})$ of the auxiliary JBSDE \eqref{auxiliaryJBSDE} by 
\begin{align*}
	\tilde{Y}= Y-\rho \mathbb{E}[x_0]- \rho \int_0^\cdot  \Pi(\bar{\theta} +\theta^B)_{s}  \ddd \widehat{W}_s,\quad 
	\tilde{Z}= Z-\rho \  \Pi(\bar{\theta} +\theta^B) \quad \text{and}\quad 
	\tilde{U}=U.
\end{align*}
    
Finally, we combine the previous results to conclude the proof of \Cref{theorem:mainTheorem}.

\begin{proof}[Proof of \Cref{theorem:mainTheorem}]
    First, by \Cref{lemma:UmformungSingleAgentProbInMFG} the unique optimal strategy $\theta^B$ to the single-agent optimization problem \eqref{SingleAgentProbForB} is given by $\theta^B=Z^B+\frac{1}{\alpha}\varphi$, where $(Y^B,Z^B,U^B)\in \mathbb{S}^\infty(\mathbb{P})\times \mathbb{H}_{\text{BMO}}^2(\mathbb{P})\times \mathcal{L}_{\boldsymbol{\nu}}^2(\mathbb{P})$ with $U^B$ bounded is the solution to the JBSDE \eqref{trueCharBSDE} for $\beta:=B$. Next, by \Cref{lemma:SoluMKVbsde}, 
    the McKean-Vlasov JBSDE \eqref{CharFBSDE} has a unique solution $(Y,Z, U)\in \mathbb{S}^2(\widehat{\mathbb{P}})\times \mathbb{H}_{\text{BMO}}^2(\widehat{\mathbb{P}})\times \mathcal{L}_{\widehat{\boldsymbol{\nu}}}^2(\widehat{\mathbb{P}})$ such that $U$ is bounded and it is given by 
    \begin{align*}
        (Y,Z, U)=\left(\rho \mathbb{E}[x_0]+\int_0^\cdot \frac{\rho}{1-\mathbb{E}[\rho]}\Pi(\theta^B)_s \ddd\widehat{W}_s\;,\; \frac{\rho}{1-\mathbb{E}[\rho]}\Pi(\theta^B)\;, \;0\right).
    \end{align*}
    By \Cref{lemma:ExistenceUniquenessOfSingleBSDE}, the assumptions of \Cref{proposition:CharUnderAssFBSDEhasSolu} are satisfied. Hence, \Cref{proposition:CharUnderAssFBSDEhasSolu} yields existence and uniqueness of the MFE $\widetilde{\theta}^{\text{aux}}$ for the auxiliary MFG \eqref{MFG3} and that it is given by $\widetilde{\theta}^{\text{aux}}=Z$. By \Cref{proposition:MFGgleichMFG3}, the unique MFE $\widetilde{\theta}$ to the MFG \eqref{MFG} is then given by $\widetilde{\theta}=\theta^B+\widetilde{\theta}^{\text{aux}}=\theta^B+Z=\frac{1}{\alpha} \varphi + Z^B + Z= \theta^B+ \frac{\rho}{1-\mathbb{E}[\rho]} \Pi(\theta^B)$.
\end{proof}

\section{Examples, computations and illustration of MFE} \label{sec:PDE}

The aim of this section is to illustrate our general (non-Markovian) MFE results through concrete examples, using PDEs to characterize and compute the MFE and its comparative statics. In \Cref{subsec:PDEcharMFE} we characterize the MFE by a PDE system (see \eqref{SystemOfPDEs}) in a setting with a finite-state Markov chain process and Markov-modulated asset dynamics. In \Cref{subsec:IllustrationMFE} this PDE characterization is used to numerically illustrate the MFE (\Cref{fig:MFE}) and its decomposition (\Cref{fig:DecompositionOfMFEstrat}) according to \Cref{theorem:mainTheorem}, as well as the impact of the common noise on the MFE (\Cref{fig:ImpactCommonNoise}) in a concrete example. 

For simplicity, let us assume for this section that the market contains a single risky asset and a riskless numeraire asset (with unit price one), i.e.\ a savings account with zero interest. Let $L$ be a Markov chain with finite state space $K$, independent of the initial $\sigma$-field $\mathcal{F}_0$ and the Brownian motion $W$. Of course, the asset dynamics of the financial market should be common information. The evolution of the (discounted) price of a single risky asset (say, a market index) is described by the SDE
\begin{align*}
    \ddd S_t = \sigma^{L_{t-}} S_t (\varphi^{L_{t-}} \ddd t + \ddd W_t), \quad t\in [0,T],  \quad S_0=1,
\end{align*}
where, for each $i\in K$, $\sigma^i>0$ and $\varphi^i\in \R$ are constants representing the volatility and market price of risk prevailing while the Markov chain $L$ is in state $i$. Since the SDE coefficients may depend on the Markov chain $L$, even European put and call options on $S$ are generally non-replicable. Thus, the market is generically incomplete even in absence of a contingent claim liability $B$ which may comprise further non-hedgeable risk factors.

\subsection{Characterization of MFE by solution of system of PDEs}\label{subsec:PDEcharMFE}

In this subsection we construct the MFE of the MFG \eqref{MFG} in the setting with a Markov-switching SDE dynamic via solutions of a PDE system \eqref{SystemOfPDEs}; see constructions \eqref{ThetaPDEoptimalStratSingleAgentProb}-\eqref{PDEmfeGeneral}. To this end, we show that the solution to the JBSDE \eqref{trueCharBSDE}, for $\beta=B$, 
admits a representation via the solution to the semilinear parabolic PDE system \eqref{SystemOfPDEs} and hence, by \Cref{theorem:mainTheorem}, the MFE is expressed in terms of the PDE solution. 

The transition rates from the Markov chain $L$ from state $i$ to state $j\in K$, $i\neq j$, are given by $q^{i\rightarrow j}\in [0,\infty)$. We denote by $\Delta L_t:=L_t-L_{t-}$ the jumps of the Markov chain $L$. As before, $B\in L^\infty(\mathcal{F}_T)$ is a bounded random variable and we assume that there are bounded continuous functions $h^i:(0,\infty)\rightarrow \R$, $i\in K$, such that $h^{L_T}(S_T)=B-\rho \mathbb{E}\left[B \vert \mathcal{F}_T^0\right]$.

In the sequel, we consider the system of interacting PDEs of reaction-diffusion type 
\begin{align}
    \begin{cases}
        v_t^i + \frac{1}{2} (\sigma^i)^2 s^2 v_{ss}^i - \frac{1}{2\alpha} \tabs{\varphi^i}^2 + \sum_{j\neq i} q^{i\rightarrow j} \cdot \frac{1}{\alpha} {(e^{\alpha (v^j - v^i)} - 1)} = 0, \quad \text{ and}\\
        v^i(T,s) = h^i(s), \quad \text{ for } s\in (0,\infty),\, t\in [0,T),\,  i\in K,
    \end{cases}\label{SystemOfPDEs}
\end{align}
which has a unique bounded classical solution $v$ in the space $C([0,T]\times (0,\infty), \R^{\tabs{K}})\cap C^{1,2}([0,T)\times (0,\infty), \R^{\tabs{K}})$ by \cite[Thm.2.4]{becherer2005classical}.

Next, we establish a link between the solution $v$ of the system of PDEs \eqref{SystemOfPDEs} and solutions of the JBSDE \eqref{trueCharBSDE} for $\beta=B$, which characterize by \Cref{theorem:mainTheorem} the optimal strategy of the single-agent optimization problem \eqref{SingleAgentProbForB}. To this end, we identify the states $i\in K$ of the Markov chain $L$ with unit vectors $e_i\in \R^{\tabs{K}}\backslash \{0\}=:E$ with jumps $\Delta L$ taking values $e^j-e^i$, $i,j\in K$ (cf.\ \cite{cohen2010comparisons}, \cite[Ex.2.1.4]{becherer2019monotone}). 
The associated integer-valued random measure for its jumps is then given by\\
$\mu(\omega, \ddd t, \ddd e)
= \sum_{s>0: \Delta L_s \neq 0} \delta_{(s,\Delta L_s)}(dt, de) 
= \sum_{i, j \in K: i \neq j}\sum_{ s>0: \Delta L_s = e_j-e_i} \delta_s(\ddd t)\delta_{e_j-e_i}(\ddd e)$ and its compensator $\nu$ given by $\nu(\omega, \ddd t, \ddd e)=\zeta(\omega,t,e)\lambda(\ddd e)\ddd t$ with
\begin{align*}
    \zeta(\omega, t, e)= \sum_{i,j\in {K}: i\neq j} q^{i\rightarrow j} \Indi_{\{L_{t-}(\omega)=i\}}\Indi_{\{e=e^j-e^i\}} \quad \text{and} \quad 
    \lambda(\ddd e)= \sum_{i,j\in {K}:j\neq i} \delta_{\{e^j-e^i\}}(\ddd e).
\end{align*}
The connection between this integer-valued random measure $\mu$ and the Markov chain $L$ is simply given by $L_t-L_0=\sum_{s \le t: \Delta L_s \neq 0} \Delta L_s= \int_0^t\int_E e \mu(\ddd s, \ddd e)$. We define processes
\begin{align}
    \begin{cases}
        Y_t^B := v^{L_t}(t, S_t),\\
        Z_t^B := \partial_s v^{L_{t-}}(t, S_t) \cdot \sigma^{L_{t-}} S_t \Indi_{[0,T)}(t),\\
        U_t^B(e):=u(e,v(t,S_t)):=\sum_{i,j\in K: i\neq j}\left(v^j(t,S_t)-v^i(t,S_t)\right)\Indi_{\{e^j-e^i\}}(e)
    \end{cases}\label{PDEBSDEidentification}
\end{align} 
for $t\in [0,T]$, and we are going to show that $(Y^B,Z^B,U^B)$ solves the JBSDE \eqref{trueCharBSDE} for $\beta=B$. Using Itô's formula, one obtains for $t\in [0,T)$ that
\begin{align}
    \ddd v^{L_t}(t,S_t)&= \left(v_t^{L_{t-}}(t, S_t)+ v_s^{L_{t-}}(t,S_t)\sigma^{L_{t-}}\varphi^{L_{t-}} S_t + \frac{1}{2} v_{ss}^{L_{t-}}(t,S_t)(\sigma^{L_{t-}})^2 S_t^2\right) \ddd t \notag \\
    &\quad + v_s^{L_{t-}} (t, S_t) \sigma^{L_{t-}} S_t \ddd W_t  + \left\{ v^{L_t}(t,S_t)-v^{L_{t-}}(t,S_t) \right\}\Indi_{\{\Delta L_t\neq 0\}}.
    \label{dvIto}
\end{align} 
For the solution to the PDE system \eqref{SystemOfPDEs}, we rewrite the equation in \eqref{SystemOfPDEs} as 
\begin{align*}
    v_t^i + \frac{1}{2} (\sigma^i)^2 s^2 v_{ss}^i = \frac{1}{2\alpha} \tabs{\varphi^i}^2 - \sum_{j\neq i} q^{i\rightarrow j} \cdot \frac{1}{\alpha} {(e^{\alpha (v^j - v^i)} - 1)}, \quad (t,s)\in [0,T)\times (0,\infty),\ i\in K.
\end{align*}
Using this and substituting \eqref{PDEBSDEidentification} into \eqref{dvIto}, noting that the Markov chain $L$ is quasi left-continuous and thus has a.s.\ no jump at the deterministic time $T$, one obtains 
\begin{align}
    \ddd Y_t^B &= Z_t^B \varphi^{L_{t-}} \ddd t + \left(v_t^{L_{t-}}(t, S_t)+ \frac{1}{2} v_{ss}^{L_{t-}}(t,S_t)(\sigma^{L_{t-}})^2 S_t^2\right) \ddd t + Z_t^B \ddd W_t\notag\\
    & \quad + \int_E U_t^B(e) \tilde{\mu} (\ddd t, \ddd e) + \int_E U_t^B(e) \nu (\ddd t, \ddd e), \notag\\
    &=(Z_t^B \varphi_t + \frac{1}{2\alpha}\vert \varphi_t \vert^2) \ddd t - \int_E \frac{\exp(\alpha U_t^B(e))-1-\alpha U_t^B(e)}{\alpha}\zeta(t,e)\lambda(\ddd e) \ddd t\notag\\
	&\quad ~ + Z_t^B \ddd W_t + \int_E U_t^B(e)\tilde{\mu}(\ddd t, \ddd e), \quad t\in [0,T],\label{YDynFromPDE}
\end{align}
writing $\varphi_t$ for $\varphi^{L_{t-}}$. As the classical solution $v$ to PDE system \eqref{SystemOfPDEs} is bounded, the processes $Y^B$ and $U^B$ defined via \eqref{PDEBSDEidentification} are  bounded, hence $(Y^B, U^B)$ is in $ \mathbb{S}^2(\mathbb{P})\times\mathcal{L}_\nu^2(\mathbb{P})$ with $U^B$ being bounded. Using arguments as in the proof of \cite[Thm.4.4]{becherer2005classical}, one can also show that $Z^B$ defined in \eqref{PDEBSDEidentification} is in $\mathbb{H}_{\text{BMO}}^2(\mathbb{P})$. Hence, $(Y^B, Z^B,U^B)\in \mathbb{S}^2(\mathbb{P})\times\mathbb{H}_{\text{BMO}}^2(\mathbb{P})\times \mathcal{L}_\nu^2(\mathbb{P})$ defined in \eqref{PDEBSDEidentification} is the solution of the JBSDE \eqref{trueCharBSDE} for $\beta=B$. 

The optimal strategy $\theta^B$ for the single-agent optimization problem \eqref{SingleAgentProbForB} is, according to \Cref{theorem:mainTheorem}, and noting that \eqref{PDEBSDEidentification} identifies $Z^B$, given by
\begin{align}
    \theta_t^B=Z_t^B+\frac{1}{\alpha} \varphi_t = v_s^{L_{t-}}(t,S_t)\sigma^{L_{t-}} S_t \Indi_{[0,T)}(t)+ \frac{1}{\alpha} \varphi^{L_{t-}}, \quad t\in [0,T].\label{ThetaPDEoptimalStratSingleAgentProb}
\end{align}
The $\theta^B$ from \eqref{ThetaPDEoptimalStratSingleAgentProb}
corresponds to the strategy
\begin{align} 
    \vartheta_t^B= v_s^{L_{t-}}(t,S_t) \Indi_{[0,T)}(t) + \frac{\varphi^{L_{t-}}}{\alpha \sigma^{L_{t-}} }\cdot \frac{1}{S_t}, \quad t\in [0,T] \label{VarthetaPDEoptimalStratSingleAgentProb},
\end{align} 
which describes the number of shares being held  over time. For this correspondence, we recall the relation \eqref{parametrizationVarthetaTheta} (with $\sigma^{L_{t-}}$ for $\sigma_t$), which is also used frequently in the sequel. By \Cref{theorem:mainTheorem}, the MFE $\widetilde{\theta}$\ ($\widetilde{\vartheta}$) to the MFG \eqref{MFG} is given by
\begin{align}
    &\widetilde{\theta}= \theta^B +\frac{\rho}{1-\mathbb{E}[\rho]}\Pi(\theta^B) \quad \text{and} \quad 
    \widetilde{\vartheta}= \vartheta^B +\frac{\rho}{1-\mathbb{E}[\rho]}\Pi(\vartheta^B). \label{PDEmfeGeneral}
\end{align}

\subsection{Illustration of the MFE and common noise effects} \label{subsec:IllustrationMFE}

\begin{figure}[t] 
    \centering
    \small Mean-field equilibrium strategy $\widetilde{\vartheta}$\par
    \vspace{0.5ex}

    \begin{minipage}[t]{0.23\textwidth}
        \centering
        \includegraphics[width=\linewidth]{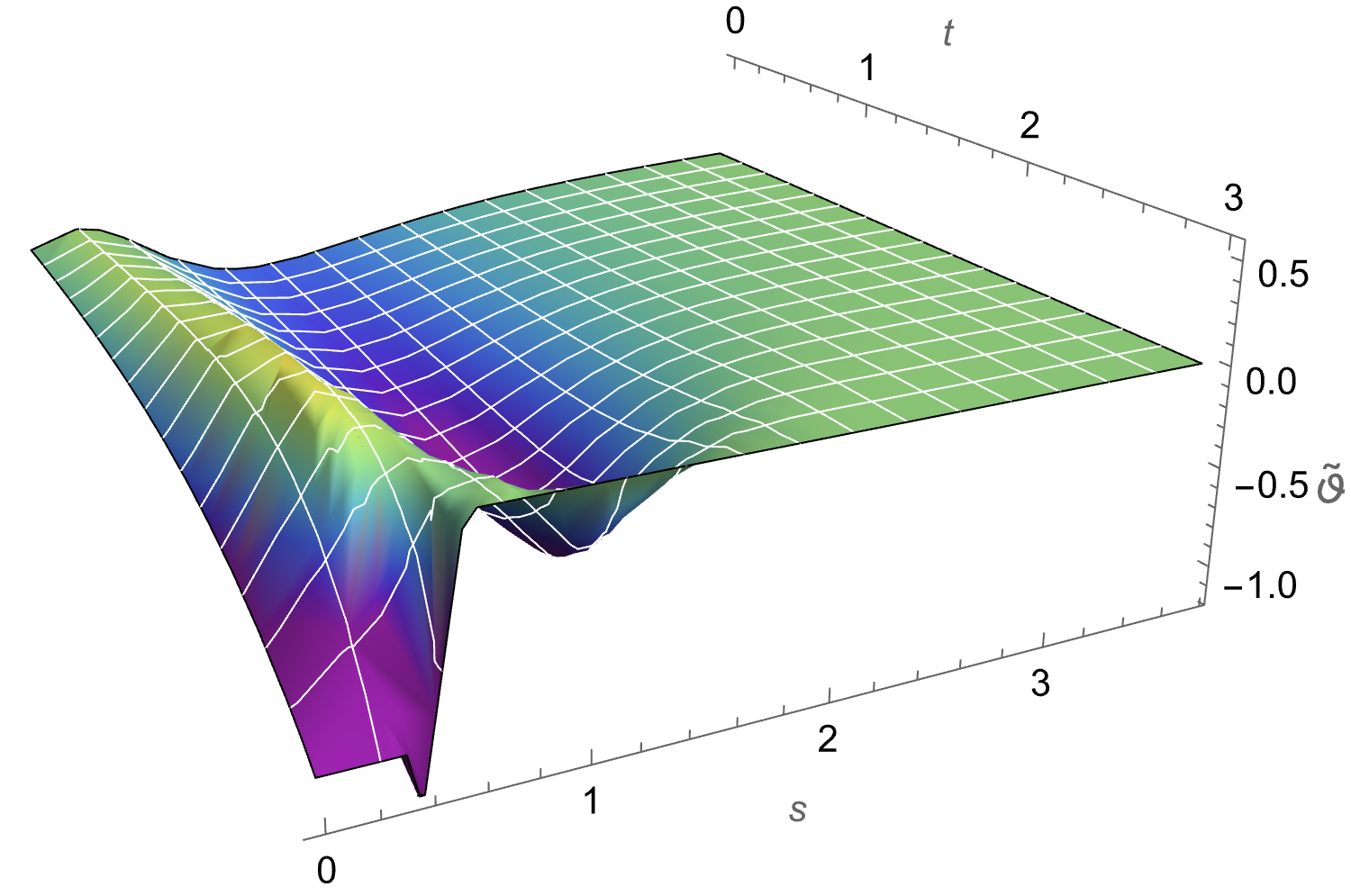}\par
        \vspace{0.3ex}
        \scriptsize $L_{t-}^B=(0,0)$
    \end{minipage}\hfill
    \begin{minipage}[t]{0.23\textwidth}
        \centering
        \includegraphics[width=\linewidth]{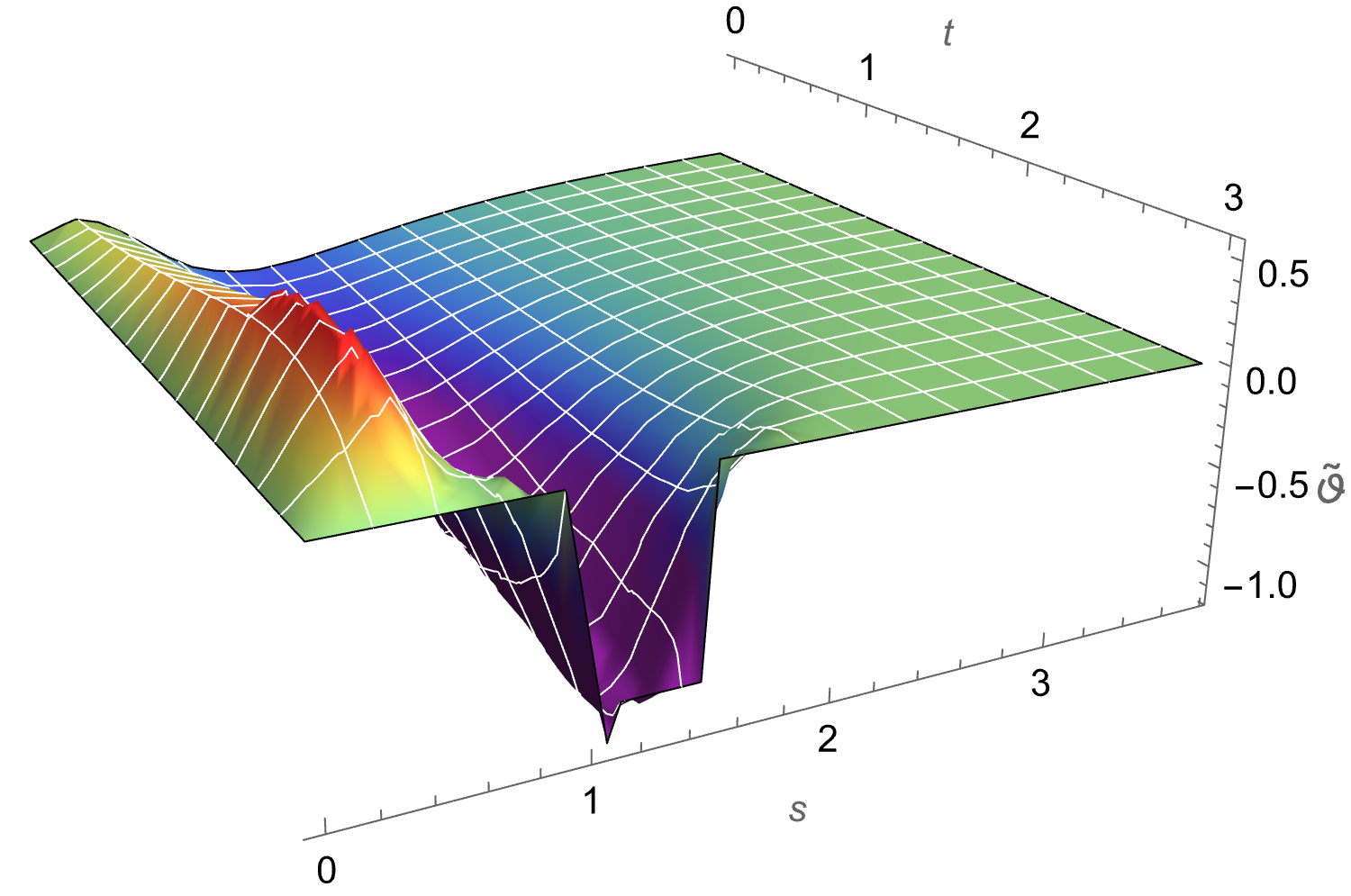}\par
        \vspace{0.3ex}
        \scriptsize $L_{t-}^B=(0,1)$
    \end{minipage}\hfill
    \begin{minipage}[t]{0.23\textwidth}
        \centering
        \includegraphics[width=\linewidth]{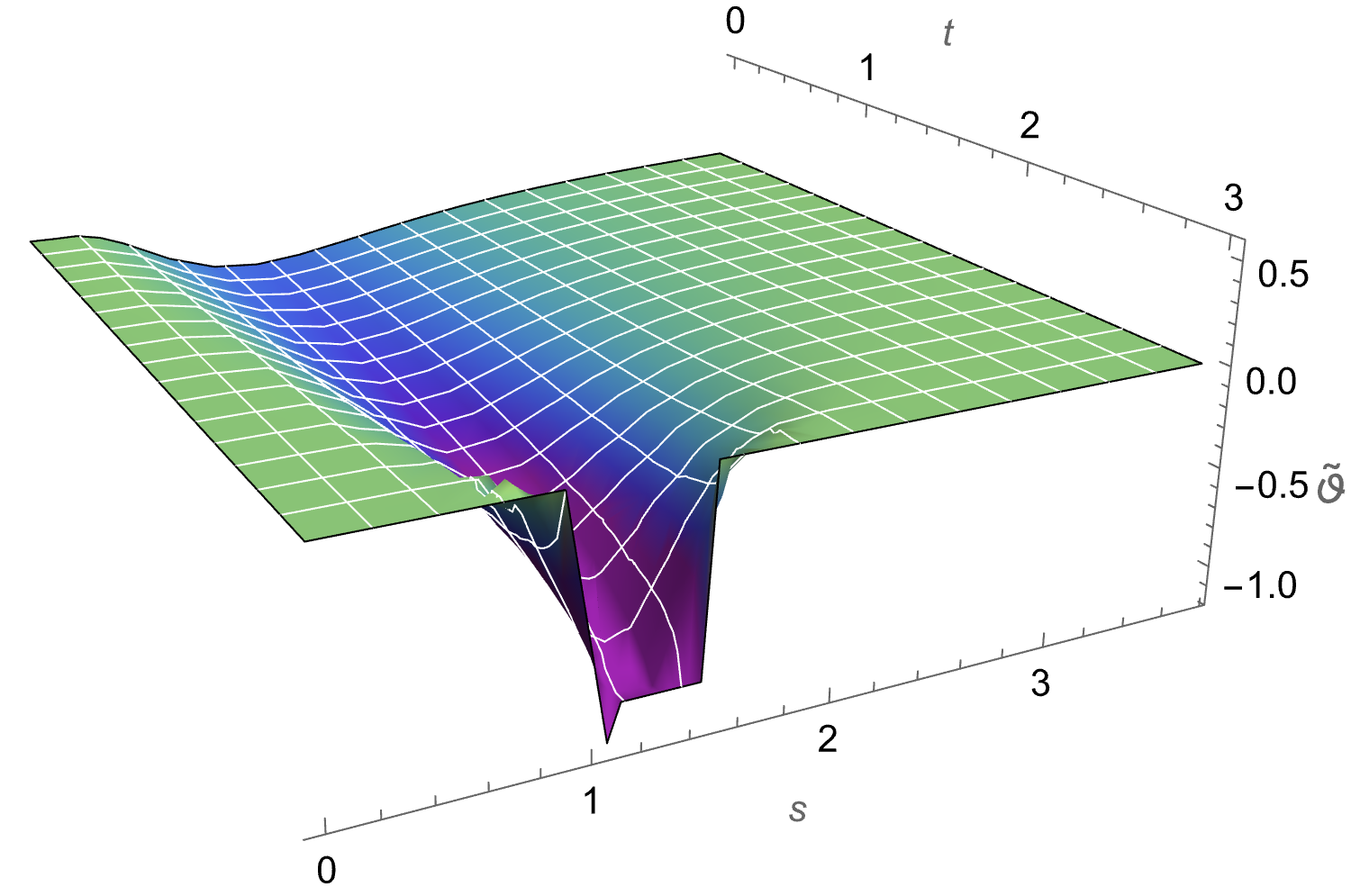}\par
        \vspace{0.3ex}
        \scriptsize $L_{t-}^B=(1,0)$
    \end{minipage}\hfill
    \begin{minipage}[t]{0.23\textwidth}
        \centering
        \includegraphics[width=\linewidth]{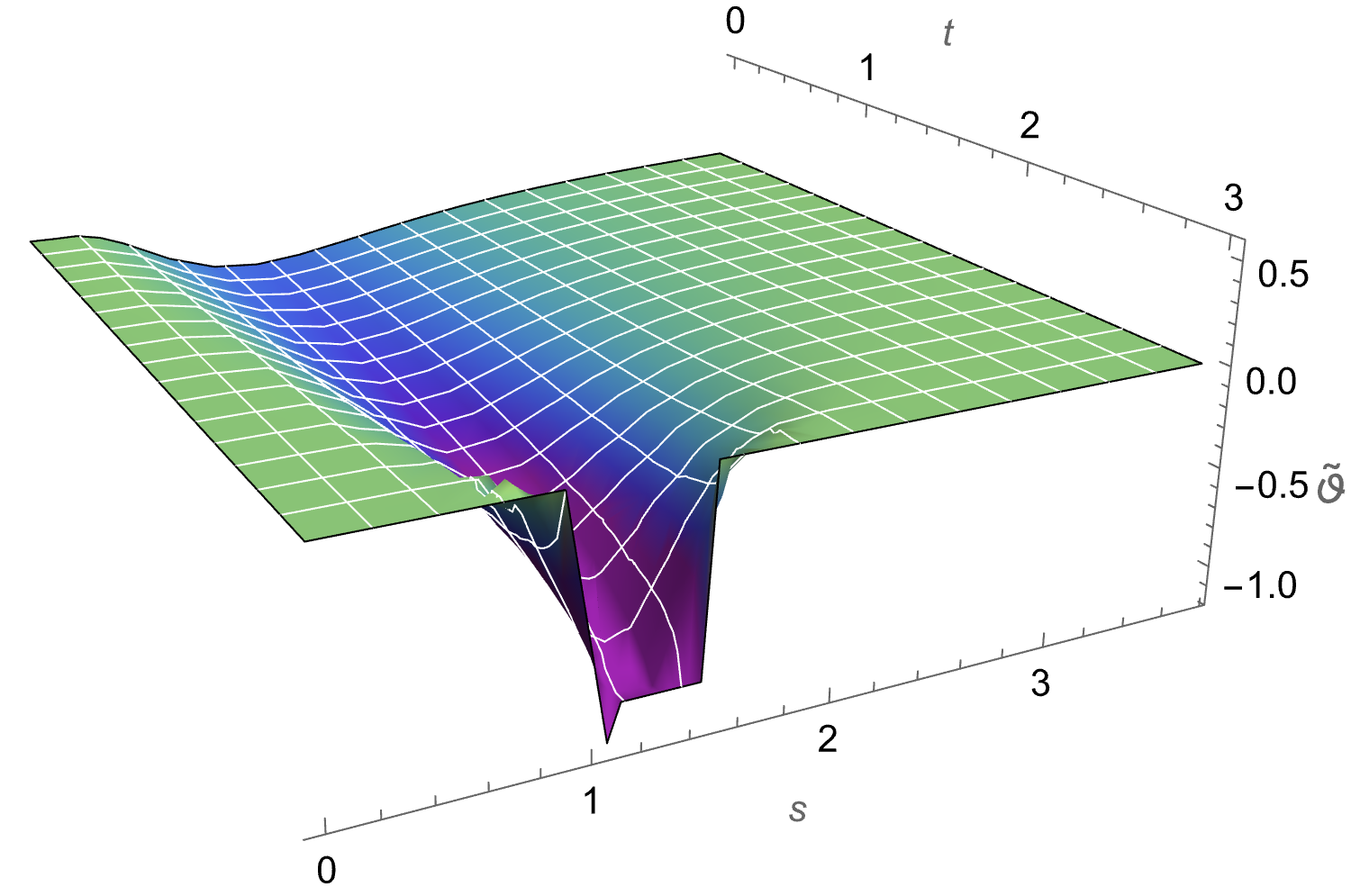}\par
        \vspace{0.3ex}
        \scriptsize $L_{t-}^B=(1,1)$
    \end{minipage}
    \caption{MFE strategy $\widetilde{\vartheta}$ (recalling relation \eqref{parametrizationVarthetaTheta}) for the MFG \eqref{MFG} with claim $B$ from \eqref{DefFinancialStopLossContractForPDE}. The MFE strategy $\widetilde{\vartheta}$ is given by a function $\widetilde{\vartheta}_t=\widetilde{\vartheta}(t,S_t,L_{t-})$ of time $t$, the asset price $S_t=s$, and the state $L_{t-}^B$ of the Markov chain defined in \eqref{DefLtB}. Parameters: $S_0=1,\sigma=0.3, \varphi = 0, C=1, K_1=0.5, K_2=0.5, \alpha=2.5, \rho=0.9, T=3, \lambda^0=0.9$ and $\lambda^1=1-\lambda^0$}\label{fig:MFE}
\end{figure}

\begin{figure}[t]
    \centering
    \small Decomposition of MFE strategy $\widetilde{\vartheta}$ 
    \par
    \begin{minipage}[t]{0.23\textwidth}
        \centering
        \includegraphics[width=\linewidth]{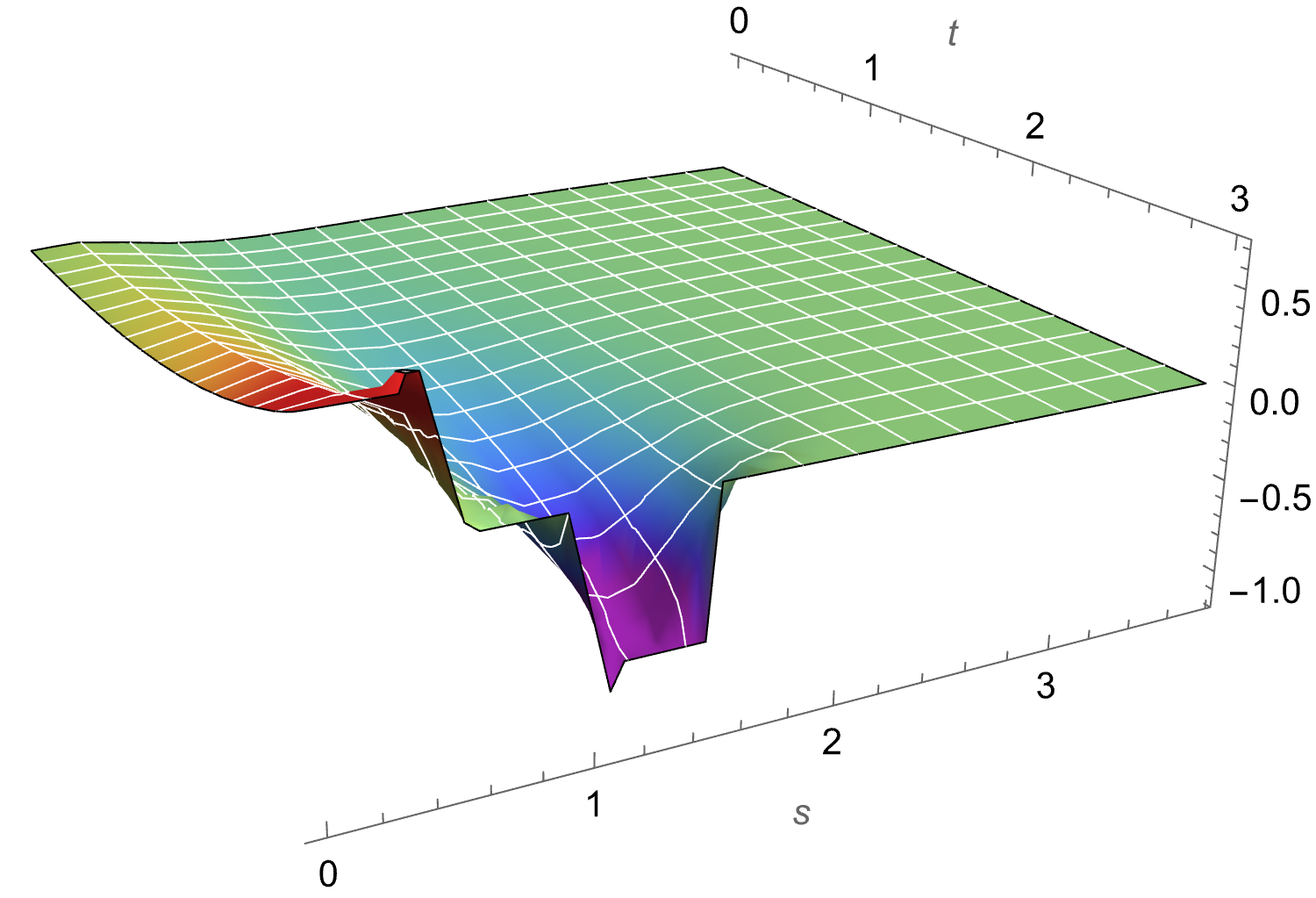}\par
        \vspace{0.3ex}
        \scriptsize hedging part $(\rho=0.9)$
    \end{minipage}\hfill
    \begin{minipage}[t]{0.23\textwidth}
        \centering
        \includegraphics[width=\linewidth]{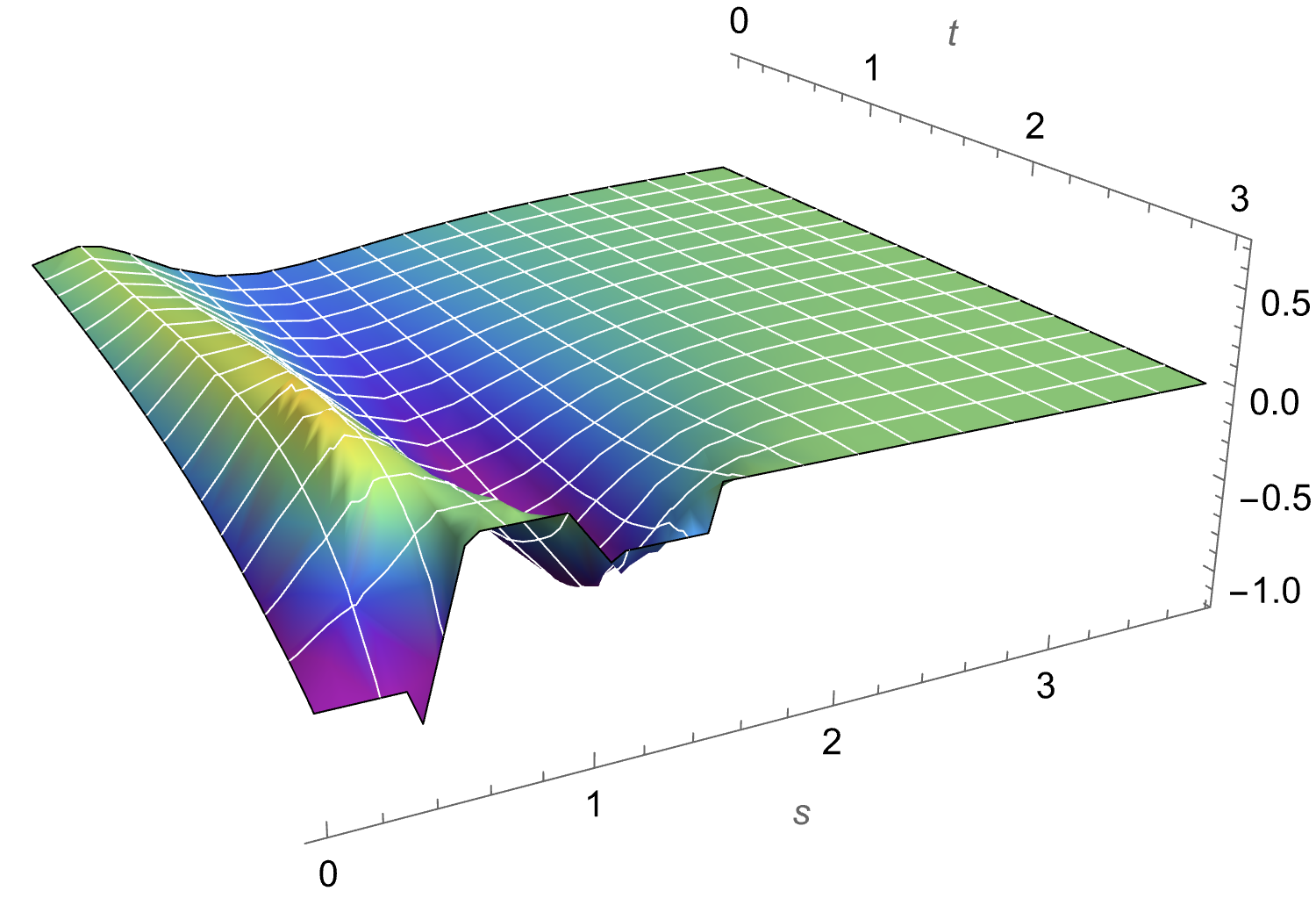}\par
        \vspace{0.3ex}
        \scriptsize mean-field part $(\rho=0.9)$
    \end{minipage}\hfill
     \begin{minipage}[t]{0.23\textwidth}
        \centering
        \includegraphics[width=\linewidth]{OderR09L09V01.png}\par
        \vspace{0.3ex}
        \scriptsize  MFE strategy $(\rho=0.9)$
    \end{minipage} \hfill
    \vspace{0.5ex}
     \begin{minipage}[t]{0.23\textwidth}
        \centering
        \includegraphics[width=\linewidth]{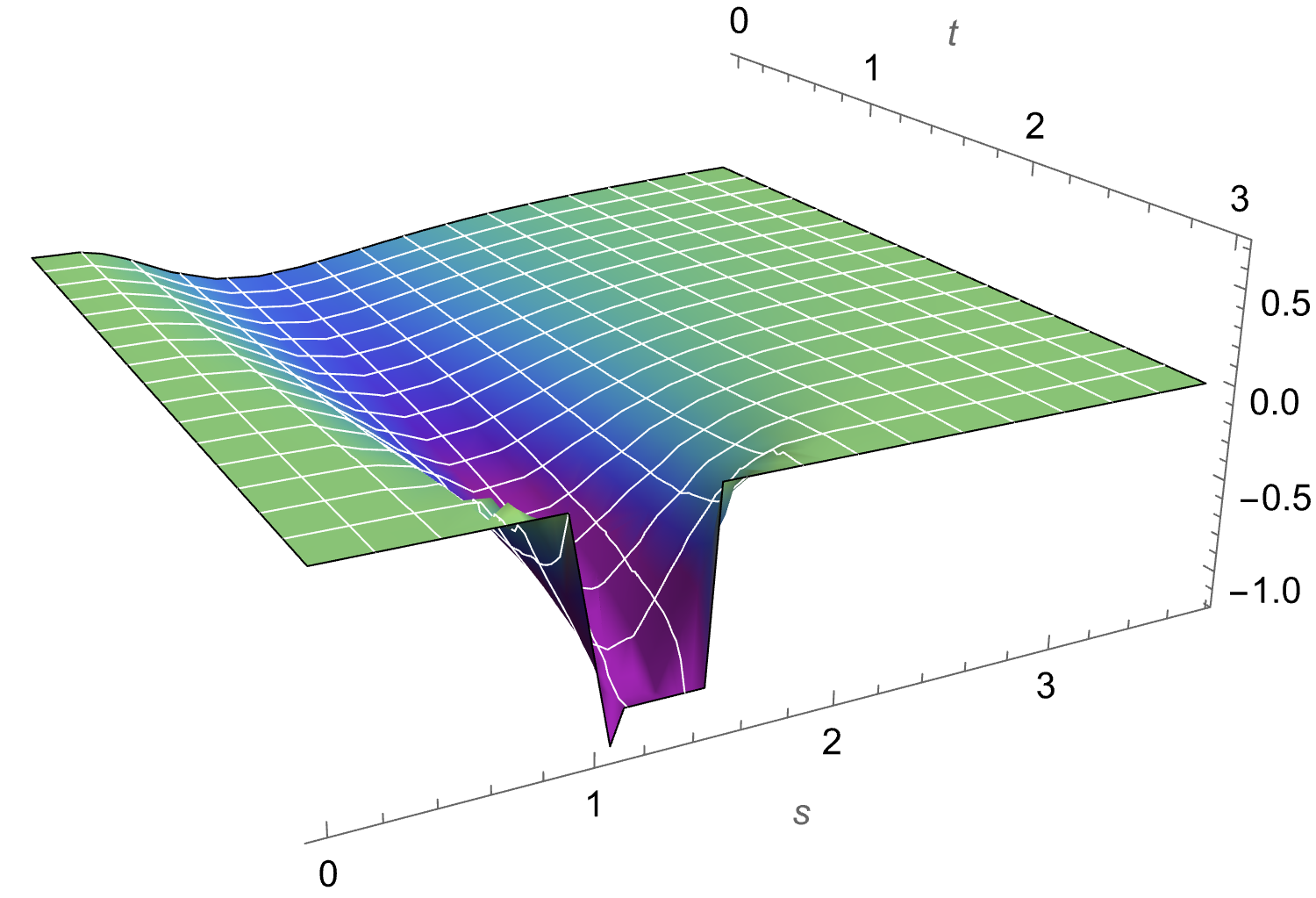}\par
        \vspace{0.3ex}
        \scriptsize hedging part=MFE strategy $(\rho=0)$ 
    \end{minipage}\hfill
    \caption{
        Decomposition of the MFE strategy $\widetilde{\vartheta}$ in  state $L_{t-}^B=(0,1)$ for $\rho=0.9$ and $\rho=0$. The remaining parameters are as in \hyperref[fig:MFE]{Fig. \ref*{fig:MFE}}.
    }\label{fig:DecompositionOfMFEstrat}
\end{figure}

\begin{figure}[t]
    \centering
    \small Impact of common noise intensity $\lambda^0$ on MFE strategy $\widetilde{\vartheta}$\par
    \vspace{0.5ex}

    \begin{minipage}[t]{0.23\textwidth}
        \centering
        \includegraphics[width=\linewidth]{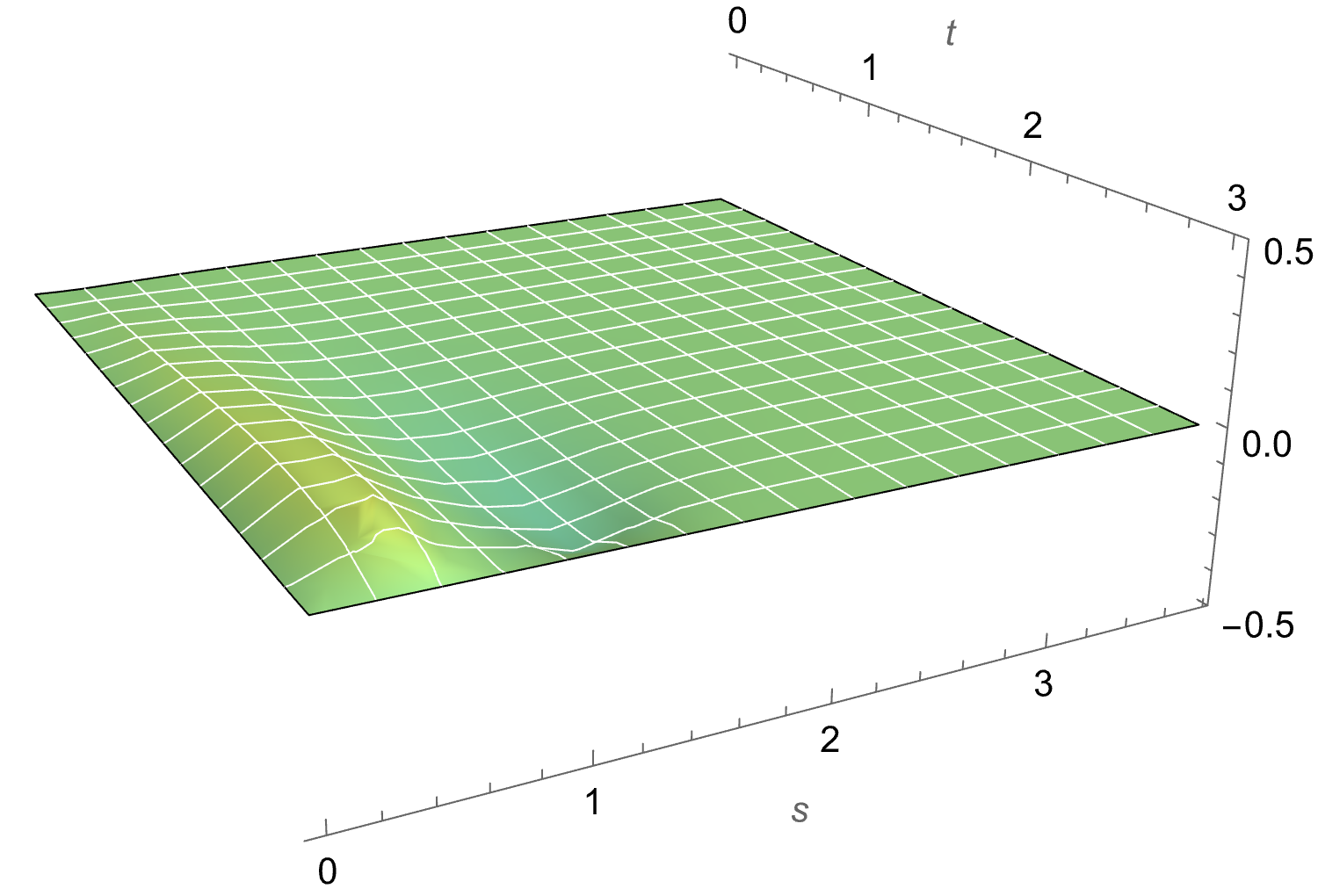}\par
        \vspace{0.3ex}
        \scriptsize $\widetilde{\vartheta}_{\lambda^0=0.3}-\widetilde{\vartheta}_{\lambda^0=0.1}$
    \end{minipage}\hfill
    \begin{minipage}[t]{0.23\textwidth}
        \centering
        \includegraphics[width=\linewidth]{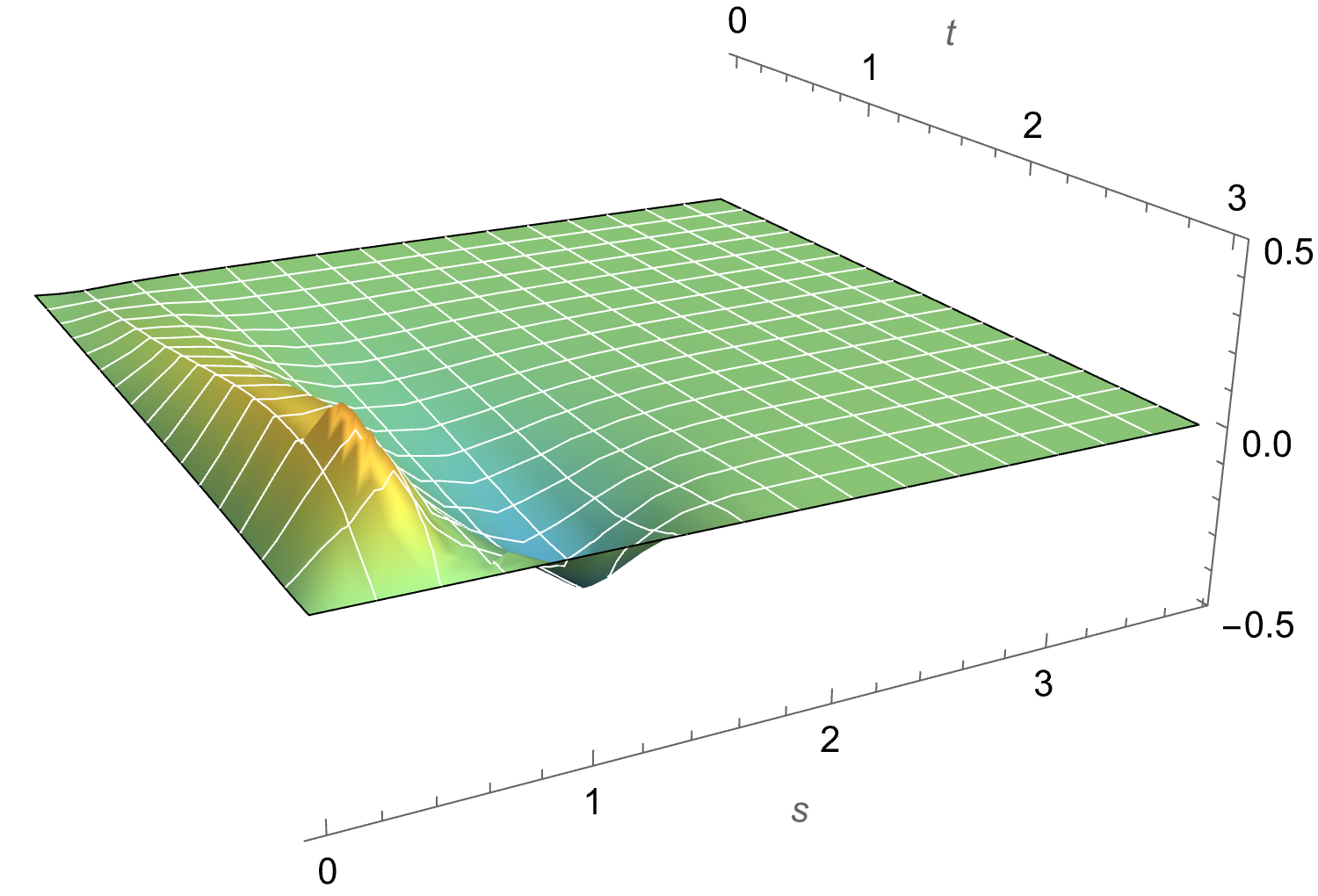}\par
        \vspace{0.3ex}
        \scriptsize $\widetilde{\vartheta}_{\lambda^0=0.5}-\widetilde{\vartheta}_{\lambda^0=0.1}$
    \end{minipage}\hfill
    \begin{minipage}[t]{0.23\textwidth}
        \centering
        \includegraphics[width=\linewidth]{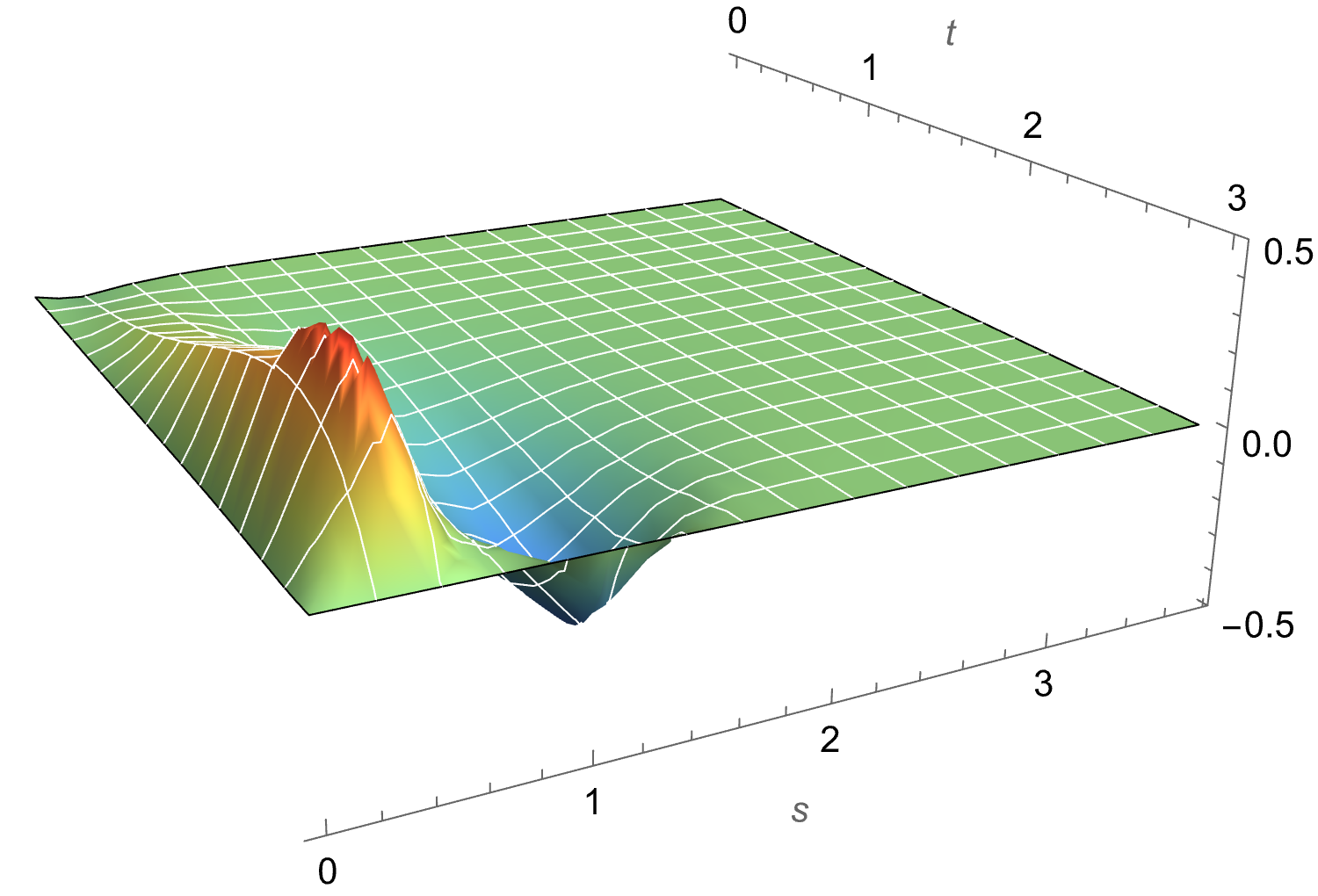}\par
        \vspace{0.3ex}
        \scriptsize $\widetilde{\vartheta}_{\lambda^0=0.7}-\widetilde{\vartheta}_{\lambda^0=0.1}$
    \end{minipage}\hfill
    \begin{minipage}[t]{0.23\textwidth}
        \centering
        \includegraphics[width=\linewidth]{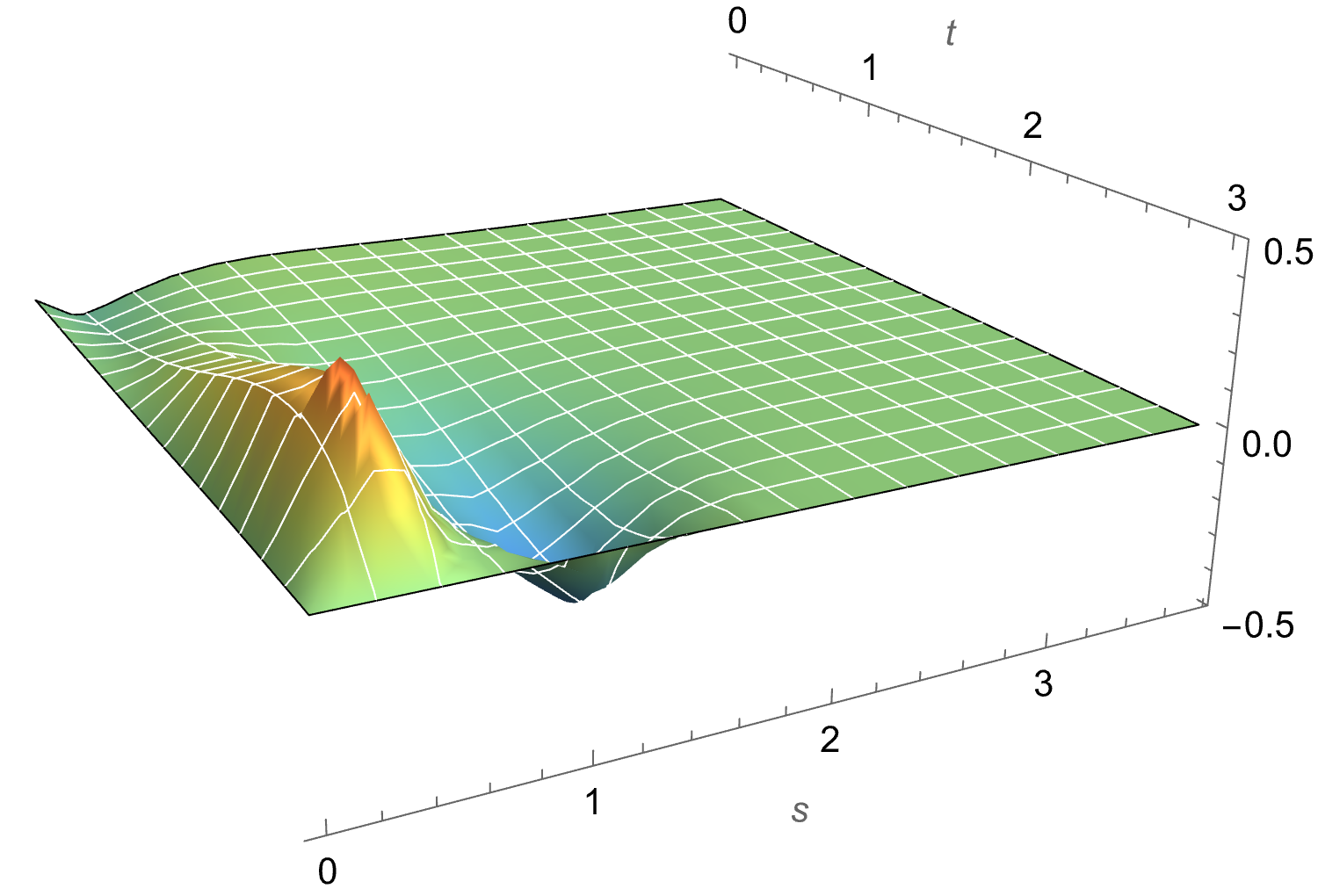}\par
        \vspace{0.3ex}
        \scriptsize $\widetilde{\vartheta}_{\lambda^0=0.9}-\widetilde{\vartheta}_{\lambda^0=0.1}$
    \end{minipage}
    \caption{
        Differences between the MFE strategies $\widetilde{\vartheta}$ for $\lambda^0\in \{0.3,0.5,0.7,0.9\}$ and the benchmark strategy for $\lambda^0=0.1$ in state $L_{t-}^B=(0,1)$. Here, $\lambda^1=1-\lambda^0$. The remaining parameters are as in \hyperref[fig:MFE]{Fig. \ref*{fig:MFE}}.
    }\label{fig:ImpactCommonNoise}
\end{figure}

Based on the characterizations \eqref{VarthetaPDEoptimalStratSingleAgentProb}-\eqref{PDEmfeGeneral} of the MFE $\widetilde{\vartheta}$ to the MFG \eqref{MFG} in the setting with a Markov modulated SDE for $S$, the present subsection illustrates in a concretely specified example the MFE (\Cref{fig:MFE}), its decomposition  (\Cref{fig:DecompositionOfMFEstrat}) as well as common noise effects on the MFE (\Cref{fig:ImpactCommonNoise}). 

For simplicity, let $S$ follow classical Black-Scholes price dynamics with the volatility and market price of risk being constants, that means $\sigma^i=\sigma$, $\varphi^i=\varphi$ do not vary with $i\in K$. We further take $\rho\neq 1$ and $\alpha>0$ as constants, and the claim $B$, against which the player wants to hedge and compare herself with the competition, as a financial stop-loss contract (cf. \Cref{remark:ExampleMoeller}). Specifically, here, we consider independent common and idiosyncratic Poisson processes $N^0,N^1$ with intensities $\lambda^0,\lambda^1$; the insurance loss is $\text{InsuranceLoss}_T=C\Indi_{\{N_T^0+N_T^1>0\}}$, i.e., of size $C>0$ if a common or idiosyncratic Poisson event has occurred up to time $T$, and the financial loss is $\text{FinancialLoss}_T=S_0-S_T$. Hence, the financial stop-loss contract $B$ in this section is given by
\begin{align}
    B:=(C \Indi_{\{N_T^0+N_T^1> 0\}}+ S_0-S_T-K_1)^+ \land K_2,\label{DefFinancialStopLossContractForPDE}
\end{align}
where $K_1,K_2\ge 0$ are constants. Let the Markov chain $L^B$ be described by
\begin{align}
    &L_t^B:=\begin{cases}
        (0,0), &\text{if } N_t^0 = N_t^1= 0, \\
        (0,1), &\text{if } N_t^0 = 0, N_t^1 > 0,\\
        (1,0), &\text{if } N_t^0 > 0, N_t^1 = 0,\\
        (1,1), &\text{if } N_t^0 > 0, N_t^1 > 0.\\
    \end{cases}\label{DefLtB}
\end{align}
The transition rates $q^{i\rightarrow j}$ for the Markov chain $L$, as formulated in the setting of \Cref{subsec:PDEcharMFE}, can now be specified using the intensities $\lambda^0, \lambda^1$.

We recall that the terminal conditions $h^i$ in the system of PDEs \eqref{SystemOfPDEs} have to satisfy $h^{L_T^B}(S_T)=B-\rho \mathbb{E}\left[B \vert \mathcal{F}_T^0\right]$. Thus, the functions $h^i$ are given by
\begin{align}
    &h^{(0,0)}(s)= (1-\rho e^{-\lambda^1 T}) \left(( S_0-s-K_1)^+ \land K_2\right) -\rho(1- e^{-\lambda^{1}T}) \left( (C + S_0-s-K_1)^+ \land K_2 \right),\notag\\
    &h^{(0,1)}(s)= (1-\rho(1- e^{-\lambda^1 T})) \left((C+S_0-s-K_1)^+ \land K_2\right) - \rho e^{-\lambda^1T} \left((S_0-s-K_1)^+ \land K_2\right),\notag\\
    &h^{(1,0)}(s)= h^{(1,1)}(s) = (1-\rho) \left((C+S_0-s-K_1)^+ \land K_2\right).\label{DefH}
\end{align}
Since $h$ is continuous, the corresponding system of PDEs \eqref{SystemOfPDEs} has a unique classical solution $v$ in $C([0,T]\times (0,\infty), \R^4)\cap C^{1,2}([0,T)\times (0,\infty), \R^4)$ by \cite[Thm.2.4]{becherer2005classical}.

Using the constructions \eqref{VarthetaPDEoptimalStratSingleAgentProb}-\eqref{PDEmfeGeneral}, we obtain the MFE $\widetilde{\vartheta}$ to the MFG \eqref{MFG} with claim $B$ \eqref{DefFinancialStopLossContractForPDE} in terms of the solutions to the system of PDEs \eqref{SystemOfPDEs} with $h$ defined in \eqref{DefH}. Indeed, we compute $\Pi((v_s^{L_{t-}^B}(t,S_t)\Indi_{[0,T)}(t))_{t\in [0,T]})$, which turns out to be 
\begin{align*}
    &\bigg(\bigg\{\bigg(\Indi_{\{N_{t-}^0=0\}}v_s^{(0,0)}(t,S_t)\bigg) e^{-\lambda^{1}\cdot t}+ \bigg(\Indi_{\{N_{t-}^0=0\}}v_s^{(0,1)}(t,S_t)\bigg)\left(1-e^{-\lambda^{1}\cdot t}\right)\\
    & \hspace{0.2cm} + \bigg(\Indi_{\{N_{t-}^0>0\}}v_s^{(1,0)}(t,S_t)\bigg) e^{-\lambda^{1}\cdot t} + \bigg(\Indi_{\{N_{t-}^0>0\}}v_s^{(1,1)}(t,S_t)\bigg) \left(1-e^{-\lambda^{1}\cdot t}\right)\bigg\}\Indi_{[0,T)}(t)\bigg)_{t\in [0,T]} 
\end{align*}
and thus, using equations \eqref{VarthetaPDEoptimalStratSingleAgentProb}-\eqref{PDEmfeGeneral}, the MFE strategy $\widetilde{\vartheta}$ (at time $t$) is equal to
\begin{align}
    &\frac{1}{1-\rho} \cdot \frac{\varphi}{\alpha\sigma S_t} + v_s^{L_{t-}}(t,S_t)\Indi_{[0,T)}(t) \notag\\
    & \ + \frac{\rho}{1-\rho} \cdot \bigg\{\bigg(\Indi_{\{N_{t-}^0=0\}}v_s^{(0,0)}(t,S_t)\bigg) e^{-\lambda^{1}\cdot t}+ \bigg(\Indi_{\{N_{t-}^0=0\}}v_s^{(0,1)}(t,S_t)\bigg)\left(1-e^{-\lambda^{1}\cdot t}\right)\notag\\
    & \ + \bigg(\Indi_{\{N_{t-}^0>0\}}v_s^{(1,0)}(t,S_t)\bigg) e^{-\lambda^{1}\cdot t} + \bigg(\Indi_{\{N_{t-}^0>0\}}v_s^{(1,1)}(t,S_t)\bigg) \left(1-e^{-\lambda^{1}\cdot t}\right)\bigg\}\Indi_{[0,T)}(t). \label{MFECalculatedExample}
\end{align} 
Finally, we note that the respective events of the Poisson processes can be related to the states of the Markov chain $L^B$ according to \Cref{DefLtB}.

Based on the representation \eqref{MFECalculatedExample} of the MFE, \Cref{fig:MFE} illustrates the MFE $\widetilde{\vartheta}$ to the MFG \eqref{MFG} with $B$ from \eqref{DefFinancialStopLossContractForPDE}, recalling the relation \eqref{parametrizationVarthetaTheta}. In the case where, e.g., up to time $t=2.5$ no common jump has occurred ($N_{2.5-}^0=0$) but an idiosyncratic jump has occurred ($N_{2.5-}^1>0$), we consider the second plot ($L_{2.5-}^B=(0,1)$) in \Cref{fig:MFE}. For the stock price being at $S_{2.5}=s=0.3$ at that time, \Cref{fig:MFE} shows, e.g., that the MFE strategy $\widetilde{\vartheta}$ holds approximately $\widetilde{\vartheta}_{2.5}\approx0.5$ units of the stock at $t=2.5$.

In the remainder of this section, we focus on trajectories with $N_{t-}^0=0$ and $N_{t-}^1>0$, i.e., $L_{t-}^B=(0,1)$. 

\Cref{fig:DecompositionOfMFEstrat} illustrates the decomposition
of the MFE strategy given in \eqref{charMFE} and financially interpreted in
\hyperref[remark:InterpretationDecompMFE]{\Cref*{remark:InterpretationDecompMFE}.1}. By the reparameterization \eqref{parametrizationVarthetaTheta}
of trading strategies, this decomposition takes the form $\widetilde{\vartheta}=\frac{1}{\sigma S}(\frac{1}{\alpha} \varphi + Z^B + Z)$ in the present one-dimensional setting. Since in the illustrated examples the market price of risk satisfies $\varphi=0$, the Merton component  $\frac{1}{\sigma S}\cdot \frac{1}{\alpha} \varphi$ is zero, so we omit the corresponding plot. 

For competition weight $\rho=0.9$, the first plot in \Cref{fig:DecompositionOfMFEstrat} shows the hedging component $\frac{1}{\sigma S}Z^B$ of the claim $B-\rho\mathbb{E}[B \vert \mathcal{F}_T^0]$. The second plot shows the mean-field component $\frac{1}{\sigma S}Z=\frac{1}{\sigma S} \cdot \frac{\rho}{1-\rho} \Pi(\theta^B)$, which arises from the mean-field interaction and can be interpreted, as explained in \Cref{remark:MKVfbsdeToMKVbsde}, as the optimal trading strategy associated with the endogenous terminal liability $\rho \mathbb{E}[X_T\vert \mathcal{F}_T^0]$. Here, the $\Pi$-projection averages the optimal single-agent strategies corresponding to the states of the Markov chain $L^B$ conditional on the common information. Since $N_{t-}^0=0$ in the present case, it averages  the strategies corresponding to $L_{t-}^B=(0,0)$ and $L_{t-}^B=(0,1)$. 
The third plot shows the resulting MFE strategy $\widetilde{\vartheta}$. We highlight that the competition weight $\rho$ affects not only the mean-field component, but also the hedging component through the claim $B-\rho\mathbb{E}[B \vert \mathcal{F}_T^0]$. Thus, for $\rho=0$, not only the mean-field component vanishes but also the hedging
component changes as the claim reduces to $B$. The fourth plot in \Cref{fig:DecompositionOfMFEstrat}
therefore shows both the hedging component and the MFE strategy, which coincide in this case.

We next study how common noise affects the MFE strategy by varying the intensity $\lambda^0$ of the common Poisson process $N^0$, while keeping $\lambda^0+\lambda^1=1$ and all other parameters fixed. \Cref{fig:ImpactCommonNoise} reveals complex and non-trivial MFE behavior, showing e.g., that the effect of $\lambda^0$ is non-monotonic, as seen from the difference plots $\widetilde{\vartheta}_{\lambda^0}-\widetilde{\vartheta}_{\lambda^0=0.1}$, where $\lambda^0=0.1$ denotes the  benchmark intensity. The top of the `hill' in the difference surface first increases as $\lambda^0$ rises from $0.3$ to $0.7$ and then decreases for $\lambda^0=0.9$. This is notable also by the changing color at the top of the `hill' (passing from green/yellow to orange to red and then back to orange).

We also briefly comment on competition weight effects on the MFE strategy, but we omit the corresponding plots. For the claim $B$ defined in  \Cref{DefFinancialStopLossContractForPDE}, we observed a monotone effect on the MFE  (for fixed $t$ and $s$). But for other payoffs, e.g., $B:=(C \Indi_{\{N_T^0> 0\}}\Indi_{\{N_T^1> 0\}} + S_0-S_T-K_1)^+ \land K_2$, we have also observed non-monotonicity. Thus, the choice of the contingent claim may affect whether the impact of the competition weight on the MFE strategy is monotonic or non-monotonic.

\section{Limiting mean-field game for 
vanishing risk aversion} \label{sec:AlphaTo0}

After considering in the previous \Cref{sec:PDE} the special case of a setting with a finite-state Markov chain, we now return to the general framework of \Cref{sec:Preliminaries,sec:FormulationMFG}. The present section shows, how for vanishing risk aversion $\alpha$ (i.e., as risk tolerance tends to infinity) one obtains convergence to a limiting MFG \eqref{MFGAlpha0} of a different type. This new mean-field portfolio game is about hedging under relative performance concerns with respect to a quadratic loss criterion, instead of exponential utility. 

Recall the probability measure $\widehat{\mathbb{P}}^0$ from \Cref{changeOfMEasure} with $\beta=0$, which plays a central role in this section, and does not depend on $\alpha$ (cf.\ \Cref{remark:measureHatP0}). 
Under $\widehat{\mathbb{P}}^0$, $\widehat{\boldsymbol{\nu}}^0$ denotes the compensator of the integer-valued random measure $\boldsymbol{\mu}$ (cf.\ \Cref{remark:PopertiesChangeOfMeasure}) and $\widehat{\tilde{\boldsymbol{\mu}}}^0:=\boldsymbol{\mu}-\widehat{\boldsymbol{\nu}}^0$ is the corresponding compensated random measure. The next theorem establishes convergence to the limiting MFG \eqref{MFGAlpha0} for vanishing risk aversion $\alpha$. 

\begin{theorem}\label{theorem:LimitMFG}
    Let $\mathbb{Q}\in \{\mathbb{P}, \widehat{\mathbb{P}}^0\}$ and let $\widetilde{\theta}^{B,\alpha}$ and $\widetilde{\theta}^{0,\alpha}$ be the mean-field equilibrium strategies to the mean-field game \eqref{MFG} with and without claim $B$. 
    Then we have for the difference $\widetilde{\psi}^{B,\alpha}:=\widetilde{\theta}^{B,\alpha}-\widetilde{\theta}^{0,\alpha}$ of mean-field equilibrium strategies
    \begin{align}
        \lim_{\alpha \rightarrow 0} \widetilde{\psi}^{B, \alpha} = \widetilde{\theta}^{L^2} \text{ in } \mathbb{H}_{\text{BMO}}^2(\mathbb{Q}), \label{convergenceAlphaTo0}
    \end{align}
    where $\widetilde{\theta}^{L^2}$ is the unique mean-field equilibrium strategy in $\mathbb{H}_{\text{BMO}}^2(\widehat{\mathbb{P}}^0)$ to the mean-field game
    \begin{equation}
	   \begin{cases}
            \text{1. Fix a real-valued random variable } F \text{ and}\\
            \text{2. }\text{find } \widetilde{\theta}\in \underset{\theta \in \mathbb{H}_{\text{BMO}}^2(\widehat{\mathbb{P}}^0)}{\argmin} \ \mathbb{E}^{\widehat{\mathbb{P}}^0}\left[\abs{X_T^\theta-B-\rho F}^2\right],\\
            \text{\phantom{2. }} \text{for wealth process $ X^{\theta}$ given by }\ddd X_t^{\theta}= \theta_t (\varphi_t \ddd t+ \ddd W_t), X_0^\theta=x_0.\\
            \text{3. Find a fixed point such that } F=\mathbb{E}\left[X_T^{\widetilde{\theta}}-B\middle \vert\mathcal{F}_T^{0}\right],\\
            \text{\phantom{3. }where $X_T^{\widetilde{\theta}}$ is the optimal wealth from step 2.}
        \end{cases}\label{MFGAlpha0}
    \end{equation}
    Moreover, there exists a constant $C$ such that
    \begin{equation*}
        \norm{\widetilde{\psi}^{B, \alpha} - \widetilde{\theta}^{L^2}}_{\mathbb{H}_{\text{BMO}}^2(\mathbb{Q})}\le \alpha C \quad\text{holds for all $\alpha\in (0,1]$.}
    \end{equation*}
\end{theorem}

Note that, at the single-agent level, it is well known that, as the risk-aversion parameter $\alpha$ tends to zero, the exponential utility-indifference hedging strategy converges to the risk-minimizing strategy, characterized under the relevant entropy-minimizing equivalent  martingale measure $\widehat{\mathbb{P}}^0$ by the Galtchouk-Kunita-Watanabe decomposition of the underlying claim; see \cite{becherer2006bounded,mania2005dynamic} and references therein. Under this measure $\widehat{\mathbb{P}}^0$, the limiting strategy also minimizes the expected squared terminal hedging error. \Cref{theorem:LimitMFG} provides an analogue statement at the mean-field level. It identifies the appropriate mean-field quantities ($\widetilde{\psi}^{B, \alpha}=\widetilde{\theta}^{B,\alpha}-\widetilde{\theta}^{0,\alpha}$) that converge to the MFE strategy of the MFG \eqref{MFGAlpha0} in which the representative agent minimizes under the measure $\widehat{\mathbb{P}}^0$ an expected squared terminal hedging error. We note that, if the market price of risk is zero, i.e.\ $\varphi=0$, the MFE strategy $\widetilde{\theta}^{B,\alpha}$ of the exponential MFG itself converges to the MFE strategy $\widetilde{\theta}^{L^2}$ of the MFG \eqref{MFGAlpha0}. Indeed, in this case the optimal strategy $\theta^{0,\alpha}$ for the single-agent problem \eqref{SingleAgentProbForB} without a claim, i.e., for $B=0$, is identically zero. The explicit construction \eqref{charMFE} of the MFE strategy from the corresponding single-agent strategy therefore yields $\widetilde{\theta}^{0,\alpha}=0$. The asserted convergence then follows directly from \Cref{theorem:LimitMFG}.

\begin{proof}[Proof of \Cref{theorem:LimitMFG}]
    The main idea for the proof is, instead of trying to establish convergence directly at the mean-field level, to leverage the explicit construction by \Cref{theorem:mainTheorem} for MFE strategies to MFGs of type \eqref{MFG} in terms of the optimal strategies for associated single-agent problems. Building on known convergence properties at the single-agent level, this construction allows us to lift the convergence results to the mean-field level. 

    An analogous construction by \Cref{proposition:L2HedgingMFG} for the MFE $\widetilde{\theta}^{L^2}$ to the limiting MFG \eqref{MFGAlpha0} in terms of a solution to a suitable corresponding single-agent problem then enables us to identify the limit as the MFE $\widetilde{\theta}^{L^2}$ of the limiting MFG.
 
    Let $C$ denote a generic constant that may change from line to line but does not depend on $\alpha$. 
    Recall $\widetilde{\psi}^{B, \alpha}$ being defined in the theorem. Let  $\theta^{B,\alpha}$ and $\theta^{0,\alpha}$ denote the optimal strategies to the single-agent optimization problems \eqref{SingleAgentProbForB} with and without a claim $B$. By \Cref{theorem:mainTheorem},  
    \begin{equation}
        \widetilde{\psi}^{B, \alpha}=\psi^{B, \alpha} + \frac{\rho}{1-\mathbb{E}[\rho]}\Pi\left(\psi^{B, \alpha}\right), \label{ReprHatThetaBAlpha}
    \end{equation}
    holds with $\psi^{B, \alpha}:=\theta^{B,\alpha}-\theta^{0,\alpha}$ being the difference of the
    optimal strategies with and without claim $B$ to the single-agent optimization problem \eqref{SingleAgentProbForB}.
    By \cite[Thm.4.6.]{becherer2006bounded}, there exists some constant such that for any difference of such optimal single-agent strategies the inequality
    \begin{align*}
        \sup_{t\in [0,T]}\norm{\mathbb{E}^{\widehat{\mathbb{P}}^0}\left[\int_t^T \abs{\psi_s^{B, \alpha} - Z_s^{B,L^2}}^2 \ddd s\bigg \vert \mathcal{F}_t\right]}_{L^\infty}
        = \norm{\psi^{B, \alpha} - Z^{B,L^2}}_{\mathbb{H}_{\text{BMO}}^2(\widehat{\mathbb{P}}^0)}^2\le \alpha^2 C
    \end{align*}
    holds true for all $\alpha\in (0,1]$, and for $(Y^{B,L^2}, Z^{B,L^2}, U^{B,L^2})\in \mathbb{S}^\infty(\widehat{\mathbb{P}}^0)\times \mathbb{H}_{\text{BMO}}^2(\widehat{\mathbb{P}}^0)\times \mathcal{L}_{\widehat{\boldsymbol{\nu}}^0}^2(\widehat{\mathbb{P}}^0)$ being the unique solution to the linear JBSDE (with trivial generator) 
    \begin{align}
        \left\{\arraycolsep=1pt \def\arraystretch{1.2}
		\begin{array}{ll}
    		\ddd Y_t^{B,L^2}&= Z_t^{B,L^2}\ddd \widehat{W}_t + \int_{E} U_t^{B,L^2}(e)\widehat{\tilde{\boldsymbol{\mu}}}^0(\ddd t, \ddd e),\\
			Y_T^{B,L^2}&=B-\rho \mathbb{E}\left[B \left\vert \mathcal{F}_T^0\right.\right].
        \end{array} \label{SingleBSDEForAlpha0} 
        \right. 
    \end{align}
    Equivalence of the norms $\tnorm{\cdot}_{\mathbb{H}_{\text{BMO}}^2 (\widehat{\mathbb{P}}^0)}$ and $\tnorm{\cdot}_{\mathbb{H}_{\text{BMO}}^2 (\mathbb{P})}$ (see \hyperref[remark:PopertiesChangeOfMeasure]{\Cref*{remark:PopertiesChangeOfMeasure}.3}) yields 
    \begin{equation} \label{singleAlpha0}
        \norm{\psi^{B, \alpha} - Z^{B,L^2}}_{\mathbb{H}_{\text{BMO}}^2(\mathbb{P})}^2\le \alpha^2 C 
        \quad   \text{for all $\alpha\in (0,1]$.}
    \end{equation}
     According to \Cref{proposition:L2HedgingMFG}, the unique MFE strategy $\widetilde{\theta}^{L^2}$ to the MFG \eqref{MFGAlpha0} is given by
    \begin{equation}
        \widetilde{\theta}^{L^2} = Z^{B,L^2}+\frac{\rho}{1-\mathbb{E}[\rho]}\Pi(Z^{B,L^2}). \label{RefReprTildeThetaL2}
    \end{equation}
    Since $\tnorm{\Pi(Z^{B,L^2})}_{\mathbb{H}_{\text{BMO}}^2(\mathbb{P})}\le \tnorm{Z^{B,L^2}}_{\mathbb{H}_{\text{BMO}}^2(\mathbb{P})}$ by \Cref{lemma:zBMOozBMO}, and since $\rho$ is bounded, $\mathcal{F}_0$-measureable (since $\mathcal{A}\subseteq \mathcal{F}_0$) and $\mathbb{E}[\rho]\neq 1$, we obtain from the representations \eqref{ReprHatThetaBAlpha}, \eqref{RefReprTildeThetaL2} together with the estimate  \eqref{singleAlpha0}, the following upper bound estimate (for all $\alpha\in (0,1]$)
    \begin{align*}
        \norm{\widetilde{\psi}^{B,\alpha}-\widetilde{\theta}^{L^2}}_{\mathbb{H}_\text{BMO}^2(\mathbb{P})}
        = \norm{\psi^{B, \alpha}- Z^{B,L^2} + \frac{\rho}{1-\mathbb{E}[\rho]}\Pi\left(\psi^{B, \alpha} - Z^{B,L^2}\right)}_{\mathbb{H}_\text{BMO}^2(\mathbb{P})} \le C\alpha. 
    \end{align*}
    Norm equivalence (\hyperref[remark:PopertiesChangeOfMeasure]{\Cref*{remark:PopertiesChangeOfMeasure}.3}), ensures the same convergence rate  in $\mathbb{H}_{\text{BMO}}^2(\widehat{\mathbb{P}}^0)$. In particular, for $\alpha\to 0 $, also the convergence holds as claimed in \Cref{convergenceAlphaTo0}.
\end{proof}

We have used the following proposition, that we postponed, in the preceding proof.

\begin{proposition}\label{proposition:L2HedgingMFG}
    There exists a unique mean-field equilibrium $\widetilde{\theta}^{L^2}\in \mathbb{H}_{\text{BMO}}^2(\widehat{\mathbb{P}}^0)$ to the mean-field game \eqref{MFGAlpha0}. It is given by
    \begin{align*}
        \widetilde{\theta}^{L^2}=Z^{L^2},
    \end{align*}
    where $(Y^{L^2},Z^{L^2},U^{L^2})\in \mathbb{S}^2(\widehat{\mathbb{P}}^0)\times \mathbb{H}_{\text{BMO}}^2(\widehat{\mathbb{P}}^0)\times \mathcal{L}_{\widehat{\boldsymbol{\nu}}^0}^2(\widehat{\mathbb{P}}^0)$ is the solution to the McKean-Vlasov JBSDE
    \begin{align}
        \left\{\arraycolsep=1pt \def\arraystretch{1.2}
		\begin{array}{ll}
    		\ddd Y_t&= Z_t\ddd \widehat{W}_t + \int_{E} U_t(e)\widehat{\tilde{\boldsymbol{\mu}}}^0(\ddd t, \ddd e),\\
			Y_T&=B+\rho\mathbb{E}\left[x_0+\int_0^T Z_s \ddd \widehat{W}_s - B \bigg \vert \mathcal{F}_T^0\right].\label{BSDEAlpha0MFG}
        \end{array}\right.
    \end{align}
    Moreover, we have $Z^{L^2}=Z^{B,L^2}+ \frac{\rho}{1-\mathbb{E}[\rho]}\Pi(Z^{B,L^2})$, where $(Y^{B,L^2}, Z^{B,L^2}, U^{B,L^2}) \in \mathbb{S}^2(\widehat{\mathbb{P}}^0)\times \mathbb{H}_{\text{BMO}}^2(\widehat{\mathbb{P}}^0)\times \mathcal{L}_{\widehat{\boldsymbol{\nu}}^0}^2(\widehat{\mathbb{P}}^0)$ is given by the (unique) solution to the JBSDE \eqref{SingleBSDEForAlpha0}.
\end{proposition}

The proof of this proposition is basically analogous to that of \Cref{theorem:mainTheorem} in the exponential case, even simpler in that 1.) the quadratic hedging problem is comparably basic, and 
2.) no additional change of measure is needed  since the claim is formulated  and we can thus work already under a martingale measure. The martingale representation property (\Cref{assumption:weakPPR} with \hyperref[example:CondIndepAndWPRP]{\Cref*{example:CondIndepAndWPRP}.3}) provides a convenient setup for quadratic hedging arguments. Below, we mention  auxiliary results and then explain how they combine to yield \Cref{proposition:L2HedgingMFG}. Except for \Cref{proposition:L2Hedging}, the proofs of these auxiliary results closely mirror their counterparts in the exponential-utility case; thus we omit the details for brevity. Note that \Cref{lemma:ExiUniqueSoluYZUBE} implies, in particular,  existence and uniqueness for the (very basic) JBSDE \eqref{SingleBSDEForAlpha0} as stated; the well-posedness of the (McKean-Vlasov) JBSDEs in \Cref{lemma:Alpha0SoluMKVbsde} and \Cref{lemma:Alpha0ExistenceUniquenessOfSingleBSDE} reduces to this result via similar arguments as in  the exponential-utility case. Finally, the statement of \Cref{proposition:L2Hedging} is classical. It can be shown by adapting the proof by \cite[Thm.1]{foellmer1986hedging}. 

\begin{proposition}\label{proposition:L2Hedging}
    For a given random variable $\xi\in L^2(\mathcal{F}_T,\widehat{\mathbb{P}}^0)$ we have that if the JBSDE
    \begin{align}
    \left\{\arraycolsep=1pt \def\arraystretch{1.2}
	\begin{array}{ll}
    	\ddd Y_t^{\xi}&= Z_t^{\xi} \ddd \widehat{W}_t + \int_{E} U_t^{\xi}(e)\widehat{\tilde{\boldsymbol{\mu}}}^0(\ddd t, \ddd e),\\
		Y_T^{\xi}&=\xi\label{BSDEPropL2Single}
    \end{array}\right.
    \end{align}
    has a solution $(Y^{\xi}, Z^{\xi}, U^{\xi})\in \mathbb{S}^2(\widehat{\mathbb{P}}^0)\times \mathbb{H}_{\text{BMO}}^2(\widehat{\mathbb{P}}^0)\times \mathcal{L}_{\widehat{\boldsymbol{\nu}}^0}^2(\widehat{\mathbb{P}}^0)$, then there exists a unique strategy $\theta^*$ for the optimization problem 
    \begin{align*}
        \text{minimize } \mathbb{E}^{\widehat{\mathbb{P}}^0}\left[\abs{X_T^\theta-\xi}^2\right] \text{ over } \theta\in \mathcal{L}_{T}^2(\widehat{\mathbb{P}}^0).
    \end{align*}
    It is given by $\theta^*=Z^\xi\in \mathbb{H}_{\text{BMO}}^2(\widehat{\mathbb{P}}^0)$.
\end{proposition}

\begin{proposition}\label{proposition:Alpha01zu1ZshgMgeMKVJBSDE}
    If for each mean-field equilibrium $\bar{\theta}^{L^2}$ of the MFG \eqref{MFGAlpha0} the JBSDE
    \begin{align}
    \left\{\arraycolsep=1pt \def\arraystretch{1.2}
		\begin{array}{ll}
    		\ddd Y_t^{\bar{\theta}^{L^2}}&= Z_t^{\bar{\theta}^{L^2}}\ddd \widehat{W}_t + \int_{E} U_t^{\bar{\theta}^{L^2}}(e)\widehat{\tilde{\boldsymbol{\mu}}}^0(\ddd t, \ddd e),\\
			Y_T^{\bar{\theta}^{L^2}}&= B-\rho \mathbb{E}\left[B \left\vert \mathcal{F}_T^0\right.\right] + \rho \mathbb{E}\left[x_0+\int_0^T 
            \bar{\theta}_s^{L^2} \ddd \widehat{W}_s \bigg\vert \mathcal{F}_T^0
            \right]
        \end{array}\label{Alpha0BSDEForGivenMFE}
        \right.
    \end{align}
    has a solution $(Y^{\bar{\theta}^{L^2}}, Z^{\bar{\theta}^{L^2}}, U^{\bar{\theta}^{L^2}})\in \mathbb{S}^2(\widehat{\mathbb{P}}^0)\times \mathbb{H}_{\text{BMO}}^2(\widehat{\mathbb{P}}^0)\times \mathcal{L}_{\widehat{\boldsymbol{\nu}}^0}^2(\widehat{\mathbb{P}}^0)$. Then there exists a solution $(Y^{L^2},Z^{L^2},U^{L^2})\in \mathbb{S}^2(\widehat{\mathbb{P}}^0)\times \mathbb{H}_{\text{BMO}}^2(\widehat{\mathbb{P}}^0)\times \mathcal{L}_{\widehat{\boldsymbol{\nu}}^0}^2(\widehat{\mathbb{P}}^0)$ to the McKean-Vlasov JBSDE \eqref{BSDEAlpha0MFG} if and only if the MFG \eqref{MFGAlpha0} has a mean-field equilibrium $\widetilde{\theta}^{L^2}$; It is given by $\widetilde{\theta}^{L^2}=Z^{L^2}$.
\end{proposition}

\begin{lemma}\label{lemma:Alpha0SoluMKVbsde}
    There exists a unique solution $(Y^{L^2},Z^{L^2},U^{L^2})\in \mathbb{S}^2(\widehat{\mathbb{P}}^0)\times \mathbb{H}_{\text{BMO}}^2(\widehat{\mathbb{P}}^0)\times \mathcal{L}_{\widehat{\boldsymbol{\nu}}^0}^2(\widehat{\mathbb{P}}^0)$ to the McKean-Vlasov JBSDE \eqref{BSDEAlpha0MFG}. Moreover, $Z^{L^2}=Z^{B,L^2}+ \frac{\rho}{1-\mathbb{E}[\rho]}\Pi(Z^{B,L^2})$, for $(Y^{B,L^2}, Z^{B,L^2}, U^{B,L^2}) $ in $ \mathbb{S}^\infty(\widehat{\mathbb{P}}^0)\times \mathbb{H}_{\text{BMO}}^2(\widehat{\mathbb{P}}^0)\times \mathcal{L}_{\widehat{\boldsymbol{\nu}}^0}^2(\widehat{\mathbb{P}}^0)$ with $U^{B,L^2}$ bounded being the solution to the JBSDE \eqref{SingleBSDEForAlpha0}.
\end{lemma}

\begin{lemma}\label{lemma:Alpha0ExistenceUniquenessOfSingleBSDE}
	For any strategy $\bar{\theta}\in \mathbb{H}_{\text{BMO}}^2(\widehat{\mathbb{P}}^0)$, the corresponding JBSDE \eqref{Alpha0BSDEForGivenMFE} has a unique solution $(Y^{\bar{\theta}^{L^2}}, Z^{\bar{\theta}^{L^2}}, U^{\bar{\theta}^{L^2}})\in \mathbb{S}^2(\widehat{\mathbb{P}}^0)\times \mathbb{H}_{\text{BMO}}^2(\widehat{\mathbb{P}}^0)\times \mathcal{L}_{\widehat{\boldsymbol{\nu}}^0}^2(\widehat{\mathbb{P}}^0)$.
\end{lemma}

\begin{proof}[Proof of \Cref{proposition:L2HedgingMFG}.]
    By \Cref{lemma:Alpha0SoluMKVbsde}, the McKean-Vlasov JBSDE \eqref{BSDEAlpha0MFG} has a unique solution $(Y^{L^2},Z^{L^2},U^{L^2})\in \mathbb{S}^2(\widehat{\mathbb{P}}^0)\times \mathbb{H}_{\text{BMO}}^2(\widehat{\mathbb{P}}^0)\times \mathcal{L}_{\widehat{\boldsymbol{\nu}}^0}^2(\widehat{\mathbb{P}}^0)$ and the $Z$ component is given by $Z^{L^2}=Z^{B,L^2}+\frac{\rho}{1-\mathbb{E}[\rho]}\Pi(Z^{B,L^2})$, where $(Y^{B,L^2}, Z^{B,L^2}, U^{B,L^2}) \in \mathbb{S}^2(\widehat{\mathbb{P}}^0)\times \mathbb{H}_{\text{BMO}}^2(\widehat{\mathbb{P}}^0)\times \mathcal{L}_{\widehat{\boldsymbol{\nu}}^0}^2(\widehat{\mathbb{P}}^0)$ is the solution to the JBSDE \eqref{SingleBSDEForAlpha0}. Furthermore, by \Cref{lemma:Alpha0ExistenceUniquenessOfSingleBSDE}, the assumptions of \Cref{proposition:Alpha01zu1ZshgMgeMKVJBSDE} are satisfied. Hence, by \Cref{proposition:Alpha01zu1ZshgMgeMKVJBSDE} it follows that a unique MFE $\widetilde{\theta}^{L^2}$ exists for the MFG \eqref{MFGAlpha0} and it is given by $\widetilde{\theta}^{L^2}=Z^{L^2}=Z^{B,L^2}+\frac{\rho}{1-\mathbb{E}[\rho]}\Pi(Z^{B,L^2})$.
\end{proof}


\bibliographystyle{alpha}
\bibliography{references.bib}


\appendix
\section{Appendix} \label{sec:appendix}

\begin{lemma}\label{lemma:solutionYB}
	For $\beta \in \{0,B\}$ with bounded $B\in L^\infty(\mathbb{P},\mathcal{F}_T)$, the JBSDE \eqref{trueCharBSDE} has a unique solution $(Y^\beta,Z^\beta,U^\beta)$ in $\mathbb{S}^\infty(\mathbb{P})\times \mathcal{L}_T^2(\mathbb{P})\times \mathcal{L}_{\boldsymbol{\nu}}^2(\mathbb{P})$ with $U^\beta$ being bounded (i.e., a bounded representative for $U^\beta\in \mathcal{L}_{\boldsymbol{\nu}}^2(\mathbb{P})$ can be chosen). Moreover, $Z^\beta$ is in $\mathbb{H}_{\text{BMO}}^2(\mathbb{P})\subseteq \mathcal{L}_T^2(\mathbb{P})$. 
\end{lemma}

\begin{proof}
	Let $\beta \in \{0,B\}$ with $B\in L^\infty(\mathbb{P},\mathcal{F}_T)$ bounded. Since $\beta$ and $\rho$ are bounded, the terminal condition of the JBSDE \eqref{trueCharBSDE} is bounded. Since the market price of risk $\varphi$ is bounded, $z\mapsto z\varphi_t$ is Lipschitz and since $\alpha$ is also bounded away from 0, $\tabs{\varphi}^2/(2\alpha)$ is bounded. Further, $u\mapsto (\exp(\alpha u)-1-\alpha u)/\alpha$ is absolutely continuous and locally bounded from above, since $\alpha$ is bounded. Thus, by \cite[Prop.4.3]{becherer2019monotone} the JBSDE \eqref{trueCharBSDE} has a unique solution $(Y^\beta,Z^\beta, U^\beta)\in \mathbb{S}^\infty(\mathbb{P})\times \mathcal{L}_T^2(\mathbb{P})\times \mathcal{L}_{\boldsymbol{\nu}}^2(\mathbb{P})$ (with $U^\beta$ bounded according to \cite[Lem.2.2]{becherer2019monotone}). Using the same arguments as in the proof of Theorem 4.1 in \cite{becherer2006bounded}, it follows that $Z^\beta\in \mathbb{H}_{\text{BMO}}^2(\mathbb{P})$.
\end{proof}

\begin{lemma}\label{lemma:La}
	For any process $\eta\in \mathcal{L}_T^2(\mathbb{P})$, we have
	\begin{align*}
		\mathbb{E}\left[\int_0^t \eta _s \cdot \ddd W_s \vert \mathcal{F}_t^0\right]= \int_0^t \Pi(\eta)_s\cdot \ddd W_s \quad \text{and} \quad 
		\mathbb{E}\left[\int_0^t \eta _s \ddd s \vert \mathcal{F}_t^0\right]= \int_0^t \Pi(\eta)_s \ddd s 
	\end{align*}
	with the notation $\Pi(\eta)$ as defined in \Cref{rn:defPi}. Cf.\ \cite[Lem.B.1]{lacker2022superposition}.
\end{lemma}
 
\begin{proof}
	To prove the first statement, for $\eta$, $W$ being multi-dimensional, it suffices to show the claim coordinate-by-coordinate. The coordinate-wise claim for bounded $\eta$ is obtained by the monotone class theorem. The claim extends to general $\eta\in \mathcal{L}_T^2(\mathbb{P})$ by approximating with bounded $\eta^m:=\Indi_{\{\tabs{\eta}\le m\}}\eta$. The second equality claimed is proven analogously.
\end{proof}

\begin{lemma}\label{lemma:zBMOozBMO}
	For $z\in \mathbb{H}_{\text{BMO}}^2(\mathbb{P})$ we have $\norm{\Pi(z)}_{\mathbb{H}_{\text{BMO}}^2(\mathbb{P})}\le \norm{z}_{\mathbb{H}_{\text{BMO}}^2(\mathbb{P})}$. In particular, $\Pi(z)$ is in $\mathbb{H}_{\text{BMO}}^2(\mathbb{P})$.
\end{lemma}

\begin{proof}
	We have $\mathbb{E}\big[\int_t^T \tabs{\Pi(z)}_s^2 \ddd s \big\vert \mathcal{F}_t\big]\le \mathbb{E}\big[\int_t^T \Pi(\tabs{z}^2)_s \ddd s \big\vert \mathcal{F}_t\big]$ for $z\in \mathbb{H}_{\text{BMO}}^2(\mathbb{P})$ and $t\le T$, by Jensen's inequality. Now, it suffices to show
	\begin{align}
		\mathbb{E}\left[\int_t^T \Pi(\tabs{z}^2)_s \ddd s \bigg\vert \mathcal{F}_t\right]\le \norm{z}_{\mathbb{H}_{\text{BMO}}^2(\mathbb{P})}^2.\label{IneqaulityPiProof}
	\end{align}	
	Since $\int_t^T \Pi(\tabs{z}^2)_s \ddd s$ is $\mathcal{F}_T^0$-measurable and $\mathcal{F}_t$ and $\mathcal{F}_T^{0}$ are conditionally independent given $\mathcal{F}_t^{0}$, we have according to \cite[Sect.3.2 Prop.13.(iv)]{rao2006probability}
	\begin{align}
		\mathbb{E}\left[\int_t^T \Pi(\tabs{z}^2)_s \ddd s \bigg\vert \mathcal{F}_t\right]
		=\mathbb{E}\left[\int_t^T \Pi(\tabs{z}^2)_s \ddd s \bigg\vert \mathcal{F}_t^0\right].\label{ProofPiBMOCondIndep}
	\end{align}
	Since for any $A_t\in \mathcal{F}_t^0$ the set $A_t\times (t,T]$ is in $ \mathcal{P}(\mathbb{F}^0)$ and $\Pi$ is defined as the conditional expectation on $\mathcal{P}(\mathbb{F}^0)$ under the finite measure $\mathbb{P}\otimes \ddd t$, we obtain the equalities
	\begin{align*}
		\mathbb{E}\left[\int_t^T \Pi(\tabs{z}^2)_s \ddd s \bigg\vert \mathcal{F}_t^0\right]
		=\mathbb{E}\left[\int_t^T \tabs{z}_s^2 \ddd s \bigg\vert \mathcal{F}_t^0\right]
		=\mathbb{E}\left[\mathbb{E}\left[\int_t^T \tabs{z}^2_s \ddd s \bigg\vert \mathcal{F}_t\right]\bigg\vert \mathcal{F}_t^0\right],
	\end{align*} 
	with the right side being dominated by $\norm{z}_{\mathbb{H}_{BMO}^2(\mathbb{P})}^2$. By \eqref{ProofPiBMOCondIndep} this yields \eqref{IneqaulityPiProof}.
\end{proof}

\begin{example}\label{example:IntInL2}
    $\mathbb{E}[\int_0^T Z_s \ddd \widehat{W}_s \vert \mathcal{F}_T^0]$ is in $L^2(\mathcal{F}_T,\widehat{\mathbb{P}})$ for $Z\in\mathbb{H}_{\text{BMO}}^2(\mathbb{P})$. This follows from \Cref{lemma:La}, \Cref{lemma:zBMOozBMO} and \hyperref[remark:PopertiesChangeOfMeasure]{\Cref*{remark:PopertiesChangeOfMeasure}.3}.
\end{example}

\begin{lemma}\label{lemma:ExiUniqueSoluYZUBE}
    There exists a solution $(Y^{B,L^2}, Z^{B,L^2}, U^{B,L^2})\in \mathbb{S}^\infty(\widehat{\mathbb{P}}^0)\times \mathbb{H}_{\text{BMO}}^2(\widehat{\mathbb{P}}^0)\times \mathcal{L}_{\widehat{\boldsymbol{\nu}}^0}^2(\widehat{\mathbb{P}}^0)$ with $U^{B,L^2}$ bounded to the JBSDE \eqref{SingleBSDEForAlpha0}, being unique in $\mathbb{S}^2(\widehat{\mathbb{P}}^0)\times \mathcal{L}_T^2(\widehat{\mathbb{P}}^0)\times \mathcal{L}_{\widehat{\boldsymbol{\nu}}^0}^2(\widehat{\mathbb{P}}^0)$.
\end{lemma}

\begin{proof}
    We start with the existence. 
    Since $B$ and $\rho$ are bounded, we obtain according to \cite[Prop.4.3]{becherer2019monotone}, that there exists a solution $(Y^{B,L^2}, Z^{B,L^2}, U^{B,L^2})\in \mathbb{S}^\infty(\widehat{\mathbb{P}}^0)\times \mathcal{L}_T^2(\widehat{\mathbb{P}}^0)\times \mathcal{L}_{\widehat{\boldsymbol{\nu}}^0}^2(\widehat{\mathbb{P}}^0)$ to the JBSDE \eqref{SingleBSDEForAlpha0} and $Y^{B,L^2}$ is bounded by $(1+\tabs{\rho}_\infty)\tabs{B}_\infty$. With \cite[Lem.2.2]{becherer2019monotone} we obtain by the boundedness of $Y$ that $U$ can be chosen as a bounded representative of $\mathcal{L}_{\widehat{\boldsymbol{\nu}}^0}^2(\widehat{\mathbb{P}}^0)$ and since we have no generator, we obtain by \cite[Lem.2.3]{becherer2019monotone} that $Z$ is in $\mathbb{H}_{\text{BMO}}^2(\widehat{\mathbb{P}}^0)$. Thus, we have a solution $(Y^{B,L^2}, Z^{B,L^2}, U^{B,L^2})\in \mathbb{S}^\infty(\widehat{\mathbb{P}}^0)\times \mathbb{H}_{\text{BMO}}^2(\widehat{\mathbb{P}}^0)\times \mathcal{L}_{\widehat{\boldsymbol{\nu}}^0}^2(\widehat{\mathbb{P}}^0)$ with $U^{B,L^2}$ bounded to the JBSDE \eqref{SingleBSDEForAlpha0}. 
    The uniqueness of the solution $(Y^{B,L^2}, Z^{B,L^2}, U^{B,L^2})$ in $\mathbb{S}^2(\widehat{\mathbb{P}}^0)\times \mathcal{L}_T^2(\widehat{\mathbb{P}}^0)\times \mathcal{L}_{\widehat{\boldsymbol{\nu}}^0}^2(\widehat{\mathbb{P}}^0)$ follows by \cite[Prop.3.2]{becherer2006bounded} because the terminal condition $Y_T^{B,L^2}$ of the JBSDE \eqref{SingleBSDEForAlpha0} is in $L^2(\widehat{\mathbb{P}}^0)$ and the JBSDE \eqref{SingleBSDEForAlpha0} has a generator, which is constant $0$ and thus Lipschitz continuous.
\end{proof}

\end{document}